\documentclass[12pt,a4paper,reqno]{amsart}
\usepackage{amsmath}
\usepackage{physics}
\usepackage{amssymb}
\numberwithin{equation}{section}
\usepackage{graphicx}
\usepackage{amssymb}
\usepackage{amsthm}
\usepackage{enumitem}
\usepackage{tikz}
\usetikzlibrary{shapes}
\usepackage{mathrsfs}
\usepackage{hyperref}
\usetikzlibrary{arrows}
\usepackage[utf8]{inputenc}
\usetikzlibrary{calc}
\numberwithin{equation}{section}
\usepackage{tikz-cd}
\usepackage[utf8]{inputenc}
\usepackage[T1]{fontenc}

\usepackage{tikz-cd}
\usepackage{txfonts}
\usepackage{lipsum} 
\let\oldpart\part
\renewcommand\part{\newpage\oldpart}

\usepackage[%
shortalphabetic, 
abbrev, 
backrefs, 
]{amsrefs}
\hypersetup{%
	linktoc=all, 
	colorlinks=true,
	citecolor=magenta,
}

\NewDocumentCommand\Crefnameitem { m m m O{\textup} O{(\roman*)}} {%
	\Crefname{#1enumi}{#2}{#3} 
	\AtBeginEnvironment{#1}{%
		\crefalias{enumi}{#1enumi}%
		\setlist[enumerate,1]{
			label={#4{#5}.},
			ref={#5}
		}%
	}  
}

\newtheorem{thm}{Theorem}[section]
\newtheorem*{theorem*}{Theorem}

\newtheorem{prop}[thm]{Proposition}
\newtheorem{lm}[thm]{Lemma}

\newtheorem{coro}[thm]{Corollary}

\newcommand{\nc}{\newcommand}

\nc{\delete}[1]{{}}

	\nc{\mlabel}[1]{\label{#1}}  
	\nc{\mcite}[1]{\cite{#1}}  
	\nc{\mref}[1]{\ref{#1}}  
	\nc{\meqref}[1]{\eqref{#1}}  
	\nc{\mbibitem}[1]{\bibitem{#1}} 

\delete{
	\nc{\mlabel}[1]{\label{#1} {{\small\tt{{\ }\ (#1)}}}}                
	\nc{\mcite}[1]{\cite{#1}{{\small\tt{{\ }(#1)}}}}  
	\nc{\mref}[1]{\ref{#1}{{\small\tt{{\ }(#1)}}}}  
	\nc{\meqref}[1]{\eqref{#1}{{\small\tt{{\ }(#1)}}}}  
	\nc{\mbibitem}[1]{\bibitem[\bf #1]{#1}} 
}

\newcommand\numberthis{\addtocounter{equation}{1}\tag{\theequation}}

\newtheorem{innercustomgeneric}{\customgenericname}
\providecommand{\customgenericname}{}
\newcommand{\newcustomtheorem}[2]{%
	\newenvironment{#1}[1]
	{%
		\renewcommand\customgenericname{#2}%
		\renewcommand\theinnercustomgeneric{##1}%
		\innercustomgeneric
	}
	{\endinnercustomgeneric}
}
\newcustomtheorem{customthm}{Theorem}

\theoremstyle{definition}
\newtheorem{example}[thm]{Example}
\newtheorem{df}[thm]{Definition}
\newtheorem{remark}[thm]{Remark}

\allowdisplaybreaks

\newcommand{\N}{\mathbb{N}}
\newcommand{\Z}{\mathbb{Z}}
\newcommand{\Q}{\mathbb{Q}}

\def \ra {\rightarrow}

\def \C {\mathbb{C}}

\def\la{\lambda}

\def \al{\alpha}

\def \om{\omega}
\def\hr{\hookrightarrow}
\def \ga {\gamma}
\def \ka{\kappa}

\def \b {\beta}
\def \op {\oplus}

\def \ssq{\subseteq}
\def \vac {\mathbf{1}}

\def \g {\mathfrak{g}}
\def \h {\mathfrak{h}}

\def \Hom {\mathrm{Hom}}

\def \End {\mathrm{End}}

\def\Id{\mathrm{Id}}
\def\Ind{\mathrm{Ind}}
\def\Res{\mathrm{Res}}
\def \wt {\mathrm{wt}}
\def\ep{\epsilon}

\def \bs {\backslash}

\def\L{\mathcal{L}}

\def\ds{\dots}
\def\o{\otimes}

\def\spn{\mathrm{span}}
\def\Om{\Omega}
\def\scrU{\mathscr{U}}
\def\fp{\mathfrak{p}}

\def\<{\langle}
\def\>{\rangle}
\def\Ind{\mathrm{Ind}}

\def\rN{\mathrm{N}}
\def\L{\mathsf{L}}
\def\U{\mathsf{U}}
\def\Mod{\mathsf{Mod}}

\def\fA{\mathfrak{A}}

\def\rN{\mathrm{N}}
\def\L{\mathsf{L}}
\def\sR{\mathsf{R}}
\def\U{\mathsf{U}}
\def\A{\mathsf{A}}

\def\adm{\mathsf{Adm}}
\def\sl{\mathfrak{sl}}

\def\ad{\mathrm{ad}}
\def\fin{\mathrm{fin}}
\def\LL{\mathsf{L}}
\def\Ord{\mathrm{Ord}}

\begin {document}

\title{Finite induction functor for vertex operator algebras}
\author{Jianqi Liu}
\address{Department of Mathematics, University of Pennsylvania, Philadelphia, PA, 19104}
\email{jliu230@sas.upenn.edu}
\maketitle

\begin{abstract}
In this paper, we introduce a new induction functor $\Ind^V_U$ between module categories corresponding to an embedding of vertex operator algebras (VOAs) $U \hookrightarrow V$. This induction functor is essentially defined at the level of the finite (Zhu) algebras, which we call the \emph{finite induction functor}. Under suitable conditions on $U$ and $V$, we prove that this functor satisfies the usual properties of induction functors, such as Frobenius reciprocity, functorial property for compositions, and an analogue of Artin’s induction theorem for certain associated characters. To better understand the effect of this functor, we explicitly determine the finite induction of irreducible modules for standard subVOAs of the rank-one lattice/affine VOA $V_{A_1}$, as well as the finite induction of irreducible modules over a parabolic-type subVOA $V_P$ of the rank-two lattice/affine VOA $V_{A_2}$.

\end{abstract}

\tableofcontents

\allowdisplaybreaks

\section{Introduction}

	Due to the complexity of the structure of VOAs and the fact that a VOA has infinitely many products, it is extremely difficult to study induced modules for a VOA embedding $U \hookrightarrow V$. In their foundational work on conformal field theory \cite{MS89}, Moore and Seiberg discovered a close analogy between group theory and the theory of chiral algebras (VOAs). Since induction plays a central role in the representation theory of finite groups, this naturally suggests the existence of an induction functor for VOAs. Understanding how modules behave under VOA embeddings is fundamental for the study of extensions, orbifolds, and coset constructions in conformal field theory. Dong and Lin made early progress in this direction \cite{DLin96}. 
	However, concrete examples for standard VOA embedding of such inductions remain unknown. 
	
	Zhu’s associative algebra $A(V)$~\cite{Z} provides a significant simplification of the representation theory of VOAs and, in turn, potentially offers an alternative approach to constructing induced modules for VOAs.
Inspired by the author's previous work on classical Yang-Baxter equations on VOAs \cite{BGL25,BGLW23} and the Borel- and parabolic-type VOAs \cite{Liu24}, in this paper, we introduce a induced module functor for VOAs on the Zhu algebra \cite{Z} level, which we call the {\em finite induction functor} for VOAs. This functor can be explicitly determined as long as we have enough knowledge about Zhu algebras $A(U)$, $A(V)$ for a VOA embedding $U\hr V$, as well as the kernel of the algebra homomorphism $\pi: A(U)\ra A(V)$ induced by the VOA embedding. Here, we do not need the VOA embedding $U\hr V$ to be conformal. i.e., the Virasoro element $\om_U$ need not be the same as  $\om_V$. 
	
The finite induction functor for VOAs is defined by the following commutative diagram of functors between categories: 
		\begin{equation}\label{eq:inddiagram0}
			\begin{tikzcd}[scale=1.5,row sep=large, column sep = large]
				\mathsf{Adm}(U)\arrow[r,dashed,shift left=2pt,"\mathrm{Ind}^V_U"]\arrow[d,shift left=2pt,"\Om_U"]& \mathsf{Adm}(V)\arrow[l,dashed,shift left=2pt,"\mathrm{Res}^V_U"]\arrow[d,shift left=2pt,"\Om_V"]\\
				\mathsf{Mod}(A(U))\arrow[r,shift left=3pt,"\mathrm{Ind}^{A(V)}_{A(U)}"]\arrow[u,shift left=2pt,"\Phi_U^\LL"]& \mathsf{Mod}(A(V)), \arrow[l,shift left=2pt,"\mathrm{Res}^{A(V)}_{A(U)}"]\arrow[u,shift left=2pt,"\Phi_V^\LL"]
			\end{tikzcd}
		\end{equation}
where $\adm(U)$ and $\adm(V)$ are the admissible (or $\N$-gradable) module categories of VOAs, $	\mathsf{Mod}(A(U))$ and $\mathsf{Mod}(A(V))$ are the module category of Zhu algebras, see Definition~\ref{df:degreezeroind}. 


		\subsubsection*{Finite versus affine algebras} We first explain the term ``finite-induction functor''. 
	In the theory of Kac-Moody algebras \cite{Kac90}, the finite-dimensional simple Lie algebra $\g$ and its affinization $\hat{\g}$ are referred to as the {\em finite} and {\em affine} Kac-Moody algebras, respectively. In the representation theory of VOAs, one can view a VOA $V$ itself or its enveloping algebra $\scrU(V)$ as an affine-type algebra similar to $\hat{\g}$, with its Zhu algebra $A(V)$ playing the role of the finite-type algebra $\g$ in the representation theory of $V$ \cite{FZ92}. This viewpoint was also adopted by De Sole and Kac in their study of $W$-algebras \cite{DSK06,DSKV16}.  
	The following diagram illustrates examples of finite and affine algebras in the VOA setting:  
	\begin{equation}\label{eq:finitvsaffine}
		\begin{tikzcd}
			\mathrm{Affine\ level\ VOA}\ V: \arrow[d,"A(-)"] & W(\g,f)\arrow[d,"A(-)"] & V_{\hat{\g}}(k,0) \arrow[d,"A(-)"] & L_{\g}(k,0) \arrow[d,"A(-)"]\\
			\mathrm{Finite\ level\ algebra}\ A(V): & W^{\mathrm{fin.}}(\g,f) & U(\g) & U(\g)/\<e_\theta^{k+1}\>.
		\end{tikzcd}
	\end{equation}
$W(\g,f)$ is the affine $W$-algebra associated with a nilpotent element $f\in \g$ \cite{A15,DSK06}, whose Zhu algebra is the finite-type $W$-algebra $W^{\mathrm{fin.}}(\g,f)$ \cite{Lo10}. 

This finite-versus-affine perspective motivates our construction of induced modules for VOAs at the level of finite algebras, as a complement to the construction of (co)induced modules at the affine-algebra level~\cite{DLin96}.

	\subsection{Finite induction functor for VOAs}
The embedding of VOAs $U \hookrightarrow V$ leads, by definition, to a homomorphism between their Zhu algebras:
\[
\begin{tikzcd}
	0 \arrow[r] & \ker(\pi) \arrow[r] & A(U) \arrow[r, "\pi"] & A(V).
\end{tikzcd}
\]
Unlike the situation for finite groups or Lie algebras, $\ker(\pi)$ is nonzero in most cases. We first prove that the natural restriction-of-scalars functor $\mathrm{Res}^{A(V)}_{A(U)}: \mathsf{Mod}(A(V)) \to \mathsf{Mod}(A(U))$ associated with the algebra homomorphism $\pi$ admits a left adjoint functor
\[
\Ind^{A(V)}_{A(U)}: \mathsf{Mod}(A(U)) \to \mathsf{Mod}(A(V)), \quad 
\Ind^{A(V)}_{A(U)}(\Omega_M) := A(V) \otimes_{\pi(A(U))} (\Omega_M / \ker(\pi) \cdot \Omega_M),
\]
which satisfies the functorial property with respect to compositions, see Propositions~\ref{prop:Frobenius} and~\ref{prop:iterationind}.

\subsubsection{Definition of finite induction functor}

Given an admissible (i.e., $\mathbb{N}$-gradable) $V$-module $W = \bigoplus_{n=0}^\infty W(n)$, the space of “highest-weight vectors”
$
\Omega(W) = \{\, v \in W : a_n v = 0,\ \mathrm{wt}\, a - n - 1 < 0 \,\}
$
is a module over the Zhu algebra $A(V)$~\cite{Z}. In particular, $\Omega(\cdot)$ can be viewed as a functor $\Omega: \mathsf{Adm}(V) \to \mathsf{Mod}(A(V))$. Dong-Li-Mason constructed a left adjoint functor to $\Omega$, called the \emph{generalized Verma module functor}~\cite{DLM1}, denoted by $\bar{M}(\cdot)$. Using properties of the universal enveloping algebra of a VOA~\cite{FZ92,FBZ04}, Damiolini-Gibney-Krashen provided an alternative construction of $\bar{M}(\cdot)$, denoted by $\Phi^\LL(\cdot)$. In particular,
\begin{equation}\label{eq:introadjointpari}
	(\Phi^\LL \dashv \Omega): \mathsf{Mod}(A(V)) \rightleftarrows \mathsf{Adm}(V)
\end{equation}
is a pair of adjoint functors describing the relationship between the representation theories of finite and affine algebras~\eqref{eq:finitvsaffine}. 
It was proved in~\cite{DGK2,GGKL25,DLM1} that~\eqref{eq:introadjointpari} is an adjoint equivalence between categories if and only if $V$ satisfies the \emph{strongly unital condition} for its mode transition algebras $\mathfrak{A}_d$. This condition holds when the VOA $V$ is rational~\cite{DGK23}. Moreover, the Heisenberg VOA and certain tensor products are examples of irrational VOAs satisfying the strongly unital condition for mode transition algebras~\cite{DGK2,LS25}.

	For the VOA embedding $U\hookrightarrow V$, when both $U$ and $V$ satisfy the strongly unital condition for mode transition algebras, then the adjoint pair \eqref{eq:introadjointpari} is an adjoint equivalence between categories, and we may lift the adjoint pair of functors on the associative algebra level $(\Ind^{A(V)}_{A(U)}\dashv \Res^{A(V)}_{A(U)}):  \mathsf{Mod}(A(U))\rightleftarrows  \mathsf{Mod}(A(V))$ through the vertical equivalence of category functors in diagram \eqref{eq:inddiagram0} to an adjoint pair of functors $(\Ind^V_U\dashv \Res^V_U): \adm(U)\rightleftarrows\adm(V)$. In particular, for any $M\in \adm(U)$, we have 
\begin{equation}\label{eq:introfiniteind}
\Ind^V_U(M):=\Phi_V^\L\left(A(V)\o_{A_U}\frac{\Om_U(M)}{\ker(\pi).\Om_U(M)}\right),
\end{equation}	
	where $A_U=\pi(A(U))\leq A(V)$ is a subalgebra. We call $\Ind^V_U$ and $ \Res^V_U$ the {\em finite induction and restriction} functors for VOAs, see Definition~\ref{df:degreezeroind}. 

\subsubsection{Properties and examples of the finite induction functor}
The usual properties of induced modules, including Frobenius reciprocity and the functorial property for compositions, are satisfied by the finite induction functor; see Propositions~\ref{prop:Frobenius} and~\ref{prop:iterationind}. 

\begin{customthm}{A}\label{main:A}
	If $U$ and $V$ both satisfy the strongly unital property for their mode transition algebras, then for any $M \in \adm(U)$ and $W \in \adm(V)$, we have
	\[
	\Hom_{\adm(V)}(\Ind^V_U(M), W) \cong \Hom_{\adm(U)}(M, \Res^V_U(W)).
	\]
	Moreover, let $U_1 \hookrightarrow U_2 \hookrightarrow V$ be consecutive embeddings of VOAs. Then we have 
	\[
	\Ind^V_{U_2} \circ \Ind^{U_2}_{U_1} = \Ind^V_{U_1}, 
	\quad \text{and} \quad 
	\Res^{V}_{U_2} \circ \Res^{U_2}_{U_1} = \Res^{V}_{U_1}.
	\]
\end{customthm}

There are many interesting examples of VOA embeddings $U \hookrightarrow V$, such as the orbifold embedding $V^G \hookrightarrow V$. In this paper, we focus on examples where the Zhu algebras $A(U)$ and $A(V)$ can be explicitly determined in terms of generators and relations. The simplest nontrivial example is the rank-one lattice VOA $V = V_{A_1}$ associated with the root lattice $A_1$, which is isomorphic to the affine VOA $L_{\widehat{\mathfrak{sl}_2}}(1,0)$~\cite{FK80}. The Zhu algebra $A(V_{A_1})$ is isomorphic to $U(\mathfrak{sl}_2)/\langle e^2 \rangle$~\cite{FZ92}. Using this presentation of $A(V_{A_1})$, we determine the finite induction of irreducible modules over various subVOAs of $V_{A_1}$; see Propositions~\ref{ex:inducedmodulesforHeisenberg}, \ref{eq:indforV_B}, and~\ref{prop:rankonevir}:

\begin{enumerate}
	\item For the Heisenberg VOA embedding $M_{\widehat{\C \alpha}}(1,0) \hookrightarrow V_{A_1}$, the finite inductions of irreducible $M_{\widehat{\C \alpha}}(1,0)$-modules are given by
$
	\mathrm{Ind}_{M_{\widehat{\C \alpha}}(1,0)}^{V_{A_1}} (M_{\widehat{\C \alpha}}(1,0)) \cong V_{A_1},$	and
	\begin{equation}\label{eq:introheisenberg}
		\begin{aligned}
			&\mathrm{Ind}_{M_{\widehat{\C \alpha}}(1,0)}^{V_{A_1}} (M_{\widehat{\C \alpha}}(1,\pm \alpha/2)) \cong V_{A_1 + \frac{1}{2}\alpha},\\
			&\mathrm{Ind}_{M_{\widehat{\C \alpha}}(1,0)}^{V_{A_1}} (M_{\widehat{\C \alpha}}(1,\lambda)) = 0,
			\quad \lambda \in \C \alpha \setminus \{0, \pm \alpha/2\}.
		\end{aligned}
	\end{equation}
	
	\item For the Borel-type VOA embedding $V_B = \bigoplus_{n \ge 0} M_{\widehat{\C \alpha}}(1, n\alpha) \hookrightarrow V_{A_1}$~\cite{Liu24}, the finite inductions of irreducible $V_B$-modules are given by $
	\mathrm{Ind}_{V_B}^{V_{A_1}} (M_{\widehat{\C \alpha}}(1,0)) \cong V_{A_1},
$
	and
	\begin{equation}\label{eq:introBorel}
		\begin{aligned}
			&\mathrm{Ind}_{V_B}^{V_{A_1}} (M_{\widehat{\C \alpha}}(1, \alpha/2)) \cong V_{A_1 + \frac{1}{2}\alpha},\\
			&\mathrm{Ind}_{V_B}^{V_{A_1}} (M_{\widehat{\C \alpha}}(1, \lambda)) = 0,
			\quad \lambda \in \C \alpha \setminus \{0, \alpha/2\}.
		\end{aligned}
	\end{equation}
	
	\item For the Virasoro VOA embedding $L(1,0) \hookrightarrow V_{A_1}$, the finite inductions of irreducible $L(1,0)$-modules are given by
$
	\Ind^{V_{A_1}}_{L(1,0)} (L(1,0)) \cong V_{A_1},
$
	and
	\begin{equation}\label{eq:introVir}
		\begin{aligned}
			&\Ind^{V_{A_1}}_{L(1,0)} (L(1,1/4)) \cong V_{A_1 + \frac{1}{2}\alpha} \oplus V_{A_1 + \frac{1}{2}\alpha},\\
			&\Ind^{V_{A_1}}_{L(1,0)} (L(1,k)) = 0,
			\quad k \in \C \setminus \{0, 1/4\}.
		\end{aligned}
	\end{equation}
\end{enumerate}

To examine the relation between the finite inductions and representation of Lie algebras, we then study the finite induction for the affine VOA embedding $V_{A_1}=L_{\widehat{\sl_2}}(1,0)\hr  L_{\widehat{\sl_3}}(1,0)=V_{A_2}$, which is given by the embedding of root lattice $A_1=\Z\al\hr \Z\al\op \Z\b=A_2$, see Figure~\ref{fig1}. In this case, both $U$ and $V$ are rational VOAs, but the the VOA embedding is not conformal. To explicitly determine the finite induction, we first give a concrete description of Zhu algebra $A(L_{\widehat{\sl_3}}(1,0))\cong U(\sl_3)/\<x_{\al+\b}^2\>$. Using the Serre's relation and some basic results in Lie algebras, we found a presentation of this associative algebra in terms of generators and relations \eqref{rel1}--\eqref{eq:morerel3}. In this case, kernel of the algebra homomorphism $\pi:A(U)\ra A(V)$ is actually zero, see Proposition~\ref{prop:embeddingofA1inA2}. We have following result about the finite induction of irreducible $V_{A_1}$-modules, see Theorem~\ref{thm:indA1A2}: 
\begin{equation}\label{eq:introA1A2ind1}
\begin{aligned}
&	\Ind^{V_{A_2}}_{V_{A_1}}(V_{A_1})\cong V_{A_2}\op V_{A_2+\la_1}\op V_{A_2+\la_2},\\ 
&	\Ind^{V_{A_2}}_{V_{A_1}}(V_{A_1+\frac{1}{2}\al})\cong V_{A_2+\la_1}\op V_{A_2+\la_2},
\end{aligned}
\end{equation}
where $\la_1,\la_2$ are the fundamental dominant weights associated to the root lattice $A_2$.

In the representation theory of finite groups, the \emph{character} $\chi_V:G\ra \C$ of a representation $\rho:G\ra \mathrm{GL}(V)$, defined by $\chi_V(g)=\tr_{V}\rho(g)$, plays a fundamental role. Artin's and Brauer's induction theorems, as well as Mackey's theory of irreducible modules, can all be interpreted in terms of characters. Since the theory of VOAs and the corresponding conformal field theory (CFT) can be viewed as a generalization of group theory \cite{MS89}, it is natural to expect a notion of character for VOAs that is compatible with induced module functors. 

The \emph{formal character} \cite{FLM,Z} 
\[
Z_W(a,\tau)=\tr_{W}\big(o(a) q^{L(0)-\frac{c}{24}}\big)
=\sum_{n\in \N} \tr_{W(n)}\big(o(a)\big) q^{n+h-\frac{c}{24}}
\]
is a natural candidate. However, we believe that this character should be viewed as an ``affine-type'' character, in view of the correspondence~\eqref{eq:finitvsaffine}, while the ``finite-type'' character is given by the top-degree coefficient $\tr_{W(0)} o(a)$ of the formal character $Z_W(a,\tau)$. Indeed, $W(0)$ is a module over the Zhu algebra $A(V)$, which is the ``finite-type'' algebra. 

Therefore, we define the \emph{finite-type character} of an ordinary $V$-module $W\in \Ord(V)$ with $\dim \Om(W)<\infty$ by 
\[
\chi_W:A(V)\ra \C,\quad \chi_W([a])=\tr_{\Om(W)} o(a).
\]
The ring generated by irreducible characters is denoted by $R^{\fin}(V)$, and is called the \emph{finite-type character ring}, see Definition~\ref{def:finitecharactr} for more details. 

Given a VOA embedding $U\hookrightarrow V$, with $V$ strongly rational, we call the character of the finite induced module $\Ind^V_U(M)$ the \emph{induced character}, denoted 
$
\Ind^V_U (\chi^U_M):=\chi_{\Ind^V_U(M)}.
$
By generalizing the argument on the finite induction of Virasoro VOA in $V_{A_1}$ \eqref{eq:introVir}, we obtain an analogue of Artin's induction theorem for finite-type characters (see Theorem~\ref{thm:Artininduction}). 

\begin{customthm}{B}\label{main:BB}
	Let $V$ be a strongly rational VOA. Assume that the conformal weights $h_0,\dots,h_r\in \C$ of its irreducible modules are pairwise distinct. Then every element of $R^{\mathrm{fin}}(V)$ can be written as a $\Q$-linear combination of induced characters from subVOAs of $V$.
\end{customthm}

\subsection{Parabolic-type subVOA $V_P$ of $V_{A_2}$}
The rest of the paper focuses on a typical example of VOA embedding $V_P\hookrightarrow V_{A_2}$ which generalizes the embedding of a parabolic subalgebra into a semisimple Lie algebra and the rank-one Borel-type subVOA embedding $V_B\hr V_{A_1}$ \eqref{eq:introBorel}. Our goal is to determine the finite induction for the VOA embedding $V_P\hr V_{A_2}$. 

For the type-$A_2$ root lattice $A_2=\Z\al\op \Z\b$ and the associated lattice VOA $V_{A_2}$ \cite{FLM}, $P:=\Z\al\op\Z_{\geq 0}\b$ is a parabolic-type submonoid of $A_2$, and $V_P=\bigoplus_{\ga\in P}M_{\hat{\h}}(1,\ga)$ is a parabolic-type subVOA of $V_{A_2}$ associated to $P$ \cite[Definition 3.4]{Liu24}. Note that $V_P$ is a CFT-type, $C_1$-cofinite, non-simple, and irrational VOA, with the degree one Lie algebra $(V_P)_1=\fp$ given by the typical parabolic subalgebra of $\sl_3=(V_{A_2})_1$ spanned by the following matrices:
$$
\begin{bmatrix}
	\ast &\ast &\ast \\
	0&\ast &\ast\\
	0&\ast &\ast
\end{bmatrix}.
$$
In particular, the VOA embedding $V_P\hookrightarrow V_{A_2}$ generalizes the Lie algebra embedding $\fp\hookrightarrow \sl_3$, and the finite induction with respect to this VOA embedding generalizes the lowest-rank parabolic induction of the type-$A$ Lie algebras.

\subsubsection{Structure of the Zhu algebra $A(V_P)$}

Since our finite induction functor \eqref{eq:introfiniteind} heavily depends on the structure of Zhu algebras, we first determine the structure of the Zhu algebra $A(V_P)$ in terms of generators and relations in Section~\ref{sec:6}. Our strategy here is similar to the computation of the Zhu algebra $A(V_B)$ for the rank-one Borel-type subVOA $V_B\leq V_{\Z\al}$ in \cite[Section~6]{Liu24}. However, since $V_P$ is essentially a rank-two object, our argument is not a straightforward generalization of the Borel-type VOA $V_B$, and the structure of $A(V_P)$ turns out to be much more intricate than that of $A(V_B)$. 

We begin by giving a detailed description of $O(V_P)$. Using the definition of lattice vertex operators, we identify a subset $O$ (see~\eqref{4.7}) of $O(V_P)$ with explicitly described spanning elements, as developed in Lemmas~\ref{lm4.1}, \ref{lm4.2}, and Proposition~\ref{prop4.4}. Then, through an inductive argument, we show that $O=O(V_P)$. This constitutes the most technical part of the determination of $A(V_P)$, with the details presented in Propositions~\ref{prop4.5} and~\ref{prop4.8}. 

The following is our third main theorem, see Definition~\ref{df4.3} and Theorem~\ref{prop4.7}.

\begin{customthm}{C}\label{main:C}
	Let $V_P=V_{\Z\al\op \Z_{\geq 0} \b}$ be the parabolic-type subVOA of $V_{A_2}$. Then $O(V_P)$ is spanned by the following elements:
	\[
	\begin{cases}
		h(-n-2)u+h(-n-1)u, & u\in V_P,\ h\in \h,\ n\geq 0;\\
		\ga(-1)v+v, & v\in M_{\hat{\h}}(1,\ga),\ \ga\in \{ \al,-\al,\b,\al+\b\};\\
		\ga(-1)^2v+\ga(-1)v, & v\in M_{\hat{\h}}(1,\ga+\ga'),\ \ga,\ga'\in \{\al,-\al,\b,\al+\b\},\ \ga+\ga'\in \{\al+\b,\b\};\\
		M_{\hat{\h}}(1,m\al+n\b), & m\al+n\b \in (\Z\al\op \Z_{\geq 0}\b)\setminus\{0,\al,-\al,\b,\al+\b\};\\
		\al(-1)^3w-\al(-1)w, & w\in M_{\hat{\h}}(1,0).
	\end{cases}
	\]
\end{customthm}
With the description of $O(V_P)$ in Theorem~\ref{main:C}, we can easily show that $A(V_P)$ is generated, as an associative algebra, by the equivalence classes of five elements:
\[
\al(-1)\vac,\quad \b(-1)\vac,\quad e^{\pm \al},\quad e^\b,\quad e^{\al+\b}.
\]
The relations among these generators can again be determined from the definition of the lattice vertex operators. This leads to our fourth main theorem; see Definition~\ref{def:relations} and Theorem~\ref{thm4.10}.

\begin{customthm}{D}\label{main:D}
	There is an isomorphism of (unital) associative algebras
	\begin{equation}\label{eq:structure}
		\begin{aligned}
			F:\; A_P
			&=\C\<x,y,x_\al,x_{-\al},x_\b,x_{\al+\b}\>/R
			\;\longrightarrow\; A(V_P),\\
			x &\mapsto [\al(-1)\vac],\quad 
			y \mapsto [\b(-1)\vac],\quad 
			x_{\pm \al} \mapsto [e^{\pm \al}],\quad 
			x_\b \mapsto [e^\b],\quad 
			x_{\al+\b} \mapsto [e^{\al+\b}],
		\end{aligned}
	\end{equation}
	where $R$ is the two-sided ideal of the tensor algebra $\C\<x,y,x_\al,x_{-\al},x_\b,x_{\al+\b}\>$ generated by the relations \eqref{4.21}--\eqref{4.26}. We use the same symbols for the equivalence classes of $x,y,x_\al,x_{-\al},x_\b,$ and $x_{\al+\b}$ in the quotient space $A_P$.
\end{customthm}

The structural theorem~\eqref{eq:structure} for $A(V_P)$ also shows that 
\[
A(V_P)\cong A(V_{A_1})[y; \Id;\delta]\op J
\]
as associative algebras, where $A(V_{A_1})[y; \Id;\delta]$ denotes the skew-polynomial algebra \cite{GW04,O33} with coefficients in the Zhu algebra $A(V_{A_1})$ of the rank-one lattice VOA $V_{A_1}$, and $J\subset A(V_P)$ is a two-sided ideal satisfying $J^2=0$ (see Corollaries~\ref{coro4.11} and~\ref{coro4.15}). Hence, $A(V_P)$ is a nilpotent extension of a skew-polynomial algebra. 

\subsubsection{Classification of irreducible $V_P$-modules}
To study the finite induction for the embedding $V_P\hookrightarrow V_{A_2}$, we first classify the irreducible $V_P$-modules $W$ and determine which of them are inducible. Since we now have an explicit description of the Zhu algebra $A(V_P)$ and there is a one-to-one correspondence between irreducible $A(V_P)$-modules and irreducible admissible $V_P$-modules \cite[Theorem~2.2.2]{Z}, Theorem~\ref{main:D} enables a complete classification of the irreducible $V_P$-modules in Section~\ref{sec:7}. 

We first construct two irreducible $V_P$-modules associated to an element $\la\in (\C\al)^\perp\subset \h=\C\al\oplus \C\b$:
\[
L^{(0,\la)}=\bigoplus_{n\in \Z} M_{\hat{\h}}(1,n\al)\otimes \C e^\la,\qquad 
L^{(\frac{1}{2}\al,\la)}=\bigoplus_{n\in \Z} M_{\hat{\h}}(1,n\al+\tfrac{1}{2}\al)\otimes \C e^\la.
\]
The $V_P$-module vertex operator $Y_M$ on $L^{(0,\la)}$ and $L^{(\frac{1}{2}\al,\la)}$ is given by slight variations of the lattice vertex operators in~\cite{FLM} (see Definition~\ref{df5.3}). Using Theorem~\ref{main:D}, we show that the bottom-degree subspaces
\[
U^{(0,\la)}=\Om_{V_P}(L^{(0,\la)}),\qquad 
U^{(\frac{1}{2}\al,\la)}=\Om_{V_P}(L^{(\frac{1}{2}\al,\la)})
\]
exhaust all possible irreducible $A(V_P)$-modules as $\la$ varies in $(\C\al)^\perp\subset \h$. The following is our fifth main theorem (see Theorem~\ref{thm5.6} and Corollary~\ref{coro5.7}):

\begin{customthm}{E}\label{main:E}
	The set
$
	\Sigma(P)=\left\{ (L^{(0,\la)},Y_M),\ (L^{(\frac{1}{2}\al,\la)},Y_M)\; :\; \la\in (\C\al)^\perp\subset \h \right\}
$	forms a complete list of irreducible modules over the rank-two parabolic-type subVOA $V_P$ of $V_{A_2}$. 
\end{customthm}

\subsubsection{Finite induction of irreducible $V_P$-modules under $V_P\hr V_{A_2}$}

With this classification theorem for irreducible $V_P$-modules, we can finally determine the finite induction $\mathrm{Ind}_{V_P}^{V_{A_2}}W$ for the irreducible $V_{P}$-modules $W$.  

The generators-and-relations presentation of $A(V_P)$ in~\eqref{4.21}--\eqref{4.26}, together with that of $A(V_{A_2})$ in~\eqref{rel1}--\eqref{eq:morerel3}, leads to a description of $\ker(\pi)$ for the associative algebra homomorphism 
$\pi: A(V_P)\ra A(V_{A_2})$ \eqref{eq:extrarelations}. Using this description, we determine the structure of modules
$\Om_{V_P}(L^{(\epsilon,\lambda)})/\ker(\pi).\Om_{V_P}(L^{(\epsilon,\la)})$ over the algebra $A_{V_P}=\pi(A(V_P))\leq A(V_{A_2})$, for all $\epsilon = 0, \tfrac{1}{2}\alpha$ and $\lambda \in (\C\alpha)^\perp$, see Proposition~\ref{lm:inducibleV_Pmodules}. It turns out that these modules are zero except when $ (\epsilon,\lambda)\in \{(0,0),\,(0,\pm \lambda_2),\, (\tfrac{1}{2}\alpha,\pm \tfrac{1}{2}\lambda_2)\}, $
where $\la_1=\frac{1}{3}\al+\frac{2}{3}\b$ and $\la_2=\frac{2}{3}\al+\frac{1}{3}\b$ are the fundamental dominant weights in the weight lattice $A_2^\circ=\Z\la_1\op \Z \la_2$, see Lemma~\ref{lm:inducibleV_Pmodules}.
 Using the explicit structures of $A_{V_P}$ and $A(V_{A_2})$, we obtain the following final characterization of induced irreducible modules for the VOA embedding $V_P \hookrightarrow V_{A_2}$, generalizing~\eqref{eq:introBorel}, see Theorem~\ref{thm:inducedmodulesforV_P}:

\begin{customthm}{F}\label{main:F}
		The finite induction of irreducible $V_P$-modules under the VOA embedding $V_P \hookrightarrow V_{A_2}$ satisfies
\[
\Ind^{V_{A_2}}_{V_P} L^{(\epsilon,\lambda)} = 0 \quad 
\text{if } (\epsilon,\lambda)\notin \{(0,0),\,(0,\pm \lambda_2),\, (\tfrac{1}{2}\alpha,\pm \tfrac{1}{2}\lambda_2)\}.
\]
For the remaining pairs of $(\epsilon,\la)$, we have 
\begin{align*}
	&\Ind^{V_{A_2}}_{V_P}L^{(0,0)} \cong V_{A_2}, 
	&& \Ind^{V_{A_2}}_{V_P}L^{(0, \lambda_2)} \cong V_{A_2+\lambda_2}, 
	&& \Ind^{V_{A_2}}_{V_P}L^{(0,- \lambda_2)} \cong 0,\\
	&\Ind^{V_{A_2}}_{V_P}L^{(\frac{1}{2}\alpha, \frac{1}{2}\lambda_2)} \cong  V_{A_2+\lambda_1}, 
	&& \Ind^{V_{A_2}}_{V_P}L^{(\frac{1}{2}\alpha, -\frac{1}{2}\lambda_2)} \cong  0,
\end{align*}
where $0$ denotes the zero module. 
\end{customthm}

This paper is organized as follows: we first define the induced module functor for associative algebra homomorphism $\pi:A\ra B$ and prove the basic properties in Section~\ref{sec:2.2} then use it to define finite induction functor for VOA embedding $U\hr V$ in Section~\ref{sec:2.3}. In Section~\ref{sec:3}, we determine the finite induction for subVOAs in the rank-one lattice VOA $V_{A_1}$ and prove \eqref{eq:introheisenberg}, \eqref{eq:introBorel}, and \eqref{eq:introVir}. In Section~\ref{sec:4}, we determine the finite induction for the affine/lattice VOA embedding $V_{A_1}\hr V_{A_2}$ and prove \eqref{eq:introA1A2ind1}. In Section~\ref{sec:5}, we introduce the notion of finite-type characters for VOAs and prove the Artin's induction theorem. In Sections~\ref{sec:6}, \ref{sec:7}, and \ref{sec:8}, we focus on the parabolic-type subVOA $V_P$ of the rank-two lattice VOA $V_{A_2}$, and prove Theorems \ref{main:C}--\ref{main:F}. 

Throughout this paper, $\N$ represents all natural numbers, including $0$. All vector spaces are defined over $\C$, the field of complex numbers.


\section{Finite induction and restriction functors for the VOA embedding $U\hr V$}

Let $U$ and $V$ be two VOAs such that $U \hookrightarrow V$. Here we do not assume that the embedding is conformal; that is, $U$ and $V$ are not required to share the same Virasoro element. 

We aim to use Zhu's associative algebras $A(U)$ and $A(V)$~\cite{Z} to define a module induction functor from $\adm(U)$ to $\adm(V)$. However, for most interesting cases of VOA embeddings $U \hookrightarrow V$, there exists an algebra homomorphism $\pi: A(U) \to A(V)$, but in general the kernel $\ker(\pi)$ is nonzero. We first give a construction of an induced module for a general associative algebra homomorphism $\pi: B \to A$, and then apply this construction to the VOA setting.

\subsection{Preliminaries on VOAs}
For the general definitions of VOAs, we refer to the classical texts \cite{FLM,FHL,DL,LL,FZ92,Z}. Here we recall the notions of modules, Borcherds' Lie algebra, Zhu algebra, generalized Verma modules, and contragredient modules.
Throughout this paper, we assume a VOA $V$ is of CFT-type: $V=\bigoplus_{n=0}^\infty V_n$, with $V_0=\C\vac$.
\subsubsection{Module categories of VOAs}

\begin{df}\label{df:modulesoverVOA}
	Let $V$ be a VOA. An {\bf admissible $V$-module} is a $\N$-graded vector space $M=\bigoplus_{n=0}^\infty M(n)$, equipped with a linear map $Y_M(\cdot,z):V\ra \End(M)[\![z,z^{-1}]\!],\ Y_M(a,z)=\sum_{n\in \Z} a_nz^{-n-1}$ called the {\bf module vertex operator}, satisfying 
	\begin{enumerate}
		\item (truncation property) For any $a\in V$ and $u\in M$, $Y_M(a,z)u\in M((z))$. 
		\item (vacuum property) $Y_M(\vac,z)=\Id_M$. 
		\item (Jacobi identity for $Y_M$) for any $a,b\in V$ and $u\in M$, 	
		\begin{equation}\label{eq:formalJacobi}
			\begin{aligned}
				z_0^{-1}\delta\left(\frac{z_1-z_2}{z_0}\right) &Y_M(a,z_1)Y_M(b,z_2)u-z_0^{-1}\delta\left(\frac{-z_2+z_1}{z_0}\right)Y_M(b,z_2)Y_M(a,z_1)u\\
				&=z_2^{-1}\delta\left(\frac{z_1-z_0}{z_2}\right) Y_M(Y(a,z_0)b,z_2)u.
			\end{aligned} 
		\end{equation}
		
		\item ($L(-1)$-derivative property) $Y_M(L(-1)a,z)=\frac{d}{dz} Y_M(a,z)$ for any $a\in V$.
		\item (grading property) For any $a\in V$, $m\in \Z$, and $n\in \N$,  
		\begin{equation}\label{admcondition}
			a_mM(n)\ssq M(n+\wt a-m-1).
		\end{equation} 
		In other words, $a_m$ acts as a degree $\wt a-m-1$ operator on $M$. 
	\end{enumerate}
	We write $\deg v=n$ if $v\in M(n)$, and call it the {\bf degree} of $v$.  Submodules, quotient modules, and irreducible modules are defined in the usual categorical sense. A morphism between two admissible $V$-modules $M$ and $W$ is a linear map $f:M\ra W$ satisfying 
	\begin{equation}\label{eq:morphism}
		f(Y_M(a,z)u)=Y_W(a,z)f(u),\quad u\in M,\ a\in V.
	\end{equation}
	Denote the category of admissible modules by $\mathsf{Adm}(V)$. The following additional notions will be used later. 
	
	\begin{enumerate}
		\item 	An admissible $V$-module $M$ is called {\bf ordinary} if each degree-$n$ subspace $M(n)=M_{n+h}$ is a finite-dimensional eigenspace of $L(0)$ of eigenvalue $n+h$, where $h\in \C$ is called the {\bf conformal weight} of $M$. In particular, if we write $L(0)v=(\wt v)\cdot v$ for $v\in M(n)$, then $\wt v=\deg v+h$. Denote the category of ordinary modules by $\Ord(V)$. 
		\item 	More generally, a {\bf weak} $V$-module is vector space $M$, together with a module vertex operator $Y_M(\cdot,z)$, satisfying conditions $(1),(2)$, and $(3)$ above \cite{DLM1}. 
		\item 	$V$ is called {\bf rational} if $\adm(V)$ is semisimple \cite{Z,DLM1}.
		\item Let $C_2(V)=\spn\{ a_{-2}b:a,b\in V\}$. $V$ is said to be {\bf $C_2$-cofinite} if $\dim V/C_2(V)<\infty$. 
		\item $V$ is called {\bf strongly rational} if $V$ is of CFT-type, simple, $V\cong V'$, rational, and $C_2$-cofinite. 
	\end{enumerate}
\end{df}


\subsubsection{Universal enveloping algebra of a VOA}

\begin{df}\cite{B86,FZ92,DGK23}\label{Def:BorLie}
	The {\bf Borcherds' Lie algebra $L(V)$} associated to a VOA $V$ is defined by 
	$
	L(V)=(V\o \C[t,t^{-1}])/\nabla (V\o \C[t,t^{-1}])=\spn\{a_{[n]}:=a\o t^n+\Im \nabla: a\in V,n\in \Z \},
	$
	where $\nabla =L(-1)\o \Id+\Id\o (d/dt)$. The Lie bracket on $L(V)$ is given by 
	$$
	[a_{[m]},b_{[n]}]=\sum_{j\geq0} \binom{m}{j} (a_jb)_{[m+n-j]},\quad a,b\in V,\ m,n\in \Z.
	$$
	For a general spanning element $a_{[n]}=a\o t^n+\Im \nabla\in L(V)$, with $a\in V$ homogeneous, define 
	\begin{equation}\label{eq:degofLV}
		\deg (a_{[n]}):=\wt a-n-1. 
	\end{equation}
	Let $\U=U(L(V))$ be the universal enveloping algebra of the Lie algebra $L(V)$. It is a graded associative algebra  by the degree \eqref{eq:degofLV}: 
	$$
	\U=\bigoplus_{d\in \Z} \U_d,\quad \U_d=\spn\left\{a^1_{[n_1]}\ds a^r_{[n_r]}\in \U: \sum\limits_{i=1}^{r}\left(\wt a^i-n_i-1\right)=d \right\}.
	$$
	Let $\U_{\leq -n}=\sum_{d\leq -n} \U_d$, which makes $\U$ a split-filtered associative algebra $\U=\bigcup_{n\in \Z} \U_{\leq -n}$. Define
	\begin{equation}\label{eq:nei}
		\rN^n_\L \U=\U\cdot \U_{\leq -n}=\U\cdot L(V)_{\leq -n},\quad 	\rN^n_\sR \U=\U_{\geq n}\cdot \U=L(V)_{\geq n}\cdot \U.
	\end{equation}
	where $L(V)_{\leq -n}=\spn\{a_{[k]}\in L(V): \deg(a_{[k]})\leq -n\}$, see \cite[Lemma 2.4.2]{DGK23}. $L(V)_{\geq n}$ is defined in a similar way. Since the identity $1=\vac_{[-1]}$ of $\U$ is contained in $\U_{\leq 0}$ and $\U_{\geq 0}$, we have $\rN^n_\L \U=\U=\rN^n_\sR \U	$ if $n\leq 0$. 
	
	The left ideals $\{\rN^n_\L \U:n\in \Z_{\geq 0}\}$ is a system of neighborhood of $0$ in $\U$, which gives a canonical seminorm on $\U$, see \cite[Definition A.6.1]{DGK23}. One can restrict these seminorms to the graded parts $\U_d$ of $\U$: 
	\begin{equation}\label{eq:completion}
		\rN^n_\L \U_d:=(\U\cdot \U_{\leq -n})_d=\sum_{j\leq -n} \U_{d-j}\cdot \U_j,\quad  	\rN^n_\sR \U_d=(\U_{\geq n}\cdot \U)_d=\sum_{i\geq n} \U_i\cdot \U_{d-i}.
	\end{equation}
	In particular, $\rN^n_\L \U_d=\rN^{n+d}_\sR \U_d$ for any $d\in \Z$. Define the completion
	$$
	\widehat{\U}_d:=\varprojlim_n\frac{\U_d}{\rN^n_\L \U_d}=\varprojlim_n\frac{\U_d}{\rN^{n+d}_\sR \U_d}\quad \mathrm{and}\quad \widehat{\U}:=\bigoplus_{d\in \Z} \widehat{\U}_d. 
	$$ 
	Let $J\subset \widehat{\U}$ be the graded ideal generated by the components of the Jacobi identity \eqref{eq:formalJacobi}, and let $\bar{J}\ssq \widehat{\U}$ be the closure of $J$ with respect to the seminorm defined by the image of neighborhoods \eqref{eq:nei} in $\widehat{U}$. Define
	$$\scrU=\scrU(V):=\widehat{\U}/\bar{J}=\bigoplus_{d\in \Z} \scrU_d.$$
	Then $\scrU$ is a graded complete seminormed associative algebra with respect to the canonical seminorm induced by the image of neighborhoods \eqref{eq:nei}. $\scrU$ is called the {\bf universal enveloping algebra of the VOA $V$}. The left and right neighborhoods at $0$ of $\scrU$  are given by 
	\begin{equation}\label{eq:nebors}
		\rN^n_\L \scrU=\scrU\cdot \scrU_{\leq -n}\quad \mathrm{and}\quad  	\rN^n_\sR \scrU=\scrU_{\geq n}\cdot \scrU,
	\end{equation}
	with $	\rN^{n+1}_\L \scrU_0=\sum_{j \geq n+1} \scrU_j\cdot \scrU_{-j}=	\rN^{n+1}_\sR \scrU_0$, for any $n\geq 0$.
\end{df}

	\subsubsection{ Zhu algebra of a VOA}
	
	\begin{df}\cite{Z}\label{def:AV}
		Let $V$ be a VOA, the {\bf Zhu algebra $\A=A(V)$} is defined as a quotient space  $\A=V/O(V)$, where 
		\begin{equation}\label{eq:defofOV}
			O(V)=\spn\left\{ a\circ b=\Res_{z=0}Y(a,z)b\frac{(1+z)^{\wt a}}{z^2}: a,b\in V\right\}.
		\end{equation}
		$\A=\spn\{[a]=a+O(V):a\in V\}$ is an associative algebra with respect to product
		\begin{equation}\label{eq:prodAV}
			[a]\ast [b]=\Res_{z=0} [Y(a,z)b]\frac{(1+z)^{\wt a}}{z}=\sum_{j\geq 0}\binom{\wt a}{j} [a_{j-1}b],\quad a,b\in V.  
		\end{equation}
		Denote the category of left $\A$-modules by $\mathsf{Mod}(\A)$. 
	\end{df}
	One can show that $\A\cong \scrU_0/\rN^1_\L\scrU_0$ as associative algebras \cite{FZ92}. 
	
	\subsubsection{Generalized Verma module}

	Let $W$ be a weak $V$-module. Then the space of ``highest-weight vectors'' in $W$
	\begin{equation}\label{eq:2.7}
		\Om(W)=\spn\{v\in W: a_nw=0,\ \deg (a_n)=\wt a-n-1<0 \}
	\end{equation}
	is a left $A(V)$-module via the representation map 
	$$
	\A\ra \End(\Om(W)),\ [a]\mapsto o(a)=a_{\wt a-1},\ a\in V.
	$$
	We can view $\Om$ as a functor 
	\begin{equation}\label{eq:defOm}
		\Om: \mathsf{Adm}(V)\ra \mathsf{Mod}(\A),\quad W\mapsto \Om(W),
	\end{equation}
	which is an one-to-one correspondence between irreducible objects in these categories, see \cite[Theorem 2.2.2]{Z}. 
	
	On the other hand, Dong-Li-Mason's generalized Verma module functor $\bar{M}(\cdot)$ \cite[Theorem 6.2]{DLM1} can be identified with the following (left) induced module functor:
	\begin{equation}\label{eq:PhiL}
		\Phi^\L:  \mathsf{Mod}(\A)\ra  \mathsf{Adm}(V),\quad S\mapsto \Phi^\L(S)=(\scrU/	\rN^1_\L \scrU)\o _{\scrU_0} S. 
	\end{equation}
$\Phi^\L(S)$ is called the {\bf generalized Verma module} associated to $A(V)$-module $S$. 	
	The functors $\Om$ and $\Phi^\L$ in \eqref{eq:defOm} and \eqref{eq:PhiL} form an adjoint pair $\Phi^\L\dashv \Om$ between abelian categories by the universal property of generalized Verma modules
	\begin{equation}\label{eq:adjointfunctor}
		\Phi^\L:\mathsf{Mod}(\A)\rightleftarrows\mathsf{Adm}(V):\Om
	\end{equation}
	In other words, there exists a natural isomorphism of vector spaces:
	$$
	\Hom_{\mathsf{Adm}(V)}(\Phi^\L(S),W)\cong \Hom_{\mathsf{Mod}(\A)}(S,\Om(W)),
	$$
	see \cite[Theorem 6.2]{DLM1} and \cite[Proposition 3.1.2]{DGK23}. 
	
\subsubsection{Mode transition algebras}
	In \cite{DGK23}, a sequence of associative algebras $\fA_d$, called the {\bf mode transition algebras}, were introduced to study the smoothing property of the sheaves of VOA-conformal blocks on $\overline{\mathcal{M}}_{g,n}$. They fit into the following exact sequence of associative algebras:
	$$
	\begin{tikzcd}
	\fA_d\arrow[r] &A_d(V)\arrow[r] &A_{d-1}(V)\arrow[r]& 0,
	\end{tikzcd}
	$$
	where $A_d(V)$, with $d\geq 1$, are the higher-level Zhu algebras \cite{DLM2}. 
	
An element 	$1_d\in \fA_d$ is called a {\bf strong unit} if $1_d\ast \al=\al$ and $\b\ast 1_d=\b$, for all $\al\in \fA_{d,0}$ and $\b\in \fA_{0,-d}$. The VOA $V$ is said to satisfy the {\bf strongly unital property} if the mode transition algebras $\fA_d$ are all strongly unital. 

\begin{example}
	It was proved in \cite{DGK23} that if $V$ is a rational VOA, then it satisfies the strongly unital property.
	
Also, the Heisenberg VOA $V=M_{\hat{\h}}(k,0)$ satisfies the strongly unital property  \cite{DGK23,DGK2}. Moreover, we proved that if $V_1$ and $V_2$ both satisfy the strongly unital property, so does the tensor product VOA $V_1\o V_2$ \cite{GGKL25,LS25}. 
\end{example}

	\begin{lm}\cite{DGK2,GGKL25}\label{lm:adjointequi}
	If the VOA $V$ satisfies the strongly unital property, then the adjoint pair $(\Phi^\L \dashv \Om):\mathsf{Mod}(\A)\rightleftarrows\mathsf{Adm}(V)$ \eqref{eq:adjointfunctor} is an adjoint equivalence between categories. In this case, any admissible $V$-module $W$ is a generalized Verma module $W=\Phi^\L_V(\Om(W))$.  
		In particular, $(\Phi^\L\dashv \Om)$ is an adjoint equivalence if $V$ is rational. 
	\end{lm}

	\subsection{Induction functor for the associative algebra homomorphism}\label{sec:2.2}
	Let $A,B$ be associative unital algebras over $\C$, and let
	\begin{equation}\label{eq:pi}
		\begin{tikzcd}
			0\arrow[r]& \ker(\pi)\arrow[r]& B\arrow[r,"\pi"]& A
		\end{tikzcd}
	\end{equation}
	be an exact sequence in the category of associative algebras. Denote the categories of left $A$ (resp. $B$)-modules by $\Mod(A)$ (resp. $\Mod(B)$). 
	\subsubsection{Construction of induced modules}
	Here we give a natural definition/construction of the induced modules for the algebra homomorphism \eqref{eq:pi}: 
	
	\begin{enumerate}
		\item $\pi(B)\leq A$ is a unital subalgebra, $\ker(\pi)\lhd B$ is a two-sided ideal, and $A$ is a $(A,\pi(B))$-bimodule. 
		\item
		There is a natural {\bf restriction of scalar functor} from $\Mod(A)$ to $\Mod(B)$ via $\pi$
		$$\Res^A_B(-):\Mod(A)\ra \Mod(B),$$
		which is an exact functor between these abelian categories.  
		\item  Conversely, let $\Om_M$ be an object in $\Mod(B)$, then $\Om_M/\ker(\pi).\Om_M$ is a left module over $\pi(B)\cong B/\ker(\pi)$, with $\pi(b).\bar{v}=\overline{b.v}$, for any $b\in B$ and $v\in \Om_M$. Define
		\begin{equation}\label{eq:defassocind}
			\Ind_B^A (\Om_M):=A\o_{\pi(B)} \left(\frac{\Om_M}{\ker(\pi).\Om_M}\right). 
		\end{equation}
		Then $	\Ind_B^A (\Om_M)$ is an object in $\Mod(A)$. By the definition of tensor product modules, it is easy to show that $	\Ind_B^A (\Om_M)=A\o_B \Om_M$ if we view $A$ as a right $B$-module via $\pi$. 
		Moreover, there exists a left $B$-module homomorphism 
		$$
		\iota:\Om_M\ra \Ind_B^A(\Om_M),\quad \iota(v)=1_A\o\bar{v}.
		$$
		$\Ind^A_B(\Om_M)$ is generated by $\iota(\Om_M)$ as a left $A$-module. 
		\item Let $f:\Om_M\ra \Om_W$ be a morphism in $\Mod(B)$. Since $f$ commutes with the action of $\ker(\pi)$, the map
		\begin{equation}\label{eq:indonmorphism}
			\Ind_B^A(f)=\Id_A\o \bar{f}:	\Ind_B^A (\Om_M)\ra 	\Ind_B^A (\Om_W),\quad a\o\bar{v}\mapsto a\o \overline{f(v)}
		\end{equation}
		is a well-defined morphism in $\Mod(A)$. Then $$\Ind_B^A(-):\Mod(B)\ra\Mod(A)$$ is a functor between these abelian categories. 
	\end{enumerate}
	If $\ker (\pi)=0$ in \eqref{eq:pi}, then $\Ind_B^A (\Om_M)=A\o_B\Om_M$ recovers Higman's definition of induced modules for the associative algebra embedding $B\hr A$ \cite{Hig55}. 
	\subsubsection{The Frobenius reciprocity} The induction defined by \eqref{eq:defassocind} and \eqref{eq:indonmorphism} satisfies the usual Frobenius reciprocity. In other words, $\Ind_B^A(-)$ is a left adjoint of the functor $\Res^A_B(-)$.
	\begin{prop}\label{prop:Frobenius}
		Let $\pi: B\ra A$ be a homomorphism of associative algebras. Then 
		\begin{equation}
			\Ind_B^A(-):\Mod(B)\rightleftarrows \Mod(A): \Res_B^A(-)
		\end{equation}
		is an adjoint pair of functors. In other words, there exists a canonical isomorphism of the hom-spaces:
		\begin{equation}
			\Hom_A(\Ind_B^A(\Om_M),K)\cong \Hom_B(\Om_M,\Res^A_B(K)),
		\end{equation}
		for any left $B$-module $\Om_M$ and any left $A$-module $K$. 
	\end{prop}
	\begin{proof}
		Define $\Phi:  \Hom_B(\Om_M,\Res^A_B(K))\ra \Hom_A(\Ind_B^A(\Om_M),K)$ by
		$$\Phi(g)(a\o \bar{v}):=a.g(v),\quad a\in A, v\in \Om_M,\ g\in  \Hom_B(\Om_M,\Res^A_B(K)). $$ 
		For $b\in \ker(\pi)$, we have $a\o \overline{b.v}=0$ in $\Ind^A_B(\Om_M)$ \eqref{eq:defassocind}. Note that $\Phi(g)(a\o \overline{b.v})=a.g(b.v)=a.(\pi(b).g(v))=0$ by the definition of $\Res^A_B(K)$. Clearly, $\Phi(g)$ is a left $A$-module homomorphism. Hence $\Phi$ is well-defined. 
		
		Conversely, define $\Psi:  \Hom_A(\Ind_B^A(\Om_M),K)\ra \Hom_B(\Om_M,\Res^A_B(K))$ by 
		$$\Psi(f)(v):=f(1\o \bar{v}),\quad v\in \Om_M,\ f\in   \Hom_A(\Ind_B^A(\Om_M),K). $$
		We have $\Psi(f)(b.v)=f(1\o \overline{b.v})=f(1\o \pi(b).\bar{v})=f(\pi(b)\o \bar{v})=\pi(b).f(1\o\bar{v})=b.\Psi(f)(v)$, for any $b\in B$ and $v\in \Om_M$, by the definition of the induced and restricted modules under $\pi$ \eqref{eq:pi}. Hence $\Psi$ is well-defined. Clearly, $\Phi$ and $\Psi$ are mutually inverse to each other. 
	\end{proof}
	
	\begin{coro}
		The pair $(\Ind_B^A(\Om_M),\iota)$ satisfies the universal property: Let $K$ be a left $A$-module, and $f:\Om_M\ra \Res^A_B(K)$ be a $B$-module homomorphism, then there exists a unique $A$-module homomorphism $F:\Ind^A_B(\Om_M)\ra K$ such that $F\circ \tau=f$.  
		$$
		\begin{tikzcd}
			\Om_M\arrow[r,"\tau"]\arrow[dr,"f"']& \Ind_B^A(\Om_M)\arrow[d,"F",dashed]\\
			&K
		\end{tikzcd}
		$$
		In particular, if we define the $A$-module $\Ind^A_B(\Om_M)$ by the universal property, then it exists and is unique up to unique isomorphism. 
	\end{coro}
	
	\subsubsection{Functorial properties of the induction functor} 
	The usual iteration property of the induction is also satisfied: 
	
	\begin{prop}\label{prop:iterationind}
		Let $A,B,C$ be three associative algebras over $\C$, with homomorphisms
		$$
		\begin{tikzcd}
			C\arrow[r,"\varphi"]& B\arrow[r,"\pi"]&A.
		\end{tikzcd}
		$$
		The induction functor defined by \eqref{eq:defassocind} and \eqref{eq:indonmorphism} satisfies the composition property: $$\Ind^A_B(-)\circ \Ind^B_C(-)\cong \Ind^A_C(-)$$
		as functors from $\Mod(C)$ to $\Mod(A)$, where $\Ind^A_C(-)$ is defined via the algebra homomorphism $\pi\circ \varphi: C\ra A$. 
	\end{prop}
	\begin{proof}
		Let $\Om_M$ be an object in $\Mod(C)$. By \eqref{eq:defassocind} we have 
		\begin{equation}\label{eq:iteratedind}
			\Ind^A_B\left(\Ind^B_C(\Om_M)\right)=A\o_{\pi(B)}\left(\frac{B\o_{\varphi(C)}(\Om_M/\ker(\varphi).\Om_M)}{\ker(\pi)\o_{\varphi(C)}(\Om_M/\ker(\varphi).\Om_M)}\right).
		\end{equation}
		We need to show it is isomorphic to $\Ind_C^A(\Om_M)=A\o_{\pi\varphi(C)}(\Om_M/\ker(\pi\varphi).\Om_M)$. Since the tensor functor is right exact, we have a canonical epimorphism of left $\pi(B)$-modules: 
		$$
		\frac{B\o_{\varphi(C)}(\Om_M/\ker(\varphi).\Om_M)}{\ker(\pi)\o_{\varphi(C)}(\Om_M/\ker(\varphi).\Om_M)}\twoheadrightarrow (B/\ker(\pi))\o_{\varphi(C)} (\Om_M/\ker(\varphi).\Om_M),
		$$
		which gives rise to an epimorphism of left $A$-modules, in view of \eqref{eq:iteratedind}: 
		\begin{equation}\label{eq:deftheta}
			\begin{aligned}
				\theta: \Ind^A_B\left(\Ind^B_C(\Om_M)\right)&\ra A\o_{\pi(B)} (B/\ker(\pi))\o_{\varphi(C)} (\Om_M/\ker(\varphi).\Om_M),\\
				a\o [b\o \bar{v}]&\mapsto a\o \bar{b}\o \bar{v},\quad a\in A,\ b\in B,\ v\in \Om_M,
			\end{aligned}
		\end{equation}
		where $[b\o \bar{v}]$ is the equivalent class of $b\o \bar{v}$ modulo $\ker(\pi)\o_{\varphi(C)}(\Om_M/\ker(\varphi).\Om_M)$, and $\bar{b}$ on the right hand side is the equivalent class of $b$ in the quotient algebra  $B/\ker(\pi)$. Define
		\begin{equation}\label{eq:defkappa}
			\begin{aligned}
				\ka: A\o_{\pi(B)} (B/\ker(\pi))\o_{\varphi(C)} (\Om_M/\ker(\varphi).\Om_M)&\ra A\o_{\pi\varphi(C)}(\Om_M/\ker(\pi\varphi).\Om_M)\\
				a\o \bar{b}\o \bar{v}&\mapsto a\cdot \pi(b)\o \bar{v},
			\end{aligned}
		\end{equation}
		where $\bar{v}$ in $a\cdot \pi(b)\o\bar{v}$ is the equivalent class of $v$ in the quotient $\Om_M/\ker(\pi\varphi).\Om_M$. 
		
		To show $\ka$ is well-defined, we need to show it preserves the left-right module actions of $\pi(B)$ and $\varphi(c)$ through tensor. Indeed, by the definition of module actions on the quotient and \eqref{eq:defkappa}, 
		\begin{align*}
			\ka(a\o \pi(b_1).\bar{b}\o \bar{v})&=\ka(a\o \overline{b_1b}\o\bar{v})=a\cdot \pi (b_1b)\o \bar{v}=a\cdot \pi(b_1)\cdot \pi(b)\o \bar{v}\\
			&=\ka(a\cdot \pi(b_1)\o \bar{b}\o \bar{v}),\quad \forall b_1\in B,\\
			\ka(a\o\bar{b}.\varphi(c)\o \bar{v})&=\ka(a\o \overline{b\varphi(c)}\o \bar{v})=a\cdot \pi(b\varphi(c))\o\bar{v}=a\cdot \pi(b)\cdot \pi(\varphi(c))\o\bar{v}\\
			&=a\cdot\pi(b)\o \pi(\varphi(c)).\bar{v}=a\cdot \pi(b)\o \overline{c.v}=\ka(a\o\bar{b}\o\overline{c.v})\\
			&=\ka(a\o \bar{b}\o \varphi(c).\bar{v}),\quad \forall c\in C.
		\end{align*}
		Hence $\ka$ in \eqref{eq:defkappa} is a well-defined left $A$-module homomorphism, and 
		\begin{equation}\label{eq:definterphi}
			\begin{aligned}
				\Phi=\ka\circ \theta: \Ind^A_B\left(\Ind^B_C(\Om_M)\right)&\ra A\o_{\pi\varphi(C)}(\Om_M/\ker(\pi\varphi).\Om_M)\\
				a\o [b\o \bar{v}]&\mapsto a\cdot \pi(b)\o \bar{v}
			\end{aligned}
		\end{equation}
		is a well-defined left $A$-module homomorphism, in view of \eqref{eq:deftheta} and \eqref{eq:defkappa}. 
		
		Conversely, with the notation in \eqref{eq:iteratedind} and \eqref{eq:deftheta}, we define 
		\begin{equation}\label{eq:definterpsi}
			\begin{aligned}
				\Psi: A\o_{\pi\varphi(C)}(\Om_M/\ker(\pi\varphi).\Om_M)&\ra \Ind^A_B\left(\Ind^B_C(\Om_M)\right),\\
				a\o \bar{w}&\mapsto a\o[1_B\o\bar{w}],\quad a\in A,\ w\in \Om_M.
			\end{aligned}
		\end{equation}
		By \eqref{eq:definterpsi} and the definition of module actions, we have 
		\begin{align*}
			\Psi(a\cdot \pi(\varphi(c))\o \bar{w})&=a\cdot \pi(\varphi(c))\o [1_B\o \bar{w}]=a\o \pi(\varphi(c)).[1_B\o \bar{w}]\\
			&=a\o [\varphi(c)\cdot 1_B\o \bar{w}]=a\o [1_B\o \varphi(c).\bar{w}]=a\o [1_B\o \overline{c.w}]\\
			&=\Psi(a\o \overline{c.w})=\Psi(a\o \pi(\varphi(c)).\bar{w}),\quad \forall c\in C.
		\end{align*}
		Hence $\Psi$ is a well-defined left $A$-module homormorphism. Clearly, $\Psi$ and $\Phi$ are mutual inverse to each other, and so $\Ind^A_B(\Ind^B_C(\Om_M))\cong \Ind_C^A(\Om_M)$ as left $A$-modules. Since the definitions \eqref{eq:definterphi} and \eqref{eq:definterpsi} of $\Phi$ and $\Psi$ are canonical, they give rise to natural isomorphism between functors $\Ind^A_B(-)\circ \Ind^B_C(-)$ and $ \Ind^A_C(-)$. 
	\end{proof}
	
	\begin{remark}
		It is also clear from \eqref{eq:defassocind} and \eqref{eq:indonmorphism}  that if $\pi=\Id_A:A\ra A$ then $\Ind_A^A(-)$ is the identity functor on $\Mod(A)$. Hence our induction functor $\Ind^A_B(-)$ satisfies the functorial property for the algebras $A$ and $B$. Moreover, if $\varphi:C\ra B$ and $\pi:B\ra A$ are both embeddings, then Proposition~\ref{prop:iterationind} recovers the iteration property of the induction functor in \cite{Hig55}. 
	\end{remark}

	\subsection{Definition of finite induction and restriction functors}\label{sec:2.3}
	
	Let $(V,Y,\vac,\om_V)$ be a CFT-type VOA, and let $(U,Y,\vac,\om_U)$ be a subVOA of $V$. 
	Here we do not assume the embedding $U\hr V$ is conformal. i.e., $\om_U$ is not necessarily equal to $\om_V$. 
	
	Since $O(U)\subset O(V)$ \eqref{eq:defofOV}, we have an exact sequence of associative algebras
	\begin{equation}\label{eq:seqforA}
		\begin{tikzcd}
			0\arrow[r]& \ker(\pi)\arrow[r]& A(U)\arrow[r,"\pi"]& A(V),
		\end{tikzcd}
	\end{equation}
	where $A(U)$ and $A(V)$ are the Zhu algebra of $U$ and $V$, respectively.
	By Proposition~\ref{prop:Frobenius}, there is an adjoint pair of functors between module categories of associative algebras
	$$\Ind_{A(U)}^{A(V)}:\mathsf{Mod}(A(U))\rightleftarrows \mathsf{Mod}(A(V)):\mathrm{Res}^{A(V)}_{A(U)}.$$
	Denote the subalgebra $\pi(A(U))\leq A(V)$ by $A_U$. Then by \eqref{eq:defassocind}, for any $\Om_M\in \mathsf{Mod}(A(U))$,
	$$\Ind_{A(U)}^{A(V)}(\Om_M)=A(V)\o_{A_U}(\Om_M/\ker(\pi).\Om_M).$$
	
	\begin{df}\label{df:degreezeroind}
		Let $U\hr V$ be a VOA embedding. Define functors 
		\begin{align*}
			&	\Ind^V_U:=\Phi^\L_V\circ \Ind^{A(V)}_{A(U)}\circ \Om_U:\adm(U)\ra \adm(V),\\
			&	\Res^V_U:=\Phi^\L_U\circ \Res^{A(V)}_{A(U)}\circ \Om_V:\adm(V)\ra \adm(U).
		\end{align*}
		In particular, $	\Ind^V_U$ and $	\Res^V_U$ fit into the following diagram
		\begin{equation}\label{eq:inddiagram}
			\begin{tikzcd}[scale=1.5,row sep=large, column sep = large]
				\mathsf{Adm}(U)\arrow[r,dashed,shift left=2pt,"\mathrm{Ind}^V_U"]\arrow[d,shift left=2pt,"\Om_U"]& \mathsf{Adm}(V)\arrow[l,dashed,shift left=2pt,"\mathrm{Res}^V_U"]\arrow[d,shift left=2pt,"\Om_V"]\\
				\mathsf{Mod}(A(U))\arrow[r,shift left=3pt,"\mathrm{Ind}^{A(V)}_{A(U)}"]\arrow[u,shift left=2pt,"\Phi_U^\L"]& \mathsf{Mod}(A(V)) \arrow[l,shift left=2pt,"\mathrm{Res}^{A(V)}_{A(U)}"]\arrow[u,shift left=2pt,"\Phi_V^\L"]
			\end{tikzcd}
		\end{equation}
		More precisely, given $M\in \adm(U)$, we let 
		\begin{equation}\label{eq:ind}
			\Ind^V_U(M)=\Phi_V^\L\left(A(V)\o_{A_U}\frac{\Om_U(M)}{\ker(\pi).\Om_U(M)}\right);
		\end{equation}
		Given $W\in \adm(V)$, we let 
		\begin{equation}\label{eq:res}
			\Res^V_U(W)=\Phi_U^\L\left(\Res^{A(V)}_{A(U)}(\Om_V(W))\right).
		\end{equation}
		We call $\Ind^V_U$ (resp. $\Res^V_U$) the {\bf finite module induction (resp. restriction)} functor with respect to the VOA embedding $U\hr V$. 
	\end{df}
	
	
	\begin{prop}\label{prop:adjointpari}
		If $U$ and $V$ both satisfy the strongly unital property for their mode transition algebras $\fA_d$, then $(\Ind^V_U\dashv \Res^V_U):\adm(U)\rightleftarrows \adm(V)$ is a pair of adjoint functors between abelian categories. i.e., for any $W\in \adm(U)$ and $M\in \adm(V)$, we have
		\begin{equation}
			\Hom_{\adm(V)}(\Ind^V_U(W),M)\cong \Hom_{\adm(U)}(W,\Res^V_U(M)).
		\end{equation}
		In particular,  if $U$ and $V$ are both rational VOAs, then $(\Ind^V_U\dashv \Res^V_U)$ is an adjoint pair. 
	\end{prop}
	\begin{proof}
		If $U$ and $V$ both satisfy the strongly unital property, then $(\Phi_U^\L\dashv \Om_U):\mathsf{Mod}(A(U))\rightleftarrows\mathsf{Adm}(U)$ and $(\Phi_V^\L\dashv \Om_V):\mathsf{Mod}(A(V))\rightleftarrows\mathsf{Adm}(V)$ are adjoint equivalence between categories, see Lemma~\ref{lm:adjointequi}. Then $(\Ind^V_U\dashv \Res^V_U):\adm(U)\rightleftarrows \adm(V)$ is an adjoint pair, since it is a lifting via the vertical equivalence functors in \eqref{eq:inddiagram} of an adjoint pair $(\Ind_{A(U)}^{A(V)}\dashv \mathrm{Res}^{A(V)}_{A(U)}):\mathsf{Mod}(A(U))\rightleftarrows \mathsf{Mod}(A(V))$.   
	\end{proof}
	
	The following iteration property of the finite induction functor follows immediately from Proposition~\ref{prop:iterationind}. 
	
	\begin{prop}\label{prop:compo}
		Let $U_1\hr U_2\hr V$ be consecutive embeddings of VOAs. Then we have $\Ind^V_{U_2}\circ \Ind^{U_2}_{U_1}=\Ind^V_{U_1}$ and $\Res^{V}_{U_2}\circ \Res^{U_2}_{U_1}=\Res^{V}_{U_1}$. 
	\end{prop}

		Note that the finite restriction functor 
	$
	\Res^V_U:\adm(V)\ra \adm(U)
	$
	is {\bf not} the natural restriction of scalar functor with respect to the VOA embedding $U\hr V$
	\begin{equation}\label{eq:affineres}
		\mathrm{res}^V_U: \adm(V)\ra \adm(U),\quad M\ra \mathrm{res}^V_U(M),
	\end{equation}
	where $\mathrm{res}^V_U(M)=M$ as a vector space, and is viewed as a $U$-module via $$Y_M(\cdot,z)|_U:U\hr V\rightarrow \End(M)[\![z,z^{-1}]\!].$$ 
	
	\begin{example}
		Indeed, assume $U$ and $V$ are simple VOAs and $U\cong \Phi^\L_{U}(\C\vac)$ as $U$-modules. For instance, if $U=M_{\hat{\h}}(1,0)$ is the rank-one Heisenberg VOA, then it satisfies $U\cong \Phi^\L_{U}(\C\vac)$, see Lemma~\ref{lm:adjointequi}. $U$ is embedded in a rank-one lattice VOA $V=V_{L}$, which is a simple VOA. 
		
		Then $\Om_V(V)=\C\vac$ and $\Res^V_U(V)=\Phi^\L_U(\C\vac)\cong U$. However, $\mathrm{res}^V_U(V)=V$ as a vector space. Hence $\Res^V_U(V)\neq \mathrm{res}^V_U(V)$. 
	\end{example}	
	
	We will study the left-adjoint functor $\mathrm{ind}^V_U$ of $\mathrm{res}^V_U$ \eqref{eq:affineres} in a separate work.

	\section{Finite induction for the subVOAs in the rank-one lattice VOA $V_{A_1}$}\label{sec:3}
	For the general theory of lattice VOAs, we refer to \cite{FLM}. 
	Let $V=V_{A_1}$ be the rank-one lattice VOA associated to the root lattice $A_1=\Z\al$, with $(\al|\al)=2$. Write $\sl_2=\C e+\C h+\C f$. Then $V_{A_1}$ is also isomorphic to the level-one affine VOA $L_{\widehat{\sl_2}}(1,0)$, with $e^\al\mapsto e(-1)\vac,\al(-1)\vac\mapsto h(-1)\vac, e^{-\al}\mapsto f(-1)\vac,$ see \cite{FK80,FZ92}. 
	
	Hence $A(V_{A_1})=A(L_{\widehat{\sl_2}}(1,0))\cong U(\sl_2)/\<e^2\>$ is a $5$-dimensional semisimple algebra \cite{DLM1,DLM97}, where $\<e^2\>$ is the two-sided ideal of $A(L_{\widehat{\sl_2}}(1,0))$ generated by $e^2$, see \cite{FZ92}. 
	Then $A(V_{A_1})\cong U(\sl_2)/\<e^2\>$ is spanned by $\{1,e,f,h,h^2\}$ subject to the relations
	\begin{equation}\label{3.26}
		eh+e=0;\quad h^2-h-2fe=0;\quad fh-f=0;\quad e^2=f^2=0,
	\end{equation}
	where we use the same notations for the elements in $U(\sl_2)$ and the quotient $U(\sl_2)/\<e^2\>$. 
	
	$A(V_{A_1})$ has two irreducible highest-weight modules up to isomorphism. One is the trivial module $\C\vac$, the other is the two-dimensional irreducible $\sl_2$-module
	$$ L(1)=\C e^{\frac{1}{2}\al}\op \C e^{-\frac{1}{2}\al},\quad \mathrm{with}\quad e.e^{\frac{1}{2}\al}=0,\ f.e^{\frac{1}{2}\al}=e^{-\frac{1}{2}\al},\ h.e^{\frac{1}{2}\al}=e^{\frac{1}{2}\al}.$$
	
	\subsection{The rank-one Heisenberg embedding $M_{\widehat{\C\al}}(1,0)\hr V_{A_1}$} Note that this embedding is conformal. 
	
	\begin{lm}\label{lm:3.1}
		$A(M_{\widehat{\C\al}}(1,0))\cong \C[x]$ the polynomial algebra. The homomorphism \eqref{eq:seqforA} $\pi: \C[x]\ra U(\sl_2)/\<e^2\>$ is given by $\pi (f(x))=f(h)$, and $\ker(\pi)=\<x^3-x\>$. 
	\end{lm}
	\begin{proof}
		The first claim is well-known \cite{FZ92}, with $[\al(-1)\vac]\in A(M_{\widehat{\C\al}}(1,0))$ corresponds to $x$.
		The VOA embedding $M_{\widehat{\C\al}}(1,0)\hr V_{A_1}\cong L_{\widehat{\sl_2}}(1,0), \al(-1)\vac\mapsto h(-1)\vac$ induces the homomorphism of Zhu algebras $\pi: \C[x]\ra U(\sl_2)/\<e^2\>,\ x\mapsto h$. Clearly, $x^3-x\in \ker(\pi)$ since it follows from \eqref{3.26} that $h^3-h=0$. Moreover, $\C[x]/\ker(\pi)$ is isomorphic to the subalgebra of $U(\sl_2)/\<e^2\>$ generated by $h$, which is $3$-dimensional. Hence $\ker (\pi)=\<x^3-x\>$. 
	\end{proof}
	
	The irreducible modules in $\adm(M_{\widehat{\C\al}}(1,0))$ are $\{ M_{\widehat{\C\al}}(1,\la):\la\in \C\al\}$ \cite{FLM}; The irreducible modules in $\adm(V_{A_1})$ are $\{V_{A_1}, V_{A_1+\frac{1}{2}\al}\}$ \cite{D,FZ92}.

	\begin{prop}\label{ex:inducedmodulesforHeisenberg}
		For the VOA embedding $M_{\widehat{\C\al}}(1,0)\hr V_{A_1}$, the finite induction and restriction form an adjoint pair of functors $$\left(\mathrm{Ind}_{M_{\widehat{\C\al}}(1,0)}^{V_{A_1}}\dashv \Res^{V_{A_1}}_{M_{\widehat{\C\al}}(1,0)}\right): \adm(M_{\widehat{\C\al}}(1,0))\rightleftarrows \adm(V_{A_1}).$$
		Moreover, they have the following effect on irreducible modules in these categories: 
		\begin{equation}\label{eq:rankoneHinLind}
			\begin{aligned}
				&\Res^{V_{A_1}}_{M_{\widehat{\C\al}}(1,0)}(V_{A_1})\cong M_{\widehat{\C\al}}(1,0),\\
				&\Res^{V_{A_1}}_{M_{\widehat{\C\al}}(1,0)}(V_{A_1+\frac{1}{2}\al})\cong M_{\widehat{\C\al}}(1,\al/2)\op M_{\widehat{\C\al}}(1,-\al/2),\\
				&\mathrm{Ind}_{M_{\widehat{\C\al}}(1,0)}^{V_{A_1}} (M_{\widehat{\C\al}}(1,0))\cong V_{A_1},\\
				&\mathrm{Ind}_{M_{\widehat{\C\al}}(1,0)}^{V_{A_1}} (M_{\widehat{\C\al}}(1,\pm \al/2))\cong V_{A_1+\frac{1}{2}\al},\\
				&\mathrm{Ind}_{M_{\widehat{\C\al}}(1,0)}^{V_{A_1}} (M_{\widehat{\C\al}}(1,\la))=0,\quad \la\in \C\al\bs\{0,\pm\al/2\}. 
			\end{aligned}
		\end{equation}
	\end{prop}	
	\begin{proof}
		Since the Heisenberg VOA $M_{\widehat{\C\al}}(1,0)$ satisfies the strongly unital property for mode transition algebras \cite[Proposition 7.2.1]{DGK23} and the lattice VOA $V_{A_1}$ is rational \cite{D}, then $\mathrm{Ind}_U^V\dashv \Res^V_U$ by Proposition~\ref{prop:adjointpari}. Moreover, any admissible $U$ (or $V$) module $M$ is a generalized Verma module $M=\Phi^\L_U(\Om(M))$. Note that $\Om_{V_{A_1}}(V_{A_1})=\C\vac$ is the trivial module over $A(M_{\widehat{\C\al}}(1,0))\cong \C[x]$,  and $\Om_{V_{A_1}}(V_{A_1+\frac{1}{2}\al})=\C e^{\frac{1}{2}\al}\op \C e^{-\frac{1}{2}\al}$ is a direct sum of irreducible $\C[x]$-modules. Then by \eqref{eq:res}, 
		\begin{align*}
			&\Res^{V_{A_1}}_{M_{\widehat{\C\al}}(1,0)}(V_{A_1})=\Phi^\L_{M_{\widehat{\C\al}}(1,0)}(\C\vac)\cong M_{\widehat{\C\al}}(1,0),\\
			&\Res^{V_{A_1}}_{M_{\widehat{\C\al}}(1,0)}(V_{A_1+\frac{1}{2}\al})=\Phi^\L_{M_{\widehat{\C\al}}(1,0)}(\C e^{\frac{1}{2}\al})\op \Phi^\L_{M_{\widehat{\C\al}}(1,0)}(\C e^{-\frac{1}{2}\al})\cong M_{\widehat{\C\al}}(1,\al/2)\op M_{\widehat{\C\al}}(1,-\al/2).
		\end{align*}
		On the other hand, for any $\la\in \C\al$, recall that $\Om_{M_{\widehat{\C\al}}(1,0)}(M_{\widehat{\C\al}}(1,\la))=\C e^\la$ is a module over $A(M_{\widehat{\C\al}}(1,0))\cong \C[x]$ via $x.e^\la=(\la|\al)e^\la$. Since $\ker(\pi). e^\la=\C ((\la|\al)^3-(\la|\al))\cdot e^\la$, we have the following characterization of the left $A_U=\C[h]$-module in \eqref{eq:ind}: 
		\begin{equation}\label{eq:AVA1mod}
			\frac{\Om_{M_{\widehat{\C\al}}(1,0)}(M_{\widehat{\C\al}}(1,\la))}{\ker(\pi).\Om_{M_{\widehat{\C\al}}(1,0)}(M_{\widehat{\C\al}}(1,\la))}=\frac{\C e^\la}{\C ((\la|\al)^3-(\la|\al))\cdot e^\la}=\begin{cases}
				\C e^\la& \mathrm{if}\ \la=0\ \mathrm{or}\ \pm\al/2,\\
				0&\mathrm{if}\ \la\in \C\al\bs\{0,\pm\al/2\}. 
			\end{cases}
		\end{equation}
		For $\la=0$, the left $\C[h]$-module \eqref{eq:AVA1mod} is $\C \vac$ and we have $e\o \vac=-eh\o \vac=-e\o h.\vac=0$ and $f\o \vac=fh\o\vac=f\o h.\vac=0$  in $A(V_{A_1})\o_{\C[h]}\C\vac$, in view of \eqref{3.26}. It follows from \eqref{eq:ind} that 
		$$\mathrm{Ind}_{M_{\widehat{\C\al}}(1,0)}^{V_{A_1}} (M_{\widehat{\C\al}}(1,0))=\Phi_{V_{A_1}}^\L(A(V_{A_1})\o_{\C[h]}\C\vac)=\Phi_{V_{A_1}}^\L(\C 1\o\vac)\cong V_{A_1}.$$
		Similarly, when the $\C[h]$-module \eqref{eq:AVA1mod} is $\C e^{\pm \frac{1}{2}\al}$, we can show that $A(V_{A_1})\o_{\C[h]}\C e^{\frac{1}{2}\al}=\C 1\o e^{\frac{1}{2}\al}+\C f\o e^{\frac{1}{2}\al}$ and $A(V_{A_1})\o_{\C[h]}\C e^{-\frac{1}{2}\al}=\C 1\o e^{-\frac{1}{2}\al}+\C e\o e^{-\frac{1}{2}\al}$ are both isomorphic to the two-dimensional left $A(V_{A_1})$-module $\Om_{V_{A_1}}(V_{A_1+\frac{1}{2}\al})$. Then 
		$$\mathrm{Ind}_{M_{\widehat{\C\al}}(1,0)}^{V_{A_1}} (M_{\widehat{\C\al}}(1,\pm \al/2))=\Phi_{V_{A_1}}^\L(A(V_{A_1})\o_{\C[h]}\C e^{\pm \frac{1}{2}\al})=\Phi_{V_{A_1}}^\L(\Om_{V_{A_1}}(V_{A_1+\frac{1}{2}\al}))\cong V_{A_1+\frac{1}{2}\al}.$$
		This shows \eqref{eq:rankoneHinLind}. 
	\end{proof}
	
	\subsection{The Borel-type subVOA embedding $V_B\hr V_{A_1}$}  Consider the rank-one Borel-type subVOA of the lattice VOA $V_{A_1}$ \cite{Liu24}:
	$$V_B=\bigoplus_{n=0}^\infty M_{\widehat{\C\al}}(1,n\al)\leq V_{A_1}.$$ 
	The following results were proved in \cite[Section 6]{Liu24}.
	
	\begin{lm}\label{lm:VB}
		$O(V_B)$ \eqref{eq:defofOV} in the definition of $A(V_B)$ is spanned by the following elements:
		\begin{equation}\label{eq:spanovb}
			\begin{cases}
				\al(-n-2)u+\al(-n-1)u,\quad  u\in V_B,\ n\geq 0,\\ 
				\al(-1)v+v,\quad v\in M_{\widehat{\C\al}}(1,m\al),\ m\geq 1,\\
				M_{\widehat{\C\al}}(1,k\al),\quad k\geq 2.
			\end{cases}
		\end{equation}
		Consequently, $A(V_B)\cong \C[x]\op \C y$ as associative algebras, where $y^2=0$, $xy=y$, and $yx=-y$. The identification map is given by $[\al(-1)\vac]\mapsto x$ and $[e^\al]\mapsto y$. Moreover, $O(V_{A_1})$ is spanned by $O(V_B)\cup \{\al(-1)^3w-\al(-1)w: w\in  M_{\widehat{\C\al}}(1,0)\}$.
		
		Finally, the irreducible $V_B$-modules are in one-to-one correspondence with irreducible Heisenberg modules $M_{\widehat{\C\al}}(1,\la)$ on which the positive-part $\bigoplus_{n\geq 1}M_{\widehat{\C\al}}(1,n\al)\leq V_B$ acts as zero.
	\end{lm}

	\begin{lm}\label{lm:VB1}
		For the VOA embedding $V_B\hr V_{A_1}$, the homomorphism $\pi$ \eqref{eq:seqforA} is given by
		\begin{equation}\label{eq:VBpi}
			\pi: \C[x]\op \C y\ra U(\sl_2)/\<e^2\>,\quad \pi (f(x))=f(h),\ \pi(y)=e.
		\end{equation}
		Moreover, $\ker(\pi)=\<x^3-x\>=\spn\{f(x)\cdot (x^3-x):f(x)\in \C[x]\}$. 
	\end{lm}
	\begin{proof}
		The format of $\pi$ follows from the diagram:
		$$
		\begin{tikzcd}
			V_B\arrow[r,hook] \arrow[d]& V_{A_1}\cong L_{\widehat{\sl_2}}(1,0)\arrow[d]\\
			\C[x]\op \C y\arrow[r,"\pi"] & U(\sl_2)/\<e^2\>
		\end{tikzcd}
		\quad 
		\begin{tikzcd}
			\al(-1)\vac \arrow[r,mapsto]\arrow[d,mapsto] & h(-1)\vac\arrow[d,mapsto]\\
			x\arrow[r,mapsto,"\pi"] & h
		\end{tikzcd}
		\quad 
		\begin{tikzcd}
			e^\al \arrow[r,mapsto]\arrow[d,mapsto] & e(-1)\vac\arrow[d,mapsto]\\
			y\arrow[r,mapsto,"\pi"] & e.
		\end{tikzcd}
		$$
		By Lemma~\ref{lm:VB} and \eqref{eq:spanovb}, $\ker(\pi)=(O(V_{A_1})\cap V_B)/O(V_B)$ is spanned by $[\al(-1)^3w-\al(-1)w]=[w]\ast x^3-[w]\ast x$ in $A(V_B)$, where $w\in M_{\widehat{\C\al}}(1,0)$ with $[w]=f(x)\in A(V_B)$.  
	\end{proof}
	
	Note that $V_B$ does not satisfy the strongly unital property for mode transition algebras since it is not a simple VOA. 
	
	By Lemma~\ref{lm:adjointequi}, the adjoint pair $\Phi_{V_B}^\L:\mathsf{Mod}(A(V_B))\rightleftarrows\mathsf{Adm}(V_B):\Om_{V_B}$ is not an adjoint equivalence. Although this is not quite satisfactory, we can still determine the finite induced modules of the irreducible $V_{B}$-modules using the previous Lemma. 
	
	\begin{prop}
		For the VOA embedding $V_B\hr V_{A_1}$, the finite induction functor has the following effect on irreducible $V_B$-modules:  
		\begin{equation}\label{eq:indforV_B}
			\begin{aligned}
				& \mathrm{Ind}_{V_B}^{V_{A_1}} (M_{\widehat{\C\al}}(1,0)) \cong V_{A_1},\\
				& \mathrm{Ind}_{V_B}^{V_{A_1}}(M_{\widehat{\C\al}}(1,\al/2))\cong V_{A_1+\frac{1}{2}\al},\\
				& \mathrm{Ind}_{V_B}^{V_{A_1}}(M_{\widehat{\C\al}}(1,\la))=0,\quad \la\in \C\al\bs\{0,\al/2\}.
			\end{aligned}
		\end{equation}
	\end{prop}
	\begin{proof}
		By \eqref{eq:VBpi}, $A_{V_B}=\pi(A(V_B))$ is the subalgebra of $U(\sl_2)/\<e^2\>$ generated by $h$ and $e$. Given irreducible $V_B$-module $M=M_{\hat{\h}}(1,\la)$. Since $ V_B=M_{\widehat{\C\al}}(1,0)\op V_+$, with $V_+$ acts as zero on $M$, we have $y\in A(V_B)$ acts as zero on $\Om_{V_B}(M)=\C e^\la$. Similar to  \eqref{eq:AVA1mod}, by Lemma~\ref{lm:VB1}, we have 
		$$
		\frac{\Om_{V_B}(M_{\widehat{\C\al}}(1,\la))}{\ker(\pi). \Om_{V_B}(M_{\widehat{\C\al}}(1,\la))}=
		\frac{\C e^\la}{\ker(\pi). e^\la}
		=\begin{cases}
			\C e^\la& \mathrm{if}\ \la=0\ \mathrm{or}\ \pm\al/2,\\
			0&\mathrm{if}\ \la\in \C\al\bs\{0,\pm\al/2\}. 
		\end{cases} 
		$$
		However, when $\la=-\al/2$, since $h^2-h-2fe=0$ in $U(\sl_2)/\<e^2\>$, we have the following relation in $A(V_{A_1})\o_{A_{V_B}} \C e^{-\frac{\al}{2}}$ due to the fact that $\pi(y)=e$ acts as zero on $ \C e^{-\frac{\al}{2}}$: 
		$$0=(h^2-h-2fe)\o e^{-\frac{\al}{2}}=\left((\al|-\al/2)^2-(\al|-\al/2)\right)\cdot 1\o e^{-\frac{\al}{2}}-2f\o e.e^{-\frac{\al}{2}}=2\cdot 1\o e^{-\frac{\al}{2}}. $$
		Hence $A(V_{A_1})\o_{A_{V_B}} \C e^{-\frac{\al}{2}}=0$ and $\Ind^{V_{A_1}}_{V_B}(M_{\widehat{\C\al}}(1,-\al/2))=0$. Now \eqref{eq:indforV_B} follows immediately from this observation and the rationality of the lattice VOA $V_{A_1}$. 
	\end{proof}
	
	\begin{remark}
		Since the subVOA $V_B$ does not satisfy the strongly unital condition for mode transition algebras , the irreducible $V_B$-modules $W$ cannot be completely determined by their bottom degree $\Om_{V_B}(W)$. Therefore, we omit the discussion for finite restrictions in this case. 
	\end{remark}

	\subsection{The Virasoro embedding $L(1,0)\hr V_{A_1}$}
	The Virasoro element $\om=\frac{1}{4}\al(-1)^2\vac\in V_{A_1}$ generates the Virasoro subVOA $L(1,0)$ \cite{DG98}. Any irreducible module $L(1,k)$ over the Virasoro algebra of central charge $1$ is an irreducible module over the VOA $L(1,0)$ \cite{FZ92}. 
	
	Recall that the Zhu algebra 
	$A(L(1,0))\cong \C[y]$, with $[\om]\mapsto y$ \cite{DMZ94,W93}. Moreover, we have the following characterization of the homomorphism $\pi$ \eqref{eq:seqforA}: 
	\begin{equation}\label{eq:vira}
		\begin{tikzcd}
			L(1,0)\arrow[r,hook]\arrow[d]& V_{A_1}\arrow[d]\\
			\C[y]\arrow[r,"\pi"] & U(\sl_2)/\<e^2\>
		\end{tikzcd}\quad 
		\begin{tikzcd}
			\om=\frac{1}{4}\al(-1)^2\vac \arrow[r,mapsto]\arrow[d,mapsto] & \frac{1}{4}h(-1)^2\vac\arrow[d,mapsto]\\
			y\arrow[r,mapsto,"\pi"] & \frac{1}{4}h^2.
		\end{tikzcd}
	\end{equation}
	Since $h^3-h=0$, we have $(\frac{1}{4}h^2)^2=\frac{1}{16}h^4=\frac{1}{16}h^2$. Hence $\ker(\pi)=\<y^2-\frac{1}{4}y\>$ by \eqref{eq:vira}, and
	\begin{equation}\label{eq:AL10}
		A_{L(1,0)}=\pi(\C[y])\cong \C[y]/\<y(y-(1/4))\>= \C[h^2]\leq A(V_{A_1})
	\end{equation}
	
	Given $k\in \C$, we have $\Om_{L(1,0)}(L(1,k))=\C v_{1,k}$ as an $A(L(1,0))\cong \C[y]$-module, with the action $y.v_{1,k}=k\cdot v_{1,k}$. Then $(y^2-\frac{1}{4}y).v_{1,k}=(k^2-\frac{1}{4}k)v_{1,k}=0$ if and only if $k=0$ or $\frac{1}{4}$. It follows that 
	\begin{equation}\label{eq:3.8}
		\frac{\Om_{L(1,0)}(L(1,k))}{\ker(\pi).\Om_{L(1,0)}(L(1,k))}=\frac{\C v_{1,k}}{\ker(\pi).v_{1,k}}=\begin{cases}
			\C v_{1,k}&\mathrm{if}\ k=0\ \mathrm{or}\ \frac{1}{4},\\
			0&\mathrm{if}\ k\in \C\bs\{0,\frac{1}{4} \}. 
		\end{cases}
	\end{equation}
	
	\begin{prop}\label{prop:rankonevir}
		For the VOA embedding $L(1,0)\hr V_{A_1}$, the finite induction functor has the following effect on irreducible $L(1,0)$-modules:
		\begin{align*}
			&\Ind^{V_{A_1}}_{L(1,0)}(L(1,0))\cong V_{A_1},\\
			&\Ind^{V_{A_1}}_{L(1,0)}(L(1,1/4))\cong V_{A_1+\frac{1}{2}\al}\op V_{A_1+\frac{1}{2}\al},\numberthis\label{eq:virind}\\
			&\Ind^{V_{A_1}}_{L(1,0)}(L(1,k))=0,\quad k\in \C\bs\{0,1/4\}. 
		\end{align*}
	\end{prop}
	\begin{proof}
		Clearly, if $k\neq 0,\frac{1}{4}$, we have  $\Ind^{V_{A_1}}_{L(1,0)}(L(1,k))=\Phi^\L_{V_{A_1}}(A(V_{A_1})\o_{A_{L(1,0)}}0)=0$ \eqref{eq:3.8}. Let $k=0$. Since $h^2.v_{1,0}=0$, we have the following identification of left $A(V_{A_1})$-modules: 
		$$A(V_{A_1})\o_{\C[h^2]} \C v_{1,0}\cong A(V_{A_1})/A(V_{A_1})\cdot  h^2=A(V_{A_1})/I,$$
		where $I=A(V_{A_1})\cdot  h^2$ is the left ideal of $A(V_{A_1})$ generated by $h^2$. By \eqref{3.26}, we have $e=eh^2\in I$; which implies $h=h^2-2fe\in I$; which further implies $f=fh\in I$. Hence $I=\spn\{e,f,h,h^2\}\subset A(V_{A_1})$ and $A(V_{A_1})/I\cong \C \vac$. We have 
		$$
		\Ind^{V_{A_1}}_{L(1,0)}(L(1,0))=\Phi^\L_{V_{A_1}}(A(V_{A_1})\o_{\C[h^2]} \C v_{1,0})\cong \Phi^\L_{V_{A_1}}(\C\vac)\cong V_{A_1}.
		$$
		Finally, let $k=\frac{1}{4}$. Since $\frac{1}{4}h^2.v_{1,\frac{1}{4}}=\pi (y).v_{1,\frac{1}{4}}=\frac{1}{4}v_{1,\frac{1}{4}}$ we have an identification:
		$$A(V_{A_1})\o_{\C[h^2]} \C v_{1,\frac{1}{4}}\cong A(V_{A_1})/A(V_{A_1})\cdot (h^2-1)=A(V_{A_1})/J,$$
		where $J=A(V_{A_1})\cdot (h^2-1)$. By \eqref{3.26} again, we have $e(h^2-1)=(-1)^2e-e=0$, $f(h^2-1)=0$, and $h(h^2-1)=h^3-h=0$. Since $A(V_{A_1})$ is generated by $e,h,f$ as an associative algebra,  we have $J=A(V_{A_1})\cdot (h^2-1)=\spn\{ h^2-1\}$ and $A(V_{A_1})/J=\spn\{\bar{1},\bar{e},\bar{f},\bar{h}\}.$ Consider the decomposition 
		$$A(V_{A_1})/J=M\op N,\quad M=\spn\{\bar{e},\bar{1}-\bar{h}\},\quad N=\spn\{\bar{f},\bar{1}+\bar{h} \}.$$
		Using \eqref{3.26}, together with $\bar{h}^2=\bar{1}$, it is easy to show that $M$ and $N$ are both isomorphic to the irreducible $\sl_2$-modules $L(1)$. Hence 
		$$
		\Ind^{V_{A_1}}_{L(1,0)}(L(1,1/4)=\Phi^\L_{V_{A_1}}(A(V_{A_1})\o_{\C[h^2]} \C v_{1,\frac{1}{4}})\cong \Phi^\L_{V_{A_1}}(M\op N)\cong V_{A_1+\frac{1}{2}\al}\op V_{A_1+\frac{1}{2}\al}.
		$$
		This proves \eqref{eq:virind}. 
	\end{proof}
	
	\begin{remark}
		Again, since the subVOA $L(1,0)$ does not satisfy the strongly unital condition for mode transition algebras \cite[Proposition 8.1.1]{DGK23}, the irreducible $L(1,0)$-modules $L(1,h)$ are not completely determined by their bottom degree $\Om_{L(1,0)}(L(1,h))=\C v_{1,h}$. i.e., $L(1,h)$ is not isomorphic to $ \Phi^\L_{L(1,0)}(\Om_{L(1,0)}(L(1,h)))$ in general. 
	\end{remark}

	
	\section{Finite induction for the rational VOA embedding $V_{A_1}\hr V_{A_2}$}\label{sec:4}

	The $A_2$-root system $\Phi_{A_2}=\{\pm \al, \pm \b,\pm(\al+\b)\}$, with standard basis $\{\al,\b\}$ such that 
	\begin{equation}\label{eq:paringinA2lattice}
		(\al|\al)=2,\quad (\b|\b)=2,\quad (\al|\b)=-1.
	\end{equation}
	The longest root $\theta=\al+\b$, see \cite{Hum2}. 
	
	The root lattice $A_1=\Z\al$, with $(\al|\al)=2$, can be naturally embedded into the root lattice $A_2=\Z\al \op \Z\b$. The induced embedding $\C\al\hr \C \al\op \C\b=\h$, leads to an embedding of Heisenberg algebras $$\widehat{\C\al}=\C\al \o \C[t,t^{-1}]\op \C K \hr \hat{\h}=\h\o \C[t,t^{-1}]\op \C K ,\quad \al(n)\mapsto \al(n),\ K\mapsto K,$$
	which further leads to an embedding of Verma modules $M_{\widehat{\C\al}}(1,n\al)\hr M_{\hat{\h}}(1,n\al)$ for any $n\in \Z$. Moreover, by the definition of lattice vertex operators \cite{FLM}, this leads to a {\bf non-conformal} embedding of lattice VOAs
	\begin{equation}\label{eq:A1A2emb}
		V_{A_1}=\bigoplus_{n\in \Z} M_{\widehat{\C\al}}(1,n\al)\hr \bigoplus_{m,n\in \Z}M_{\hat{\h}}(1,n\al+m\b)=V_{A_2}. 
	\end{equation}
	Since $V_{A_1}\cong L_{\widehat{\sl_2}}(1,0)$ and $V_{A_2}\cong L_{\widehat{\sl_3}}(1,0)$ as VOAs, the lattice VOA embedding \eqref{eq:A1A2emb} is equivalent to the affine VOA embedding:
	\begin{equation}\label{eq:affineA1toA2}
		\begin{aligned}
			L_{\widehat{\sl_2}}(1,0)&\hr  L_{\widehat{\sl_3}}(1,0),\\
			e(-1)\vac\mapsto x_\al(-1)\vac,\quad  f(-1)\vac&\mapsto  x_{-\al}(-1)\vac,\quad h(-1)\vac\mapsto x(-1)\vac,
		\end{aligned}
	\end{equation}
	where $\{x_{\pm \al},x\}\subset \sl_3$ forms a Lie subalgebra $\sl_2$. See Figure~\ref{fig1} for an illustration. 
	
	\begin{figure}
		\centering
		\begin{tikzpicture}
			\coordinate (Origin)   at (0,0);
			\coordinate (XAxisMin) at (-5,0);
			\coordinate (XAxisMax) at (5,0);
			\coordinate (YAxisMin) at (0,0);
			\coordinate (YAxisMax) at (0,5);
			
			\clip (-5.2,-3.8) rectangle (5.2,3.7); 
			\begin{scope} 
				\pgftransformcm{1}{0}{1/2}{sqrt(3)/2}{\pgfpoint{0cm}{0cm}} 
				\coordinate (Bone) at (0,2);
				\coordinate (Btwo) at (2,-2);
				\draw[style=help lines,dashed] (-7,-6) grid[step=2cm] (6,6);
				\foreach \x in {-4,-3,...,4}{
					\foreach \y in {1,...,4}{
						\coordinate (Dot\x\y) at (2*\x,2*\y);
						\node[draw,circle,inner sep=2pt,fill] at (Dot\x\y) {};
					}
				}
				\foreach \x in {-4,-3,...,4}{
					\foreach \y in {0}{
						\coordinate (Dot\x\y) at (2*\x,2*\y);
						\node[draw,red,circle,inner sep=2pt,fill] at (Dot\x\y) {};
					}
				}
				
				\foreach \x in {-4,-3,...,4}{
					\foreach \y in {-1,...,-4}{
						\coordinate (Dot\x\y) at (2*\x,2*\y);
						\node[draw,circle,inner sep=2pt,fill] at (Dot\x\y) {};
					}
				}
				\foreach \x in {0}{
					\foreach \y in {0}{
						\coordinate (Dot\x\y) at (2*\x,2*\y);
						\node[draw,circle,red,inner sep=2pt,fill] at (Dot\x\y) {};
						\node [xshift=0cm, yshift=-0.4cm]{$0$};
					}
				}
				
				\draw [thick,-latex] (Origin) 
				-- (Bone) node [above right]  {$\al+\b$};
				\draw [thick,-latex] (Origin)
				-- (Btwo) node [below right]  {$-\b$};
				\draw [thick,-latex,red] (Origin) node [xshift=0cm, yshift=-0.3cm]{}
				-- ($(Bone)+(Btwo)$) node [below right] {$\al$};
				\draw [thick,red,-latex] (Origin)
				-- ($-1*(Bone)-1*(Btwo)$) node [below left] {$-\al$};
				\draw [thick,-latex] (Origin)
				-- ($-1*(Btwo)$) coordinate (B3) node [above left] {$\b$};
				\draw [thick,-latex] (Origin)
				-- ($-1*(Bone)$) node [below left]  {$-\al-\b$};
			\end{scope} 
			\begin{scope}
				\pgftransformcm{1}{0}{-1/2}{sqrt(3)/2}{\pgfpoint{0cm}{0cm}} 
				\draw[style=help lines,dashed] (-6,-6) grid[step=2cm] (6,6);
			\end{scope}
		\end{tikzpicture}
		\caption{\label{fig1}}
	\end{figure}

	\subsection{Relations in the associative algebra $A(V_{A_2})$}
	Choose a standard basis of $\sl_3$:
	$$\sl_3=\spn\{x_{\pm \al}, x_{\pm \b}, x_{\pm(\al+\b)}, x,y\},$$
	where $x=h_\al$ and $y=h_\b$. The spanning elements satisfy the Serre relations:
	\begin{align*}
		& [x,x_{\pm \al}]=\pm 2x_{\pm \al},&& [x,x_{\pm \b}]=\mp x_{\pm \b},&& [y,x_{\pm \al}]=\mp x_{\pm \al}, && [y,x_{\pm \b}]=\pm 2x _{\pm \b},\\
		& [x,y]=0,&& [x_\al,x_{-\al}]=x,&& [x_\b,x_{-\b}]=y,&& [x_\al,x_\b]=x_{\al+\b},\numberthis\label{eq:serrerel}\\
		& [x_{-\al},x_{-\b}]=x_{-(\al+\b)},&& (\ad x_\al)^2(x_\b)=0,&& (\ad x_{-\al})^2(x_\b)=0,&& (\ad x_\al)(x_{-\b})=0.
	\end{align*}
	
	
	By the PBW-theorem, we have an embedding of universal enveloping algebra
	\begin{equation}\label{eq:iota}
		\iota: U(\sl_2)\hr U(\sl_3),\quad e^s h^s f^t\mapsto x_{\al}^rx^s x_{-\al}^t,\ r,s,t\in \N.
	\end{equation}
	
	\begin{prop}\label{prop:embeddingofA1inA2}
		The embedding $\iota$ \eqref{eq:iota} induces an injective homomorphism $\pi$ \eqref{eq:seqforA}  between Zhu algebras 
		\begin{equation}\label{homoA1A2}
			\begin{aligned}
				\pi: A(V_{A_1})\cong U(\sl_2)/\<e^2\>&\hr U(\sl_3)/\<x_{\al+\b}^2\>=A(V_{A_2}),\\
				e\mapsto x_\al,\quad f&\mapsto x_{-\al},\quad h\mapsto h_\al,
			\end{aligned}
		\end{equation}
		where we use the same notation for the image of elements in the quotients. 
	\end{prop}
	\begin{proof}
		Denote the two-sided ideal $\<x_{\al+\b}^2\>\lhd U(\sl_3)$ by $J$. We claim that $x_\al^2\in \iota(U(\sl_2))\cap J$. Indeed, the bracket relation among standard basis elements of $\sl_3$ satisfy $[x_{-\b},x_{\al+\b}]=x_\al$ and $[x_{-\b},x_{\al}]=0$. Then 
		\begin{equation}\label{eq:xal}
			[x_{-\b},x_{\al+\b}^2]=[x_{-\b},x_{\al+\b}]x_{\al+\b}+x_{\al+\b}[x_{-\b},x_{\al+\b}]=x_\al x_{\al+\b}+x_{\al+\b}x_\al\in J,
		\end{equation}
		and so $x_\al^2=\frac{1}{2}([x_{-\b}, x_\al x_{\al+\b}+x_{\al+\b}x_\al])\in  J$. In particular, $\iota(\<e^2\>)\ssq \iota(U(\sl_2))\cap J$ and $\pi$ \eqref{homoA1A2} is  well-defined. Note that $\pi$ agrees with the affine VOA embedding $L_{\widehat{\sl_2}}(1,0)\hr L_{\widehat{\sl_3}}(1,0)$ \eqref{eq:affineA1toA2}. 
		
		Since $x^2_\al\in J$, by applying the same Lie brackets that give rise to relations \eqref{3.26}, we can show that $x_\al x+ x_\al\in J$, $x_{-\al} x-x_{-\al}\in J$, and $x^3-x\in J$. Moreover, a similar argument as \eqref{eq:xal} shows $x_\b^2\in J$, and so 
		$$x_\b y+ x_\b\in J,\quad x_{-\b} y-x_{-\b}\in J,\quad \mathrm{and} \quad y^3-y\in J.$$
		
		To show $\pi$ is injective, we note that $\ker(\pi)$ is a two-sided ideal of the semisimple algebra $U(\sl_2)/\<e^2\>\cong M_{1}(\C)\times M_{2}(\C)$ with primitive central idempotents $h^2$ and $1-h^2$. Then the ideal $\ker(\pi)$ is either generated by $h^2$ or $1-h^2$. 
		
		If $h^2\in \ker(\pi)$, then $x^2=\iota(h^2)\in J\lhd U(\sl_3)$. It follows that $x^3\in J$ and $x=x^3-(x^3-x)\in J$. But then $\sl_3\subset J$, since $\sl_3$ is a simple Lie algebra. This indicates $U(\sl_3)=J$, which is a contradiction. On the other hand, if $1-h^2\in \ker(\pi)$, then $1-x^2\in J$. Since $[x,x_\b]=(\al|\b)x_\b=-x_\b$, we have $[x_\b,1-x^2]=x_\b x+xx_\b\in J$ and 
		$xx_\b-x_\b x+x_\b=0\in J$. Thus, 
		$$2xx_\b+x_\b=(x_\b x+xx_\b)+(xx_\b-x_\b x+x_\b)\in J.$$
		Since we also have $yx_\b -x_\b\in J$, it follows that $(2x+y)x_\b\in J$. Note that $[2x+y,x_{-\b}]=(-\b| 2\al+\b) x_{-\b}=0$, we have 
		$$4xy+2y^2=(2x+y)\cdot 2y=[(2x+y)x_\b,x_{-\b}]\in J.$$
		Since $y^3-y\in J$ and $1-x^2\in J$, we also have 
		$$2xy+y^2\in J\implies 2xy^2+y\in J\implies 2x^2 y^2+xy\in J\implies 2y^2+xy\in J. $$
		In particular, we have $2xy\in J$ and $y^2\in J$. Similar to the case when $x^2\in J$, we can derive that $\sl_3\ssq J$, which is a contradiction. Therefore, $\pi$ \eqref{homoA1A2} is injective. 
	\end{proof}
	
	The relations among generators of the algebra $U(\sl_3)/\<x_{\al+\b}^2\>$ are used multiple times in the proof of the Proposition. More generally, an algebraic calculation gives rise to all the relations in this algebra, generalizing the rank-one relations \eqref{3.26}. 
	
	\begin{prop}\label{prop:presentationofAVA2}
		The associative algebra $A(V_{A_2})=U(\sl_3)/\<x_{\al+\b}^2\>$ can be presented by generators $\{x_{\pm \al}, x_{\pm \b}, x_{\pm(\al+\b)}, x,y\}$, together with relations: 
		\begin{align*}
			&\begin{cases}
				xx_{\pm \al}=\pm x_{\pm \al},&x_{\pm \al}x=\mp x_{\pm \al},\\
				yx_{\pm \b}=\pm x_{\pm\b},& x_{\pm \b} y=\mp x_{\pm \b},\\
				(x+y)x_{\pm (\al+\b)}=\pm x_{\pm (\al+\b)},& x_{\pm (\al+\b)}(x+y)=\mp x_{\pm( \al+\b)},
			\end{cases}\numberthis\label{rel1}\\
			&\begin{cases}
				yx_{\pm \al}=x_{\pm \al}y\mp x_{\pm \al},& xx_{\pm \b}=x_{\pm \b} x\mp x_{\pm \b},\\
				xx_{\pm (\al+\b)}=x_{\pm (\al+\b)}x\pm x_{\pm(\al+\b)}
				& yx_{\pm (\al+\b)}=x_{\pm (\al+\b)}y\pm x_{\pm (\al+\b)},
			\end{cases}\numberthis\label{rel2}\\
			&\begin{cases}
				x_\al x_{-\al}=\frac{1}{2}x^2+\frac{1}{2}x,\quad & x_{\al+\b}x_{-\al-\b}=\frac{1}{2} (x+y)^2+\frac{1}{2}(x+y)\\
				x_\b x_{-\b}=\frac{1}{2}y^2+\frac{1}{2}y,& 	x_\gamma x_\theta=0\ \mathrm{if}\ \ga+\theta\notin \Phi_{A_2}=\{ \pm \al, \pm\b,\pm (\al+\b)\}
			\end{cases}\numberthis\label{rel3} \\
			&\begin{cases}
				x^3-x=0,&  (x+y)^3-(x+y)=0,\\
				y^3-y=0,& xy=yx.
			\end{cases}\numberthis\label{rel4}
		\end{align*}
		Moreover, by applying proper Lie brackets to the existing relations, the following additional relations in $A(V_{A_2})$ can be derived from \eqref{rel1}--\eqref{rel4}: 
		\begin{align*}
			&\begin{cases}
				x_\al x_\b=-x_{\al+\b}y,\ & x_\b x_\al=-x_{\al+\b}y-x_{\al+\b},\\
				x_{-\al}x_{-\b}=-x_{-(\al+\b)}x+x_{-(\al+\b)},&  x_{-\b}x_{-\al}=-x_{-(\al+\b)}x. 
			\end{cases}\numberthis\label{eq:morerel1} \\
			&\begin{cases}
				x_\b x_{-(\al+\b)}=-x_{-\al}y-x_{-\al}, \quad \ \ \ &	x_{-(\al+\b)}x_\b=-x_{-\al}y,\\
				x_{-\b}x_{\al+\b}=x_{\al}y-x_\al,&	x_{\al+\b}x_{-\b}=x_{\al}y.
			\end{cases}\numberthis\label{eq:morerel2}\\
			&\begin{cases}
				x_\al x_{-(\al+\b)}=x_{-\b}x+x_{-\b},\qquad \ \ & x_{-(\al+\b)}x_\al=x_{-\b}x,\\
				x_{-\al}x_{\al+\b}=-x_{\b}x+x_\b,& x_{\al+\b}x_{-\al}=-x_{\b}x. 
			\end{cases}\numberthis\label{eq:morerel3}
		\end{align*}
	\end{prop}
	Since the Proposition involves calculation only, we omit the details of the proof. 
	Note that some relations can be derived from the others. For instance, $ x_\b x_\al=-x_{\al+\b}y-x_{\al+\b}$ follows  from $[x_\al,x_\b]=x_{\al+\b}$, $ x_\al x_\b=-x_{\al+\b}y$, and $ x_{\al+\b}(x+y)=-x_{\al+\b}$. In other words,  \eqref{rel1}--\eqref{rel4} are not the minimal sets of the relations. We list them as such for the convenience of our later argument. 
	
	\begin{coro}\label{coro:basisofAVA2}
		$A(V_{A_2})$ is a $19$-dimensional associative algebra spanned by 
		\begin{equation}\label{eq:spanningeltsofAVA2}
			\{1,\; x_{\pm \alpha},\; x_{\pm \beta},\; x_{\pm (\alpha+\beta)},\; x,\; x^2,\; y,\; y^2,\; yx,\;  y^2x,\; yx_{\pm \alpha},\; x_{\pm \beta}x,\; x_{\pm (\alpha+\beta)}x\}.
		\end{equation}
		Moreover, the elements \eqref{eq:spanningeltsofAVA2} are linearly independent. 
	\end{coro}

	\subsection{Irreducible modules under finite induction and restriction for $V_{A_1}\hr V_{A_2}$}

	Note that both $V_{A_1}\cong L_{\widehat{\sl_2}}(1,0)$ and $V_{A_2}\cong L_{\widehat{\sl_3}}(1,0)$ are rational VOAs \cite{D,DLM97}. Then the vertical functors in the diagram below are equivalence of categories. Therefore, the finite induction and restriction forms an adjoint pair. 
	They are completely determined by the adjoint pair of functors on the Zhu algebra $\A$-level. 
	\begin{equation}\label{eq:commdiagramA1A2}
		\begin{tikzcd}[scale=1.5,row sep=huge, column sep = huge]
			\mathsf{Adm}(V_{A_1})\arrow[r,dashed,shift left=2pt,"\mathrm{Ind}^{V_{A_2}}_{V_{A_1}}"]\arrow[d,shift left=2pt,"\Om_{V_{A_1}}"]& \mathsf{Adm}(V_{A_2})\arrow[l,dashed,shift left=2pt,"\mathrm{Res}^{V_{A_2}}_{V_{A_1}}"]\arrow[d,shift left=2pt,"\Om_{V_{A_2}}"]\\
			\mathsf{Mod}(A(V_{A_1}))\arrow[r,shift left=3pt,"\mathrm{Ind}^{A(V_{A_2})}_{A(V_{A_1})}"]\arrow[u,shift left=2pt,"\Phi_{V_{A_1}}^\L"]& \mathsf{Mod}(A(V_{A_2})) \arrow[l,shift left=2pt,"\mathrm{Res}^{A(V_{A_2})}_{A(V_{A_1})}"]\arrow[u,shift left=2pt,"\Phi ^\L _{V_{A_2}}"]
		\end{tikzcd}
	\end{equation}
	Let $\la_1=\frac{2}{3}\al+\frac{1}{3}\b$ and $\la_2=\frac{1}{3}\al+\frac{2}{3}\b$ be the fundamental dominant weights for $A_2$. Then the weight lattice $A_2^\circ=\Z\la_1\op \Z \la_2$, and $A_2^\circ/A_2=A_2\sqcup (A_2+\la_1)\sqcup (A_2+\la_2)$ as right cosets. By \cite[Theorem 3.1]{D}, the VOA $V_{A_2}$ has three irreducible modules $V_{A_2},$ $V_{A_2+\la_1},$ and  $V_{A_2+\la_2},$ with bottom level
	\begin{align*}
		\Om_{V_{A_2}}(V_{A_2})&=\C \vac\cong L(0),\\
		\Om_{V_{A_2}}(V_{A_2+\la_1})&=\C e^{\la_1}+\C e^{\la_1-\al}+\C e^{\la_1-\al-\b}\cong L(\la_1),\\ 
		\Om_{V_{A_2}}(V_{A_2+\la_2})&=\C e^{\la_2}+\C e^{\la_2-\b}+\C e^{\la_2-\al-\b}\cong L(\la_2),
	\end{align*}
	where $L(\la)$ is the irreducible highest-weight module for $\sl_3$ of highest weight $\la\in \h^\ast$. 
	
	
	\begin{thm}\label{thm:indA1A2}
		For the VOA embedding $V_{A_1}\hr V_{A_2}$, the finite induction and restriction functors have the following effects on irreducible modules: 
		\begin{align}
			&	\Ind^{V_{A_2}}_{V_{A_1}}(V_{A_1})\cong V_{A_2}\op V_{A_2+\la_1}\op V_{A_2+\la_2},\label{eq:A1A2ind1}\\
			&	\Ind^{V_{A_2}}_{V_{A_1}}(V_{A_1+\frac{1}{2}\al})\cong V_{A_2+\la_1}\op V_{A_2+\la_2}, \label{eq:A1A2ind2}\\
			&	\Res^{V_{A_2}}_{V_{A_1}}(V_{A_2})\cong V_{A_1},\label{eq:A1A2res1}\\
			&	\Res^{V_{A_2}}_{V_{A_1}}(V_{A_2+\la_i})\cong V_{A_1}\op V_{A_1+\frac{1}{2}\al},\quad i=1,2. \label{eq:A1A2res2}
		\end{align}
	\end{thm}
	\begin{proof}
		Since $\pi: A(V_{A_1})\ra A(V_{A_2})$ is injective \eqref{homoA1A2}, then by \eqref{eq:commdiagramA1A2} and \eqref{eq:ind}, we have 
		$$
		\Ind^{V_{A_2}}_{V_{A_1}}(V_{A_1})=\Phi^\L_{V_{A_2}}\left(A(V_{A_2})\o _{A(V_{A_1})} \Om_{V_{A_1}}(V_{A_1})\right)=\Phi^\L_{V_{A_2}}\left(A(V_{A_2})\o _{A(V_{A_1})}\C \vac\right). 
		$$
		By Corollary~\ref{coro:basisofAVA2}, the left $A(V_{A_2})$-module $A(V_{A_2})\o _{A(V_{A_1})}\C \vac$ is spanned by the following elements: 
		\begin{equation}\label{eq:spanning1}
			\{ x_{\pm (\al+\b)}\o \vac,\; x_{\pm \b}\o \vac,\; y^2\o \vac,\; y\o \vac,\; 1\o \vac \}.
		\end{equation}
		To show these elements are nonzero, we observe that $A(V_{A_2})\o _{A(V_{A_1})}\C \vac\cong A(V_{A_2})/A(V_{A_2}).\{x_{\pm \al} ,x\}$ as a left $A(V_{A_2})$-module, since $\C \vac$ is the trivial $A(V_{A_1})$-module. Using the relations \eqref{rel1}--\eqref{eq:morerel3}, it is easy to check that 
		$$
		A(V_{A_2}).\{ x_{\pm \al},x\}=\spn\{ x_{\pm \alpha},\; x,\; x^2,\; yx,\;  y^2x,\; yx_{\pm \alpha},\; x_{\pm \beta}x,\; x_{\pm (\alpha+\beta)}x\}.
		$$
		Since the elements \eqref{eq:spanningeltsofAVA2} form a basis of $A(V_{A_2})$, it follows that $\{\overline{x_{\pm (\al+\b)}},\; \overline{x_{\pm \al}},\; \bar{y}^2,\; \bar{y},\; \bar{1} \}$ form a basis of $A(V_{A_2})/A(V_{A_2}).\{x_{\pm \al} ,x\}$. Hence \eqref{eq:spanning1} forms a basis of $A(V_{A_2})\o _{A(V_{A_1})}\C \vac$. 
		
		Using the relations \eqref{rel1}--\eqref{eq:morerel3}, it is easy to see that 
		\begin{align*}
			&L(\la_1)=\spn\{x_{\al+\b}\o \vac,\ x_\b\o \vac,\ (y^2-y)\o \vac \},\\
			& L(\la_2)=\spn\{x_{-(\al+\b)}\o \vac,\ x_{-\b}\o \vac,\ (y^2+y)\o\vac  \},\\
			& L(0)=\spn\{(1-y^2)\o \vac \}.
		\end{align*}
		Therefore, $A(V_{A_2})\o _{A(V_{A_1})}\C \vac\cong L(0)\op L(\la_1)\op L(\la_2)$. This proves \eqref{eq:A1A2ind1}.
		
		Denote $\Om_{V_{A_1}}(V_{A_1+\frac{1}{2}\al})=L(\frac{1}{2}\al)=\C e^+\op \C e^-$ as a $\sl_2$-module, with the $A(V_{A_1})\cong U(\sl_2)/\<x_\al^2\>$-module action  given by 
		\begin{align*}
			& x_\al. e^+=0,&& x_\al. e^-=e^+,&& x. e^\pm=\pm e^\pm,\\
			& x_{-\al}.e^+=e^-,&& x_{-\al}.e^-=0.&&
		\end{align*}
		The the left $A(V_{A_2})$-module $A(V_{A_2})\o _{A(V_{A_1})} \Om_{V_{A_1}}(V_{A_1+\frac{1}{2}\al})$ is spanned by 
		\begin{equation}\label{eq:spanning2}
			\{x_{\pm(\al+\b)}\o e^\pm,\ x_{\pm \b}\o e^\pm,\ y^2\o e^\pm,\ y\o e^\pm,\   1\o e^\pm \}.
		\end{equation}
		Using  \eqref{rel1}--\eqref{eq:morerel3}, we can derive the following relations among the spanning elements:
		\begin{align*}
			x_{\al+\b}\o e^+&=x_{\al+\b}\o x_\al.e^-=x_{\al+\b}x_\al\o e^-=0\o e^-=0,\\
			x_{-(\al+\b)}\o e^-&= x_{-(\al+\b)}\o x_{-\al}.e^+=x_{-(\al+\b)}x_{-\al}\o e^+=0\o e^+=0,\\
			x_{\al+\b}\o e^-&=x_{\al+\b}\o x_{-\al}.e^+=x_{\al+\b} x_{-\al}\o e^+= -x_\b x\o e^+=-x_\b\o e^+,\\
			x_{-(\al+\b)}\o e^+&= x_{-(\al+\b)}\o x_\al. e^-=x_{-(\al+\b)}x_\al\o e^-=x_{-\b}x\o e^-=x_{-\b}\o e^-,\\
			x_{\b}\o e^-&=x_\b\o x_{-\al}.e^+=x_\b x_{-\al}\o e^+=0\o e^+=0,\\
			x_{-\b}\o e^+&= x_{-\b}\o x_{\al}.e^-=x_{-\b}x_\al\o e^-=0\o e^-=0,\\
			(y^2+y)\o e^+&=2x_\b x_{-\b}\o e^+=2x_\b. (x_{-\b}\o e^+)=0,\\
			(y^2-y)\o e^-&= 2x_{-\b}x_\b\o e^-=2x_{-\b}\o (x_\b\o e^-)=0.
		\end{align*}
		There are no additional relations among the spanning elements \eqref{eq:spanning2}. Then
		$$A(V_{A_2})\o _{A(V_{A_1})} \Om_{V_{A_1}}(V_{A_1+\frac{1}{2}\al})=\spn\{1\o e^\pm, x_{\al+\b}\o e^-, x_{-(\al+\b)}\o e^+, y\o e^\pm \}.$$
		It is easy to check that  $(y+1)\o e^+$ is the highest-weight vector for $\sl_3$ of weight $\la_1$, and $x_{\al+\b}\o e^-$ is the highest-weight vector for $\sl_3$ of weight $\la_2$. Moreover, we have 
		\begin{align*}
			A(V_{A_2}).((y+1)\o e^+)&=\spn\{(y+1)\o e^+,\ y\o e^-,\ x_{-(\al+\b)}\o e^+ \}\cong L(\la_1),\\
			A(V_{A_2}).(x_{\al+\b}\o e^-)&=\spn\{x_{\al+\b}\o e^-,\ y\o e^+,\ (1-y)\o e^- \}\cong L(\la_2). 
		\end{align*}
		Then $A(V_{A_2})\o _{A(V_{A_1})} \Om_{V_{A_1}}(V_{A_1+\frac{1}{2}\al})\cong L(\la_1)\op L(\la_2)$. This proves \eqref{eq:A1A2ind2}. 
		
		By \eqref{eq:res} and the fact that the $V_{A_1}$-adjoint module $V_{A_1}$ is the generalized Verma module associated to $\C\vac$, we have 
		$$\Res^{V_{A_2}}_{V_{A_1}}(V_{A_2})=\Phi^\L_{V_{A_1}}\left(\Res^{A(V_{A_2})}_{A(V_{A_1})}(\Om_{V_{A_2}}(V_{A_2}))\right)=\Phi^\L_{V_{A_1}}(\C\vac)\cong V_{A_1}.$$
		On the other hand, obviously the $\sl_3$-modules $L(\la_1)$ and $L(\la_2)$ decompose into $\C\vac \op L(\frac{1}{2}\al)$ when viewed as $\sl_2$-modules via the embedding \eqref{eq:iota}. It follows that 
		$$\Res^{V_{A_2}}_{V_{A_1}}(V_{A_2+\la_i})=\Phi^\L_{V_{A_1}}\left(\Res^{A(V_{A_2})}_{A(V_{A_1})}(\Om_{V_{A_2}}(V_{A_2+\la_i}))\right)=\Phi^\L_{V_{A_1}}(\C\vac\op L(\al/2))\cong V_{A_1}\op V_{A_1+\frac{1}{2}\al}.$$
		This proves \eqref{eq:A1A2res1} and \eqref{eq:A1A2res2}. 
	\end{proof}


	\section{Finite-type character ring of VOAs and Artin's induction theorem}\label{sec:5}
	This section aims to address the question of whether a module over a VOA $V$ can be obtained through finite induction from modules over subVOAs $U \subset V$. To properly formulate this question, we introduce the notion of the finite-type character ring of a VOA.
	
	\subsection{Trace functions on $A(V)$} 
	
	Let $W=\bigoplus_{n=0}^\infty W(n)$ be an irreducible ordinary $V$-module of conformal weight $h$. 
	The one-point correlation function on the torus associated to $W$ is the following power series \cite{FLM,Z}:
	\begin{equation}
		Z_W(a,\tau)=\tr\nolimits_{W} o(a) q^{L(0)-\frac{c}{24}}
		=\sum_{n\in \N} \tr\nolimits_{W(n)}o(a) q^{n+h-\frac{c}{24}},
	\end{equation}
	where $c$ is the central charge of $V$ and $W(0)=\Omega(W)$. 
	Its top-degree summand $\tr\nolimits_{\Omega(W)}o(a)\, q^{h-\frac{c}{24}}$ can be viewed as a function on $A(V)$, since the map 
	\[
	o:A(V)\longrightarrow \End(\Omega(W)),\qquad [a]\mapsto o(a)
	\]
	is a representation of $A(V)$. This leads to the following definition. 
	
	\begin{df}\label{def:finitecharactr}
Let $\Ord^{\fin.}(V)$ be the subcategory of $\Ord(V)$  whose objects are ordinary $V$-modules $W$ such that $\Om(W)$ is finite-dimensional. We define the {\bf finite-type character of $W$} to be the trace function
		\begin{equation}
			\chi_{W}: A(V)\longrightarrow \C,\qquad \chi_W([a])=\tr\nolimits_{\Omega(W)} o(a). 
		\end{equation}
		Then $\chi_W\in \mathrm{SF}(A(V)):=\{f:A(V)\to \C \mid f([a]\ast [b])=f([b]\ast [a]),\ \forall [a],[b]\in A(V)\}$, 
		the space of symmetric functions on $A(V)$. Define
		\[
		R^\fin(V):=\sum_{W\in \Ord^{\fin.}(V)} \Z\cdot  \chi_W.
		\]
		Then $R^\fin(V)$ is a sub-abelian group of $\mathrm{SF}(A(V))$. 
	\end{df}
	
Note that if $W=\bigoplus_{n=0}^\infty W(n)$ is an irreducible ordinary $V$-module, with $W(n)$ being an eigenspace of $L(0)$ of eigenvalue $h+n$ for all $n\in \N$, then the bottom degree $W(0)=\Om(W)$ is finite-dimensional. In particular $W\in \Ord^{\fin.}(V)$. 
	
	\begin{lm}\label{lm:finitechar}
		Let $V$ be a strongly rational VOA with irreducible modules $\{V=W^0,\dots,W^r \}$. Then the irreducible finite characters $\{\chi_{W^0},\dots,\chi_{W^r}\}$ form a basis of $\mathrm{SF}(A(V))$. In particular,
		\begin{equation}\label{eq:finchar}
			R^\fin(V)=\bigoplus_{i=0}^r \Z \cdot \chi_{W^i}. 
		\end{equation}
		Moreover, $R^\fin (V)$ is a commutative ring with respect to the following fusion product:
		\begin{equation}\label{eq:prodonRfin}
			\chi_{W^i}\ast \chi_{W^j}:= \chi_{W^i\boxtimes_{P(z)}W^j}
			=\sum_{k=0}^r N_{ij}^k \chi_{W^k},
		\end{equation}
		where $W^i\boxtimes_{P(z)}W^j$ denotes Huang–Lepowsky’s $P(z)$-tensor product \cite{HL,H05}. 
		We call $R^\fin(V)$ the {\bf finite-type character ring} of the VOA $V$. 
	\end{lm}
	
	\begin{proof}
		Since $A(V)$ is semisimple, it was proved in \cite[Lemma 5.3.3]{Z} that $\{\chi_{W^0},\dots,\chi_{W^r}\}$ form a basis of $\mathrm{SF}(A(V))$. 
		Given any $M,N\in \Ord(V)$, we have  
		\[
		\chi_{M\oplus N}([a])=\tr\nolimits_{\Omega(M)\oplus \Omega(N)}o(a)
		=\chi_{M}([a])+\chi_N([a]),
		\]
		for any $[a]\in A(V)$. Hence $R^\fin(V)=\bigoplus_{i=0}^r \Z \cdot \chi_{W^i}$. 
		
		Finally, the fact that $R^\fin(V)$ forms a commutative ring with respect to \eqref{eq:prodonRfin} follows from the commutativity and associativity of the fusion rules \cite{H05,GL25} and the decomposition of the $P(z)$-tensor product module 
		$W^i\boxtimes_{P(z)}W^j=\bigoplus_{k=0}^r N_{ij}^k W^k$ for strongly rational VOAs \cite{HL}. 
	\end{proof}
	
	Let $U\hookrightarrow V$ be a VOA embedding, and assume that $V$ is strongly rational. 
	Then for any irreducible ordinary $U$-module $M\in \Ord(U)$, the finite induction
	\[
	\Ind^V_U(M)
	=\Phi_V^\L\!\left(A(V)\otimes_{A_U}\frac{\Omega_U(M)}{\ker(\pi)\cdot\Omega_U(M)}\right)
	=\bigoplus_{i=0}^r m_i W^i
	\]
	is a direct sum of irreducible ordinary $V$-modules. For the character $\chi^U_M\in R^\fin(U)$, we define 
	\begin{equation}
		\Ind^V_U (\chi^U_M):=\chi_{\Ind^V_U(M)}=\sum_{i=0}^r m_i \chi^V_{W^i}\in R^\fin(V),
	\end{equation}
	and call it the {\bf induced finite character} of $M$ with respect to $U\hookrightarrow V$. 
	
	If, furthermore, $U$ is also strongly rational, then $R^\fin(U)$ is a free abelian group generated by the irreducible finite characters of $U$ by Lemma~\ref{lm:finitechar}, and 
	\begin{equation}
		\Ind^V_U: R^\fin(U)\longrightarrow R^\fin(V),\quad 
		\chi^U_M\longmapsto \Ind^V_U (\chi^U_M)
	\end{equation}
	is a homomorphism of abelian groups. 
	
	\begin{remark}
		In the case of finite groups, the group algebra $\C[G]$ has a canonical basis given by the elements of the group $G$. This basis induces the Schur inner product  
		\[
		(\chi_M, \chi_N) = \frac{1}{|G|}\sum_{g \in G} \chi_M(g)\chi_N(g^{-1}),\quad M,N\in \mathrm{Irr}(G),
		\]
		on the character ring $R(G)$. Moreover, the coset representatives of a subgroup $H \leq G$ give a concrete expression for the induced character:
		\[
		\Ind_H^G(\chi)(g) = \sum_{\substack{g_0 \in G ,\  g_0^{-1} g g_0 \in H}} \chi(g_0^{-1} g g_0),\quad \chi\in R(G).
		\]
		However, in the VOA situation, the Zhu algebra $A(V)$ does not possess a canonical basis analogous to that of $\C[G]$. In fact, $A(V)$ is more closely related to the universal enveloping algebra $U(\mathfrak{g})$ than to a group algebra. It is not yet known whether the finite character ring $R^{\mathrm{fin}}(V)$ has properties similar to those of $R(G)$.
		
	\end{remark}
	
	\subsection{Artin's induction theorem for finite-type characters}
	Let $V$ be a strongly rational VOA with central charge $c$, and let $U=\<\omega\>$ be the subVOA of $V$ generated by the Virasoro element $\omega$. Then $U \hookrightarrow V$ is a conformal embedding of VOAs, and $U$ is a quotient of the universal Virasoro VOA $\bar{V}(c,0)$. Since $A(V)$ is semisimple \cite{DLM1,Z}, write 
	\[
	A(V)=\prod_{i=0}^r \End(W^i(0)) \cong \prod_{i=0}^r M_{n_i}(\C).
	\]
	Note that $[\omega]\in A(V)$ acts as the scalar $h_i \cdot \Id$ on the irreducible $A(V)$-module $W^i(0)$ for $0\leq i\leq r$, where $h_i$ is the conformal weight of $W^i$. Therefore, the element $[\omega]\in A(V)$ can be expressed as 
	\[
	[\omega]=(h_0 I_{n_0},\dots,h_r I_{n_r})\in \prod_{i=0}^r M_{n_i}(\C),
	\]
	where $I_{n}\in M_n(\C)$ is the identity matrix. It follows that 
	\[
	A_U=\pi(A(U))\cong \C[y]/\<(y-h_0)\cdots(y-h_r)\>\subseteq A(V).
	\]
	Consider the following diagram:
	\[
	\begin{tikzcd}
		\bar{V}(c,0)\arrow[r,two heads] \arrow[d]& U\arrow[r,hook]\arrow[d]& V\arrow[d]\\
		\C[y]\arrow[r,two heads]& A(U)=\C[y]/I\arrow[r,"\pi"]& A(V)
	\end{tikzcd}
	\]
	It is clear that $\ker(\pi)=(\<(y-h_0)\cdots(y-h_r)\>+I)/I$, see \eqref{eq:AL10} for a special case. In particular, for the irreducible $U$-module $L(c,h_j)$, we have $\ker(\pi)\cdot v_{c,h_j}=0$ for all $0\leq j\leq r$. 
	Note that 
	\[
	A(V)\otimes_{A_U}\frac{\Omega_U(L(c,h_j))}{\ker(\pi)\cdot\Omega_U(L(c,h_j))}\cong 
	\left(\prod_{i=0}^r M_{n_i}(\C)\right)\otimes_{\C[[\omega]]} \C v_{c,h_j}
	\cong \prod_{i=0}^r \left(M_{n_i}(\C)\otimes_{\C[[\omega]]}\C v_{c,h_j}\right).
	\]
	If $h_i\neq h_j$, then the following relation holds in $M_{n_i}(\C)\otimes_{\C[[\omega]]}\C v_{c,h_j}$:
	\[
	h_i\cdot I_{n_i}\otimes v_{c,h_j}-I_{n_i}\otimes h_j\cdot v_{c,h_j}
	=I_{n_i}\ast [\omega]\otimes v_{c,h_j}-I_{n_i}\otimes o(\omega)\cdot v_{c,h_j}=0.
	\]
	Hence $(h_i-h_j)\cdot (I_{n_i}\otimes v_{c,h_j})=0$, and so $M_{n_i}(\C)\otimes_{\C[[\omega]]}\C v_{c,h_j}=0$. On the other hand, if $h_i=h_j$, we have
	\[
	M_{n_i}(\C)\otimes_{\C[[\omega]]}\C v_{c,h_j}\cong M_{n_i}(\C)\cong n_i W^i(0)
	\]
	as left $A(V)$-modules. In particular, the induced $V$-module has the following decomposition:
	\[
	\Ind^V_{U} (L(c,h_j))\cong \bigoplus_{\substack{0\leq i\leq r \\ h_i=h_j}} n_i W^i.
	\]

	Now, if the conformal weights $h_0,\dots,h_r\in \C$ of irreducible $V$-modules are pairwise distinct, then $\Ind^V_U (L(c,h_i))\cong n_i W^i$ as $V$-modules for all $i$. In this case, the irreducible character $\chi^V_{W^i}\in R^{\mathrm{fin}}(V)$ (see \eqref{eq:finchar}) can be written as
	\[
	\chi^V_{W^i}=\frac{1}{n_i}\Ind^V_U(\chi^U_{L(c,h_i)}),\quad 0\leq i\leq r.
	\]
	
	This leads to the following analogue of Artin’s induction theorem, see Proposition~\ref{prop:rankonevir} for an example: 
	
	\begin{thm}\label{thm:Artininduction}
		Let $V$ be a strongly rational VOA. Assume that the conformal weights $h_0,\dots,h_r\in \C$ of irreducible $V$-modules are pairwise distinct. Then every element in $R^{\mathrm{fin}}(V)$ can be written as a $\Q$-linear combination of induced characters from subVOAs of $V$.
	\end{thm}
	
	
	
	
	\begin{remark}
		Other than the Virasoro subVOA generated by $\omega$, a strongly rational VOA $V$ generally does not possess many canonical subVOAs. Therefore, it is in general difficult to establish a stronger form of Brauer’s induction theorem for finite characters.
	\end{remark}



	\section{Motivation example: the rank-two parabolic-type VOA $V_P\hr V_{A_2}$}
	\label{sec:6}

	In the remaining sections, we study the example that motivates our construction of the finite module induction functor \eqref{eq:inddiagram}, namely, the rank-two parabolic-type subVOA $V_P$ \cite{Liu24} of the lattice (or affine) VOA $V_{A_2}=L_{\widehat{\mathfrak{sl}}_3}(1,0)$.
	
	Consider the type-$A_2$ root lattice $A_2=\mathbb{Z}\alpha\oplus\mathbb{Z}\beta$ \eqref{eq:paringinA2lattice}, with the standard $2$-cocycle $\ep:A_2\times A_2\ra \<\pm 1\>$ such that 
	\begin{equation}\label{eq:twococycle}
		\ep(\al,\al)=1,\quad \ep(\b,\b)=1,\quad \ep(\al,\b)=1,\quad \ep(\b,\al)=-1. 
	\end{equation}
	We refer to the additive submonoid 
	\[
	P=\mathbb{Z}\alpha\oplus\mathbb{Z}_{\geq 0}\beta \subseteq A_2
	\]
	as a parabolic-type submonoid, see \cite[Definition 2.1]{Liu24}. The submonoid $P\subseteq A_2$ gives rise to an associated subVOA of $V_{A_2}$:
	\begin{equation}\label{eq:defofVP}
		V_P=\bigoplus_{\gamma\in P=\mathbb{Z}\alpha\oplus\mathbb{Z}_{\geq 0}\beta} M_{\hat{\mathfrak{h}}}(1,\gamma)
		=\bigoplus_{m\in\mathbb{Z},\, n\in\mathbb{N}} M_{\hat{\mathfrak{h}}}(1,m\alpha+n\beta),
	\end{equation}
	which we call a {\bf rank-two parabolic-type (sub)VOA} of $V_{A_2}$. Note that the VOA embedding $V_P\hookrightarrow V_{A_2}$ is conformal.
	
	\begin{remark}
		We call $V_P$ a parabolic-type VOA because it shares many properties with a parabolic subalgebra of a semisimple Lie algebra, see \cite[Section~4]{Liu24}. In fact, the degree-one Lie subalgebra 
		\[
		\mathfrak{p}=(V_P)_1=\mathrm{span}\{x_{\pm\alpha},\, x_{\beta},\, x_{\alpha+\beta},\, x,\, y\}
		\]
		is the standard parabolic subalgebra of $\mathfrak{sl}_3=(V_{A_2})_1$, consisting of block upper-triangular matrices with respect to the standard basis \eqref{eq:serrerel}:
		\[
		\begin{bmatrix}
			\ast & \ast & \ast \\
			0 & \ast & \ast\\
			0 & \ast & \ast
		\end{bmatrix}.
		\]
		Therefore, we believe that the induction for the VOA embedding $V_P\hookrightarrow V_{A_2}$ is well motivated.
	\end{remark}


	We will determine the Zhu algebra $A(V_P)$
	in terms of generators and relations. 
	Then, we will classify all the irreducible $V_P$-modules and determine all the finite inductions from the irreducible $V_P$-modules with respect to the VOA embedding $V_P\hookrightarrow V_{A_2}$.

	
	


	\subsection{A spanning set of $O(V_{P})$}\label{Sec:8.1}
	We first give a concrete description of $O(V_P)$ in \eqref{eq:defofOV}, and then use it to determine the structure of $A(V_P)$. 
	For $A_1=\mathbb{Z}\alpha$, recall that $O(V_{A_1})$ is spanned by the following elements, see \cite[Corollary~6.14]{Liu24}:
	\begin{equation}\label{3.29}
		\begin{cases}
			&\alpha(-n-2)u+\alpha(-n-1)u, \qquad u\in V_{A_1},\ n\geq 0,\\[4pt]
			&\pm \alpha(-1)v+v, \qquad v\in M_{\hat{\mathfrak{h}}}(1,\pm\alpha),\\[4pt]
			&M_{\hat{\mathfrak{h}}}(1,\pm k\alpha), \qquad k\geq 2,\\[4pt]
			&\alpha(-1)^3w-\alpha(-1)w, \qquad w\in M_{\hat{\mathfrak{h}}}(1,0).
		\end{cases}
	\end{equation}
	
	\begin{lm}\label{lm4.1}
		Suppose $\gamma,\theta\in P$ such that $\frac{(\gamma|\gamma)}{2}=N\geq 1$ and $(\gamma|\theta)=n\geq 1$. Then $e^{\gamma+\theta}\in O(V_P)$.
	\end{lm}
	\begin{proof}
		By the definition of the lattice vertex operator $Y(e^\gamma,z)e^\theta$, we have 
		\[
		e^\gamma_{-n-1}e^\theta
		=\mathrm{Res}_{z}\, z^{-n-1}E^{-}(-\gamma,z)z^{(\gamma|\theta)}\varepsilon(\gamma,\theta)e^{\gamma+\theta}
		=\varepsilon(\gamma,\theta)e^{\gamma+\theta},
		\]
		and $e^{\gamma}_{-m}e^\theta=0$ for all $m\leq n$. Since $n\geq 1$ and $\mathrm{wt}\,e^\gamma=\frac{(\gamma|\gamma)}{2}=N\geq 1$, we obtain
		\begin{align*}
			\mathrm{Res}_{z}\, Y(e^{\gamma},z)e^{\theta} \frac{(1+z)^N}{z^{1+n}}
			&=e^\gamma_{-n-1}e^\theta+\binom{N}{1}e^\gamma_{-n}e^\theta+\dots+\binom{N}{N}e^\gamma_{-n-1+N}e^\theta\\
			&=\varepsilon(\gamma,\theta)e^{\gamma+\theta}\in O(V_P).
		\end{align*}
		This shows $e^{\gamma+\theta}\in O(V_P)$ since $\varepsilon(\gamma,\theta)\in\{\pm1\}$.
	\end{proof}
	
	\begin{lm}\label{lm4.2}
		Let 
		\[
		S:=\{e^{m\alpha+n\beta}:m\in\mathbb{Z},\,n\in\mathbb{N}\}\setminus\{e^{\pm\alpha},e^{\beta},e^{\alpha+\beta}\}.
		\]
		Then $S\subset O(V_P)$.
	\end{lm}
	\begin{proof}
		For $m\geq 1$, since $(\alpha|m\alpha+\beta)=2m-1\geq 1$, by Lemma~\ref{lm4.1} and induction on $m$, we have $e^{(m+1)\alpha+\beta}\in O(V_P)$ for all $m\geq 1$. Similarly, since $(\beta|\alpha+n\beta)=2n-1\geq 1$ for $n\geq 1$, we obtain
		\begin{equation}\label{4.3}
			e^{\alpha+(n+1)\beta}\in O(V_P),\qquad n\geq 1.
		\end{equation}
		Now let $n\geq 2$ and assume that $e^{m\alpha+n\beta}\in O(V_P)$ for all $m\geq 1$. We want to show that $e^{m\alpha+(n+1)\beta}\in O(V_P)$ for all $m\geq 1$.
		
		Indeed, since $(m\alpha+n\beta|\alpha+\beta)=m((\alpha|\alpha)+(\alpha|\beta))+n((\beta|\alpha)+(\beta|\beta))=m+n\geq 1$, by Lemma~\ref{lm4.1} we have $e^{(m+1)\alpha+(n+1)\beta}\in O(V_P)$ for all $m\geq 1$. Thus $e^{m\alpha+(n+1)\beta}\in O(V_P)$ for all $m\geq 2$, and by \eqref{4.3}, the same holds for $m=1$. This completes the induction step and shows that $e^{m\alpha+n\beta}\in O(V_P)$ for all $m\geq 1$, $n\geq 2$. Hence,
		\begin{equation}\label{4.4}
			S_{1}:=\{e^{m\alpha+n\beta}:m\geq 1,n\geq 2\}\cup\{e^{m\alpha+\beta}:m\geq 2\}\subset O(V_P).
		\end{equation}
		
		On the other hand, for any $m\geq 1$, since $(-m\alpha|\beta)=m\geq 1$, we have $e^{-m\alpha+\beta}\in O(V_P)$ for all $m\geq 1$. Similarly, since $(-\alpha+n\beta|\beta)=1+2n\geq 1$ for all $n\geq 0$, we have 
		\begin{equation}\label{4.4'}
			e^{-\alpha+(n+1)\beta}\in O(V_P),\qquad n\geq 0. 
		\end{equation}
		Using the fact that $(-m\alpha+n\beta|-\alpha+\beta)=3m+3n\geq 1$ for $m,n\geq 1$, together with \eqref{4.4'}, we can similarly show that
		\begin{equation}\label{4.5}
			S_2:=\{e^{-m\alpha+n\beta}: m\geq 1, n\geq 1\}\subset O(V_P).
		\end{equation}
		
		Finally, for any $m,n\geq 1$, since $(\alpha|m\alpha)=(-\alpha|-m\alpha)=2m>1$ and $(\beta|n\beta)=2n>1$, by Lemma~\ref{lm4.1} again we see that $e^{\pm m\alpha}\in O(V_P)$ and $e^{n\beta}\in O(V_P)$ for all $m\geq 2$, $n\geq 2$. Combining \eqref{4.4} and \eqref{4.5}, we conclude that
		$	S=S_1\cup S_2\cup\{e^{\pm m\alpha}:m\geq 2\}\cup\{e^{n\beta}:n\geq 2\}\subset O(V_P).
		$
	\end{proof}
	
	\begin{df}\label{df4.3}
		Let $O$ be the subspace of $V_P$ spanned by the following elements:
		\begin{equation}\label{4.7}
			\begin{cases}
				&h(-n-2)u+h(-n-1)u,\qquad u\in V_P,\ h\in \mathfrak{h},\ n\geq 0;\\[4pt]
				&\gamma(-1)v+v,\qquad v\in M_{\hat{\mathfrak{h}}}(1,\gamma),\ \gamma\in \{\alpha,-\alpha,\beta,\alpha+\beta\};\\[4pt]
				&\gamma(-1)^2v+\gamma(-1)v,\quad v\in M_{\hat{\mathfrak{h}}}(1,\gamma+\gamma'),\ \gamma,\gamma'\in \{\alpha,-\alpha,\beta,\alpha+\beta\},\gamma+\gamma'\in \{\alpha+\beta,\beta\};\\[4pt]
				&M_{\hat{\mathfrak{h}}}(1,m\alpha+n\beta),\quad m\alpha+n\beta\in (\mathbb{Z}\alpha\oplus\mathbb{Z}_{\geq 0}\beta)\setminus\{0,\alpha,-\alpha,\beta,\alpha+\beta\};\\[4pt]
				&\alpha(-1)^3w-\alpha(-1)w,\qquad w\in M_{\hat{\mathfrak{h}}}(1,0).
			\end{cases}
		\end{equation}
	\end{df}
	
	Note that the only possible ordered pairs $(\gamma,\gamma')$ such that $\gamma,\gamma'\in \{\alpha,-\alpha,\beta,\alpha+\beta\}$ and $\gamma+\gamma'\in \{\alpha+\beta,\beta\}$, as in \eqref{4.7}, are contained in the following set:
	\begin{equation}\label{4.8}
		\{(\alpha,\beta),(\beta,\alpha),(-\alpha,\alpha+\beta),(\alpha+\beta,-\alpha)\}.
	\end{equation}
	Hence, the elements $\gamma(-1)^2v+\gamma(-1)v$ in \eqref{4.7} can be written more explicitly as
	\begin{align}
		&\alpha(-1)^2v+\alpha(-1)v, &&\beta(-1)^2v+\beta(-1)v, &&v\in M_{\hat{\mathfrak{h}}}(1,\alpha+\beta),\label{extra1}\\
		&\alpha(-1)^2v-\alpha(-1)v, &&(\alpha+\beta)(-1)^2v+(\alpha+\beta)(-1)v, &&v\in M_{\hat{\mathfrak{h}}}(1,\beta).\label{extra2}
	\end{align}
	Moreover, we observe that $O(V_{\mathbb{Z}\alpha})\subset O$ and $O(V_{\mathbb{Z}_{\geq 0}\gamma})\subset O$ for any $\gamma\in\{\alpha,-\alpha,\beta,\alpha+\beta\}$, since $M_{\widehat{\mathbb{C}\alpha}}(1,\pm\alpha)\subset M_{\hat{\mathfrak{h}}}(1,\pm\alpha)$ and $M_{\widehat{\mathbb{C}\gamma}}(1,\gamma)\subset M_{\hat{\mathfrak{h}}}(1,\gamma)$ by \eqref{3.29}.
	
	\begin{remark}
		Figure~\ref{fig4} illustrates the definition of $O$. The black dots in the diagram represent elements in the parabolic-type submonoid $P=\mathbb{Z}\alpha\oplus\mathbb{Z}_{\geq 0}\beta$ of $A_2$. Except for the roots represented by the red vectors, the Heisenberg modules $M_{\hat{\mathfrak{h}}}(1,\gamma)$ associated to all other dots $\gamma$ are contained in the subspace $O$.
		\begin{figure}
			\centering
			\begin{tikzpicture}
				\coordinate (Origin)   at (0,0);
				\coordinate (XAxisMin) at (-5,0);
				\coordinate (XAxisMax) at (5,0);
				\coordinate (YAxisMin) at (0,0);
				\coordinate (YAxisMax) at (0,5);
				
				\clip (-5.1,-1.9) rectangle (5.1,3.7); 
				\begin{scope} 
					\pgftransformcm{1}{0}{1/2}{sqrt(3)/2}{\pgfpoint{0cm}{0cm}} 
					\coordinate (Bone) at (0,2);
					\coordinate (Btwo) at (2,-2);
					\draw[style=help lines,dashed] (-7,-6) grid[step=2cm] (6,6);
					\foreach \x in {-4,-3,...,4}{
						\foreach \y in {0,1,...,4}{
							\coordinate (Dot\x\y) at (2*\x,2*\y);
							\node[draw,circle,inner sep=2pt,fill] at (Dot\x\y) {};
						}
					}
					\draw [thick,-latex,red] (Origin) node [xshift=0cm, yshift=-0.3cm] {$0$}
					-- (Bone) node [above right]  {$\al+\b$};
					\draw [thick,-latex,red] (Origin)
					-- ($(Bone)+(Btwo)$) node [below right] {$\al$};
					\draw [thick,-latex,red] (Origin)
					-- ($-1*(Bone)-1*(Btwo)$) node [below left] {$-\al$};
					\draw [thick,-latex,red] (Origin)
					-- ($-1*(Btwo)$) coordinate (B3) node [above left] {$\b$};
				\end{scope} 
				\begin{scope}
					\pgftransformcm{1}{0}{-1/2}{sqrt(3)/2}{\pgfpoint{0cm}{0cm}} 
					\draw[style=help lines,dashed] (-6,-6) grid[step=2cm] (6,6);
				\end{scope}
			\end{tikzpicture}
			\caption{\label{fig4}}
		\end{figure}
		
	\end{remark}

	\subsection{Proof of the main theorem}
	In this subsection, we show that $O$ coincides with $O(V_P)$.
	
	\begin{prop}\label{prop4.4}
		Let $O$ be the subspace given in Definition~\ref{df4.3}. Then $O\subseteq O(V_P)$.
	\end{prop}
	
	\begin{proof}
		It is clear that $h(-n-2)u+h(-n-1)u\in O(V_P)$ for any $h\in \h$, $u\in V_P$, and $n\geq 0$.  
		By the congruence 
		\[
		h(-m)v\equiv (-1)^{m-1}v\ast (h(-1)\vac)\pmod{O(V_P)}
		\]
		and the fact that $O(V_P)$ is a two-sided ideal with respect to $\ast$, we obtain 
		$h(-m)O(V_P)\subseteq O(V_P)$ for any $h\in \h$ and $m\geq 1$.  
		Then, by Lemma~\ref{lm4.2}, we have 
		\[
		M_{\hat{\h}}(1,m\al+n\b)\subseteq O(V_P)
		\quad \text{for } m\al+n\b \in P\setminus\{\al,-\al,\b,\al+\b\}.
		\]
		
		Moreover, for $\ga\in\{\al,-\al,\b,\al+\b\}$, since $(\ga|\ga)=2$, we obtain
		\[
		\ga(-1)e^\ga+e^\ga=e^\ga_{-2}\vac+e^\ga_{-1}\vac=e^\ga\circ\vac\in O(V_P),
		\]
		and hence $\ga(-1)v+v\in O(V_P)$ for any $v\in M_{\hat{\h}}(1,\ga)$ and $\ga\in\{\al,-\al,\b,\al+\b\}$.
		
		Now suppose $\ga,\ga'\in\{\al,-\al,\b,\al+\b\}$ are such that $\ga+\ga'\in\{\b,\al+\b\}$.  
		Since $h(-n-2)u+h(-n-1)u\equiv 0\pmod{O(V_P)}$ for all $h\in\h$, $n\geq 0$, and $u\in V_P$, we have
		\begin{align*}
			0&\equiv e^{\ga}\circ e^{\ga'}=\frac{1}{2}\epsilon(\ga,\ga')\ga(-2)e^{\ga+\ga'}+\frac{1}{2!}\epsilon(\ga,\ga')\ga(-1)^2e^{\ga+\ga'}+\epsilon(\ga,\ga')\ga(-1)e^{\ga+\ga'}\\
			&\equiv \frac{1}{2}\epsilon(\ga,\ga')\big(\ga(-1)^2e^{\ga+\ga'}+\ga(-1)e^{\ga+\ga'}\big)\pmod{O(V_P)}.
		\end{align*}
		Thus, $\ga(-1)^2v+\ga(-1)v\in O(V_P)$ for any $v\in M_{\hat{\h}}(1,\ga+\ga')$.
		
		Finally, since $O(V_{\Z\al})\subset O(V_P)$, it follows from~\eqref{3.29} that 
		\[
		\al(-1)^3\vac-\al(-1)\vac\in O(V_P).
		\]
		Consequently, $\al(-1)^3w-\al(-1)w\in O(V_P)$ for any $w\in M_{\hat{\h}}(1,0)$, as $h(-m)O(V_P)\subseteq O(V_P)$ for all $h\in\h$ and $m\geq 1$.
	\end{proof}
	
	Conversely, to prove $O(V_P)\subseteq O$, we need to show that 
	\begin{equation}\label{eq:goal}
		M_{\hat{\h}}(1,\eta)\circ M_{\hat{\h}}(1,\theta)\subset O, 
		\qquad \text{for all } \eta,\theta\in P=\Z\al\op\Z_{\geq 0}\b.
	\end{equation}
	By the construction of $O$ in~\eqref{4.7} and the fact that 
	$Y(e^{\eta},z)e^{\theta}\in M_{\hat{\h}}(1,\eta+\theta)((z))$, we have 
	\[
	M_{\hat{\h}}(1,\eta)\circ M_{\hat{\h}}(1,\theta)\subset M_{\hat{\h}}(1,\eta+\theta)\subset O,
	\quad \text{if }\ \eta+\theta\in P\setminus\{0,\al,-\al,\b,\al+\b\}.
	\]
	Hence, it remains to show that 
	\begin{align}
		&M_{\hat{\h}}(1,0)\circ M_{\hat{\h}}(1,\ga)\subset O,&&
		M_{\hat{\h}}(1,\ga)\circ M_{\hat{\h}}(1,0)\subset O; \label{4.10}\\
		&M_{\hat{\h}}(1,\al)\circ M_{\hat{\h}}(1,\b)\subset O,&&
		M_{\hat{\h}}(1,\b)\circ M_{\hat{\h}}(1,\al)\subset O; \label{4.11}\\
		&M_{\hat{\h}}(1,\al+\b)\circ M_{\hat{\h}}(1,-\al)\subset O,&&
		M_{\hat{\h}}(1,-\al)\circ M_{\hat{\h}}(1,\al+\b)\subset O; \label{4.12}\\
		&M_{\hat{\h}}(1,\al)\circ M_{\hat{\h}}(1,-\al)\subset O,&&
		M_{\hat{\h}}(1,-\al)\circ M_{\hat{\h}}(1,\al)\subset O, \label{4.13}
	\end{align}
	where $\ga\in\{\al,-\al,\b,\al+\b\}$.
	
	The inclusion~\eqref{4.10} can be proved by arguments similar to those used in 
	\cite[Proposition~6.5, Lemma~6.6, and Proposition~6.7]{Liu24}, so we omit the details. 

	\subsubsection{Proof of \eqref{4.11} and \eqref{4.12}} Let $(\ga,\ga')$ be an ordered pair in the set~\eqref{4.8}.  
	Given a spanning element $u=h^1(-n_1)\cdots h^r(-n_r)e^\ga$ of $M_{\hat{\h}}(1,\ga)$ and 
	$v=h^1(-m_1)\cdots h^s(-m_s)e^{\ga'}$ of $M_{\hat{\h}}(1,\ga')$,  
	we need to show that 
	\[
	\Res_z Y(u,z)v\,\frac{(1+z)^{\wt u}}{z^{2}}\equiv 0\pmod{O}.
	\]
	For $u=e^\ga$, this congruence holds because of the following (stronger) statement.
	
	\begin{prop}\label{prop4.5}
		Let $(\ga,\ga')$ be an ordered pair in the set~\eqref{4.8}, and let $n\geq 0$. Then 
		\begin{equation}\label{4.14}
			\Res_z Y(e^{\ga},z)\big(h^1(-n_1)\cdots h^r(-n_r)e^{\ga'}\big)\frac{(1+z)}{z^{2+n}}\in O,
		\end{equation}
		where $r\geq 0$, $h^i\in \h$ for all $i$, and $n_1\geq \cdots \geq n_r\geq 1$.
	\end{prop}
	
	\begin{proof}
		It is clear from~\eqref{4.7} that $h(-m)O\subset O$ for any $m\geq 1$ and $h\in\h$.  
		We first claim that $L(-1)u+L(0)u\in O$ for all $u\in V_P$.  
		Indeed, if $u\in M_{\hat{\h}}(1,m\al+n\b)$ with $m\al+n\b\in P\setminus\{\al,-\al,\b,\al+\b\}$,  
		then $L(-1)u+L(0)u\in M_{\hat{\h}}(1,m\al+n\b)\subset O$ by~\eqref{4.7}.  
		
		Now let $u=h^1(-n_1)\cdots h^r(-n_r)e^\ga$, where $\ga\in\{\al,-\al,\b,\al+\b\}$,  
		$h^i\in\h$ for all $i$, and $n_1\geq\cdots\geq n_r\geq 1$.  
		Since $L(-1)e^\ga=\ga(-1)e^\ga$, we obtain
		\begin{align*}
			L(-1)u+L(0)u
			&=h^1(-n_1)\cdots h^r(-n_r)(\ga(-1)e^\ga+e^\ga)\\
			&\quad+\sum_{j=1}^r (h^j(-n_j-1)+h^j(-n_j))\,h^1(-n_1)\cdots \widehat{h^j(-n_j)}\cdots h^r(-n_r)e^\ga\\
			&\equiv 0\pmod{O},
		\end{align*}
		by~\eqref{4.7}. Hence $L(-1)u+L(0)u\in O$ for all $u\in V_P$.
		
		Now take $(\ga,\ga')\in\{(\al,\b),(\b,\al),(-\al,\al+\b),(\al+\b,-\al)\}$.  
		Note that $\epsilon(\ga,\ga')=-1$.
		We first use induction on $n\geq 0$ to prove that
		\begin{equation}\label{4.15}
			\Res_z Y(e^{\ga},z)e^{\ga'}\frac{(1+z)}{z^{2+n}}\in O,
		\end{equation}
		which will serve as the base case for induction on $r$ in~\eqref{4.14}.
		
		For $n=0$, we have
		\[
		\Res_z Y(e^\ga,z)e^{\ga'}\frac{(1+z)}{z^{2}}
		= e^\ga_{-2}e^{\ga'}+e^\ga_{-1}e^{\ga'}
		= -\frac{1}{2}\ga(-1)^2e^{\ga+\ga'}-\frac{1}{2}\ga(-1)e^{\ga+\ga'}
		\equiv 0\pmod{O}.
		\]
		Assume~\eqref{4.15} holds for smaller $n$. Then
		\begin{align*}
			&(n+1)(n+2)\big(e^\ga_{-n-3}e^{\ga'}+e^\ga_{-n-2}e^{\ga'}\big)
			=(n+1)(L(-1)e^\ga)_{-n-2}e^{\ga'}+(n+2)(L(-1)e^\ga)_{-n-1}e^{\ga'}\\
			&=(n+1)(\ga(-1)e^\ga)_{-n-2}e^{\ga'}+(n+2)(\ga(-1)e^\ga)_{-n-1}e^{\ga'}\\
			&=(n+1)\!\sum_{j\geq 0}\!\ga(-1-j)e^\ga_{-n-2+j}e^{\ga'}
			+(n+1)\!\sum_{j\geq 0}\! e^\ga_{-n-3-j}\ga(j)e^{\ga'}\\
			&\quad+(n+2)\!\sum_{j\geq 0}\!\ga(-1-j)e^\ga_{-n-1+j}e^{\ga'}
			+(n+2)\!\sum_{j\geq 0}\! e^\ga_{-n-2-j}\ga(j)e^{\ga'}\\
			&=(n+1)\sum_{j\geq 0}\ga(-1-j)e^\ga_{-n-2+j}e^{\ga'}
			+(n+1)(\ga|\ga')e^{\ga}_{-n-3}e^{\ga'}\\
			&\quad+(n+2)\sum_{j\geq 0}\ga(-1-j)e^\ga_{-n-1+j}e^{\ga'}
			+(n+2)(\ga|\ga')e^{\ga}_{-n-2}e^{\ga'}.
		\end{align*}
		Since $(\ga|\ga')=-1$ and $\ga(-2-j)v+\ga(-1-j)v\in O$ for any $j\geq 0$, we have \begin{align*} &(n+1)(n+3)\left(e^\ga_{-n-3}e^{\ga'}+e^\ga_{-n-2}e^{\ga'}\right)\\
			&=(n+1)\ga(-1)e^\ga_{-n-2}e^{\ga'}+(n+1)\sum_{t\geq 0}\ga(-2-t)e^\ga_{-n-1+t}e^{\ga'}\\
			&\quad +(n+2)\sum_{j\geq 0} \ga(-1-j)e^\ga_{-n-1+j}e^{\ga'}-e^\ga_{-n-2}e^{\ga'}\\ &\equiv (n+1)\ga(-1)e^\ga_{-n-2}e^{\ga'}+\sum_{j\geq 0} \ga(-1-j) e^\ga_{-n-1+j}e^{\ga'}-e^\ga_{-n-2}e^{\ga'}\pmod{O}\\ &\equiv (n+1)\ga(-1)e^\ga_{-n-2}e^{\ga'}+\sum_{j\geq 0} (-1)^j\ga(-1) e^\ga_{-n-1+j}e^{\ga'}-e^\ga_{-n-2}e^{\ga'}\pmod{O}. 
		\end{align*}
		
		Since $\ga(-1)O\subset O$, then by the induction hypothesis we have \begin{align*} \ga(-1)e^\ga_{-n-1+j}e^{\ga'}\equiv (-1)e^\ga_{-n-1+j-1}e^{\ga'}\equiv \ds\equiv (-1)^j \ga(-1)e^\ga_{-n-1}e^{\ga'}\pmod{O}, \end{align*} for any $0\leq j\leq n$, and $e^\ga_{-n-2}e^{\ga'}\equiv \ds \equiv (-1)^{n+1}e^\ga_{-1}e^{\ga'}=(-1)^{n+1}\epsilon(\ga,\ga')\ga(-1)e^{\ga+\ga'}\pmod{O}$. Moreover, we have $e^{\ga}_{0}e^{\ga'}=\epsilon(\ga,\ga')e^{\ga+\ga'}$, and $e^\ga_{m}e^{\ga'}=0$ for $m\geq 1$.
		It follows that \begin{align*} &(n+1)\ga(-1)e^\ga_{-n-2}e^{\ga'}+\sum_{j=0}^{n+1} (-1)^j\ga(-1) e^\ga_{-n-1+j}e^{\ga'}-e^\ga_{-n-2}e^{\ga'}\\ &\equiv (n+1) \ga(-1)e^{\ga}_{-n-2}e^{\ga'}+(n+1)\ga(-1)e^{\ga}_{-n-1}e^{\ga'}+ (-1)^{n+1} \ga(-1)e^{\ga}_{0}e^{\ga'}-e^\ga_{-n-2}e^{\ga'}\\ &\equiv 0+ (-1)^{n+1} \ga(-1)\epsilon(\ga,\ga')e^{\ga+\ga'}-(-1)^{n+1}\epsilon(\ga,\ga')\ga(-1)e^{\ga+\ga'}\\ &\equiv 0\pmod{O}. \end{align*}
		Hence $(n+1)(n+3)\big(e^\ga_{-n-3}e^{\ga'}+e^\ga_{-n-2}e^{\ga'}\big)\in O$, 
		completing the proof of~\eqref{4.15}.
		
		Finally, we use induction on the length $r$ to prove~\eqref{4.14}.  
		The base case $r=0$ follows from~\eqref{4.15}.  
		Assume~\eqref{4.14} holds for smaller $r\geq 1$. Then \begin{align*} &\Res_z Y(e^\ga,z)h^1(-n_1)\ds h^r(-n_r)e^{\ga'} \frac{(1+z)}{z^{2+n}}\\ &= h^1(-n_1)\ds h^r(-n_r) \Res_z Y(e^\ga,z)e^{\ga'} \frac{(1+z)}{z^{2+n}}\\ &-\sum_{j=1}^r (h^j|\ga) \Res_z h^1(-n_1)\ds h^{j-1}(-n_{j-1}) Y(e^\ga,z)h^{j+1}(-n_{j+1})\ds h^r(-n_r)e^{\ga'} \frac{(1+z)}{z^{2+n+n_j}}\\ &\equiv 0\pmod{O}, \end{align*}
		where the last congruence follows from the induction hypothesis and the fact that 
		$h(-m)O\subset O$ for all $h\in\h$ and $m\geq 1$.
	\end{proof}
	
	By a slight modification of the induction arguments in 
	\cite[Lemma~6.6 and Proposition~6.7]{Liu24}, one can show that for any $n\geq 0$,
	\begin{equation}\label{4.16}
		\Res_z Y\big(h^1(-n_1)\cdots h^r(-n_r)e^\ga,z\big)
		h^1(-m_1)\cdots h^s(-m_s)e^{\ga'}\frac{(1+z)^{n_1+\cdots+n_r+1}}{z^{2+n}}\in O.
	\end{equation}
	Note that the only properties of $O'$ used in the proof of 
	\cite[Lemma~6.6 and Proposition~6.7]{Liu24} are the relations 
	$\al(-n-2)v+\al(-n-1)v\in O'$ and the congruence analogous to~\eqref{4.14}, 
	both of which, in our rank-two parabolic case, are satisfied by~\eqref{4.7} 
	and Proposition~\ref{prop4.5}.
	
	Hence,~\eqref{4.11} and~\eqref{4.12} follow from~\eqref{4.16}.

	\subsubsection{Proof of \eqref{4.13}} Given a spanning element 
	\[
	u = h^1(-n_1)\cdots h^r(-n_r)e^\alpha \in M_{\hat{\h}}(1,\alpha)
	\quad \text{and} \quad 
	v = h^1(-m_1)\cdots h^s(-m_s)e^{-\alpha} \in M_{\hat{\h}}(1,-\alpha),
	\]
	we need to show that $u\circ v\in O$. Again, we only prove the base case when $u=e^\alpha$.
	
	\begin{prop}\label{prop4.8}
		For any $n\geq 0$, we have 
		\begin{equation}\label{4.17}
			\Res_z\, Y(e^\alpha,z)
			\big(h^1(-n_1)\cdots h^r(-n_r)e^{-\alpha}\big)
			\frac{(1+z)}{z^{2+n}}
			\in O,
		\end{equation}
		where $r\geq 0$, $h^i\in \h$ for all $i$, and $n_1\geq \cdots \geq n_r\geq 1$.
	\end{prop}
	\begin{proof}
		Again, we first prove \eqref{4.17} for $r=0$ by induction on $n\geq 0$. For the base case $n=0$, we note that 
		$$e^\al\circ e^{-\al}=e^\al_{-2}e^{-\al}+e^\al_{-1}e^{-\al}\equiv \frac{1}{6}(\al(-1)^3\vac-\al(-1)\vac)\equiv 0\pmod{O}. $$
		Suppose the conclusion holds for smaller $n\geq 1$. 
		Then by a similar calculation as Proposition~\ref{prop4.5}, with $\ga=\al$ and $\ga'=-\al$, noting that $e^\al_me^{-\al}=0$ for $m\geq 2$, we have 
		\begin{align*}
			&(n+1)(n+4) \left(e^\al_{-n-3}e^{-\al}+e^\al_{-n-2}e^{-\al}\right)\\
			&\equiv (n+1) \al(-1)e^\al_{-n-2}e^{-\al}+\sum_{j= 0}^{n+2}(-1)^j \al(-1)e^\al_{-n-1+j}e^{-\al}-2e^\al_{-n-2}e^{-\al}\pmod{O}\\
			& \equiv ((n+1) \al(-1)e^\al_{-n-2}e^{-\al}+ (n+1)\al(-1)e^{\al}_{-n-1}e^{-\al})\\
			&\ \ + (-1)^{n+1}\al(-1)e^\al_0e^{-\al}+(-1)^{n+2}\al(-1)e^\al_{1}e^{-\al}-2(-1)^{n+1}e^\al_{-1}e^{-\al}\pmod{O}\\
			&\equiv 0+ (-1)^{n+1}\epsilon(\al,-\al)\left(\al(-1)^2\vac-\al(-1)\vac-\al(-2)\vac-\al(-1)^2\vac\right)\\
			&\equiv 0\pmod{O}.
		\end{align*}
		This proves \eqref{4.17} for $r=0$.  
		The induction step for the general case $r\geq 1$ is similar to the proof of Proposition~\ref{prop4.5}, and we omit the details.
	\end{proof}
	
	With \eqref{4.10}--\eqref{4.13} and Proposition~\ref{prop4.4}, we have our final conclusion in this subsection:
	
	\begin{thm}\label{prop4.7}
		Let $P$ be the parabolic-type submonoid $\Z\al\op \Z_{\geq 0}\b$ of the root lattice $A_2=\Z\al\op \Z\b$. The subspace $O(V_P)$ of $V_P$ is equal to $O$ in Definition~\ref{df4.3}. 
	\end{thm}

	\subsection{The Zhu algebra of $V_{P}$}\label{Sec:8.2}
	With the explicit expression of $O(V_P)$ by \eqref{4.7} and Proposition~\ref{prop4.7}, we give a concrete description of Zhu algebra $A(V_P)$.

	\subsubsection{Generators and relations}By \eqref{4.7}, we have 
	\[
	M_{\hat{\h}}(1,m\alpha+n\beta)\subset O 
	\quad \text{for any } \; m\alpha+n\beta\in P\setminus\{0,\alpha,-\alpha,\beta,\alpha+\beta\}.
	\]
	It follows that 
	\[
	A(V_P)=V_P/O
	=[M_{\hat{\h}}(1,0)]+[M_{\hat{\h}}(1,\alpha)]+[M_{\hat{\h}}(1,-\alpha)]
	+[M_{\hat{\h}}(1,\beta)]+[M_{\hat{\h}}(1,\alpha+\beta)],
	\]
	where $[S]$ denotes the equivalence class of a subspace or element $S\subset V_P$. 
	Moreover, it is straightforward to verify that the relations in \eqref{4.7} imply 
	\begin{equation}\label{4.20}
		\begin{aligned}
			[M_{\hat{\h}}(1,0)]
			&=\C[[\alpha(-1)\vac],[\beta(-1)\vac]]\,\big/\,\<[\alpha(-1)\vac]^3-[\alpha(-1)\vac]\>,\\
			[M_{\hat{\h}}(1,\pm\alpha)]
			&=\spn\{[\beta(-1)^n e^{\pm\alpha}]: n\in \N\},\\
			[M_{\hat{\h}}(1,\beta)]
			&=\C[e^\beta]+\C[\alpha(-1)e^\beta]+\C[\alpha(-1)^2e^\beta],\\
			[M_{\hat{\h}}(1,\alpha+\beta)]
			&=\C[e^{\alpha+\beta}]+\C[\alpha(-1)e^{\alpha+\beta}]+\C[\alpha(-1)^2e^{\alpha+\beta}].
		\end{aligned}
	\end{equation}
	
	\begin{df}\label{def:relations}
		Let $A_P$ be the associative (unital) algebra defined by 
		\[
		A_P:=\C\<x,y,x_\alpha,x_{-\alpha},x_\beta,x_{\alpha+\beta}\>/R,
		\]
		where $\C\<x,y,x_\alpha,x_{-\alpha},x_\beta,x_{\alpha+\beta}\>$ denotes the tensor algebra on six generators 
		$x,y,x_\alpha,x_{-\alpha},x_\beta,x_{\alpha+\beta}$, 
		and $R$ is the two-sided ideal generated by the following relations:
		\begin{align}
			&xx_{\pm\alpha}=\pm x_{\pm\alpha},\quad 
			x_{\pm\alpha}x=\mp x_{\pm\alpha},\quad 
			x_\alpha x_{-\alpha}=\tfrac{1}{2}x^2+\tfrac{1}{2}x,\quad 
			x_{-\alpha}x_\alpha=\tfrac{1}{2}x^2-\tfrac{1}{2}x; \label{4.21}\\[4pt]
			&xy=yx,\quad 
			x^3-x=0,\quad 
			yx_{\pm\alpha}=x_{\pm\alpha}y\mp x_{\pm\alpha},\quad 
			x_\beta y+x_\beta=0,\quad 
			yx_\beta-x_\beta=0; \label{4.22}\\[4pt]
			&x_{\alpha+\beta}(x+y)+x_{\alpha+\beta}=0,\quad 
			(x+y)x_{\alpha+\beta}-x_{\alpha+\beta}=0; \label{4.23}\\[4pt]
			&xx_\beta-x_\beta x+x_\beta=0,\quad 
			xx_{\alpha+\beta}-x_{\alpha+\beta}x-x_{\alpha+\beta}=0; \label{4.24}\\[4pt]
			&x_\alpha x_\beta=-x_{\alpha+\beta}y,\quad 
			x_\beta x_\alpha=-x_{\alpha+\beta}y-x_{\alpha+\beta},\quad 
			x_{-\alpha}x_{\alpha+\beta}=-x_\beta x+x_\beta,\quad 
			x_{\alpha+\beta}x_{-\alpha}=-x_\beta x; \label{4.25}\\[4pt]
			&x_{\pm\alpha}^2=x_\beta^2=x_{\alpha+\beta}^2
			=x_\alpha x_{\alpha+\beta}=x_{\alpha+\beta}x_\alpha
			=x_\beta x_{\alpha+\beta}=x_{\alpha+\beta}x_\beta
			=x_\beta x_{-\alpha}=x_{-\alpha}x_\beta=0. \label{4.26}
		\end{align}
	\end{df}
	
	\noindent
	Relations in \eqref{4.21} resemble those of Zhu’s algebra $A(V_{A_1})$ for the rank-one lattice VOA $V_{A_1}$, 
	see \eqref{3.26} and \eqref{3.29}.  
	Relations in \eqref{4.22}–\eqref{4.24} describe the product rules between $\{x,y\}$ and 
	$\{x_\alpha,x_{-\alpha},x_\beta,x_{\alpha+\beta}\}$, while 
	relations in \eqref{4.25}–\eqref{4.26} encode the products among $\{x_\alpha,x_{-\alpha},x_\beta,x_{\alpha+\beta}\}$.
	
	Observe that relations \eqref{4.21}--\eqref{4.26} are contained among the relations 
	\eqref{rel1}--\eqref{rel4} for the Zhu algebra $A(V_{A_2})$. 
	
	\subsubsection{Structure of $A(V_P)$}
	We are now ready to establish the main result of this section.
	
	\begin{thm}\label{thm4.10}
		There exists an isomorphism of unital associative algebras
		\begin{equation}\label{4.27}
			\begin{aligned}
				F:\ 
				A_P
				&=\C\<x,y,x_\alpha,x_{-\alpha},x_\beta,x_{\alpha+\beta}\>/R
				\;\longrightarrow\;
				A(V_P),\\[4pt]
				x &\longmapsto [\alpha(-1)\vac],\quad 
				y \longmapsto [\beta(-1)\vac],\\
				x_{\pm\alpha} &\longmapsto [e^{\pm\alpha}],\quad 
				x_\beta \longmapsto [e^{\beta}],\quad 
				x_{\alpha+\beta} \longmapsto [e^{\alpha+\beta}],
			\end{aligned}
		\end{equation}
		where the same symbols $x,y,x_\alpha,x_{-\alpha},x_\beta,x_{\alpha+\beta}$ 
		are used to denote their equivalence classes in $A_P$.
	\end{thm}
	
	\begin{proof}
		First, we show that $F$ is well-defined, i.e., $F$ preserves the relations given by \eqref{4.21}--\eqref{4.26}. Indeed, by \eqref{3.26} and the fact that there is an algebra homomorphism $A(V_{\Z\al})\to A(V_P)$, $F$ preserves \eqref{4.21}. Note that the following relations hold in $A(V_P)$:
		\begin{align*}
			&[\al(-1)\vac]\ast[\b(-1)\vac]=[\b(-1)\al(-1)\vac]=[\al(-1)\b(-1)\vac]=[\b(-1)\vac]\ast[\al(-1)\vac],\\
			& [\b(-1)\vac] \ast [e^{\pm \al}]-[\b(-1)\vac]\ast[e^{\pm \al}]=[\b(0)e^{\pm \al}]=\mp [e^{\pm\al}],\\
			& [e^\b]\ast [\b(-1)\vac]=[\b(-1)e^\b]=-[e^\b],\quad [\b(-1)\vac]\ast [e^\b]=[(\b(0)+\b(-1))e^\b]=[e^\b],
		\end{align*}
		where the last equality follows from $\b(-1)e^\b+e^\b\in O=O(V_P)$ by \eqref{4.7} and Proposition~\ref{prop4.7}. Hence $F$ preserves \eqref{4.22}. Similarly, one can show that $F$ preserves \eqref{4.23}. The preservation of \eqref{4.24} under $F$ follows from 
		\begin{align*}
			[[\al(-1)\vac], [e^\b]]=[\al(0)e^\b]=-[e^\b],\quad [[\al(-1)\vac], [e^{\al+\b}]]=[\al(0)e^{\al+\b}]=[e^{\al+\b}];
		\end{align*} 
		and the preservation of \eqref{4.25} under $F$ follows from 
		\begin{align*}
			&[e^\al]\ast[e^\b]=[e^\b_{-1}e^\al]=[\epsilon(\b,\al)\b(-1)e^{\al+\b}]=-[e^{\al+\b}]\ast [\b(-1)\vac],\\
			&[e^\b]\ast [e^\al]=[\epsilon(\al,\b)\al(-1)e^{\al+\b}]=-[e^{\al+\b}]\ast [\b(-1)\vac]-[e^{\al+\b}],\\
			& [e^{-\al}]\ast [e^{\al+\b}]=[\epsilon(\al+\b,-\al)(\al+\b)(-1)e^{\b}]=-[e^\b]\ast[\al(-1)\vac]+[e^\b],\\
			& [e^{\al+\b}]\ast [e^{-\al}]=[\epsilon(-\al,\al+\b)(-\al(-1))e^\b]=-[e^\b]\ast[\al(-1)\vac],
		\end{align*}
		where we used the fact that $[\b(-1)e^\b]=-[e^\b]$ in $A(V_P)$. 
		
		Finally, for $\ga,\ga'\in \{ \al,-\al,\b,\al+\b\}$ such that $\ga+\ga'\notin \{0,\al,-\al,\b,\al+\b\}$, by \eqref{4.7} and Proposition~\ref{prop4.7}, we have 
		$e^\ga\ast e^{\ga'}\in M_{\hat{\h}}(1,\ga+\ga')\subset O$, and hence $[e^\ga]\ast[e^{\ga'}]=0$ in $A(V_P)$. This shows that $F$ preserves \eqref{4.26}, and therefore $F$ is well-defined.
		
		By \eqref{4.20}, it is clear that $F$ is surjective. To show that $F$ is an isomorphism, we construct its inverse. Consider the following linear map:
		\begin{equation}\label{4.28}
			\bar{(\cdot)}: \h=\C\al \op \C \b\longrightarrow A_P, \quad h=\la\al+\mu \b\longmapsto \bar{h}=\la x+\mu y,\quad \la,\mu\in\C.
		\end{equation}
		Again, we use the same symbols $x$ and $y$ to denote their images in $A_P$. Now define
		\begin{equation}\label{P-Inverse}
			\begin{aligned}
				G: V_P &\longrightarrow A_P=\C\<x,y,x_\al,x_{-\al},x_\b,x_{\al+\b}\>/R,\\
				h^1(-n_1-1)\cdots h^r(-n_r-1)e^{\ga}&\longmapsto (-1)^{n_1+\cdots+n_r} x_{\ga}\cdot \overline{h^1}\cdot \overline{h^2}\cdots \overline{h^r},\quad \ga\in \{ \al,-\al,\b,\al+\b\},\\
				h^1(-n_1-1)\cdots h^r(-n_r-1)\vac&\longmapsto (-1)^{n_1+\cdots+n_r}  \overline{h^1}\cdot \overline{h^2}\cdots \overline{h^r},\\
				M_{\hat{\h}}(1,m\al+n\b)&\longmapsto 0,\quad m\al+n\b\in P\setminus\{ 0,\al,-\al, \b,\al+\b\},
			\end{aligned}
		\end{equation}
		where $r\geq 0$, $n_1\geq \cdots \geq n_r\geq 0$, and $\overline{h^i}$ denotes the image of $h^i\in \h$ under $\bar{(\cdot)}$ in \eqref{4.28} for all $i$. 
		
		Next, we show that $G$ vanishes on $O(V_P)=O$ given by \eqref{4.7}. Clearly, $G(h(-n-2)u+h(-n-1)u)=0$ for any $h\in \h$, $u\in V_P$, and $n\geq 0$.
		
		To show that $G(\ga(-1)v+v)=0$, where $v=h^1(-n_1-1)\cdots h^r(-n_r-1)e^{\ga}\in M_{\hat{\h}}(1,\ga)$, note that 
		\begin{align*}
			&G(\ga(-1)h^1(-n_1-1)\cdots h^r(-n_r-1)e^{\ga}+h^1(-n_1-1)\cdots h^r(-n_r-1)e^{\ga})\\
			&=(-1)^{n_1+\cdots +n_r}x_\ga \cdot\overline{\ga}\cdot \overline{h^1}\cdots \overline{h^r}+(-1)^{n_1+\cdots +n_r}x_\ga \cdot \overline{h^1}\cdots \overline{h^r}\\
			&=(x_\ga\cdot \overline{\ga}+x_\ga)(-1)^{n_1+\cdots +n_r}\overline{h^1}\cdots \overline{h^r}\\
			&=0,
		\end{align*} 
		since $x_\al x+x_\al=0$, $x_{-\al}x-x_{-\al}=0$, $x_\b y+x_\b=0$, and $x_{\al+\b}(x+y)+x_{\al+\b}=0$, by \eqref{4.21}, \eqref{4.22}, \eqref{4.23}, and \eqref{4.28}. 
		
		To show $G(\ga(-1)^2v+\ga(-1)v)=0$, where $(\ga,\ga')\in \{ (\al,\b),(\b,\al),(-\al,\al+\b),(\al+\b,-\al)\}$ as in \eqref{4.8}, we claim that the following relations hold in $A_P$:
		\begin{align}
			&x_{\al+\b}x^2+x_{\al+\b}x=0\quad \text{and}\quad x_{\al+\b}y^2+x_{\al+\b}y=0,\label{4.30}\\
			&x_\b(x+y)^2+x_\b(x+y)=0\quad \text{and}\quad x_\b x^2-x_\b x=0.\label{4.31}
		\end{align}
		Indeed, by \eqref{4.23}, \eqref{4.25}, \eqref{4.24}, and \eqref{4.21}, we have 
		$$(x_{\al+\b}x+x_{\al+\b})x=-x_{\al+\b}yx=x_\al x_\b x=x_\al x x_\b+x_\al x_\b=-x_\al x_\b+x_\al x_\b=0.$$
		The second equality in \eqref{4.30} can be proved similarly, and we omit the details. Furthermore, since $x_\b y=-x_\b$ by \eqref{4.22} and $x_\b x=xx_\b+x_\b$ by \eqref{4.24}, we have 
		\begin{align*}
			x_\b(x+y)^2+x_\b(x+y)&=x_\b x^2+x_\b xy+x_\b yx+x_\b y^2+x_\b x+x_\b y\\
			&=x_\b x^2+(-x_\b)x-x_\b x+(-1)^2 x_\b +x_\b x-x_\b\\
			&=(x_\b x-x_\b)x=xx_\b x.
		\end{align*}
		On the other hand, since $x_{\al+\b}x_{-\al}=-x_\b x$ in \eqref{4.25} and $xx_{\al+\b}=x_{\al+\b}x+x_{\al+\b}$ in \eqref{4.24}, we have 
		$$xx_\b x=-xx_{\al+\b}x_{-\al}=-x_{\al+\b}(xx_{-\al})-x_{\al+\b}x_{-\al}=-x_{\al+\b}(-x_{-\al})-x_{\al+\b}x_{-\al}=0.$$
		This proves both equalities in \eqref{4.31}. 
		
		Now let $(\ga,\ga')$ be any ordered pair in \eqref{4.8}. By \eqref{4.30}, \eqref{4.31}, and \eqref{4.28}, we have $x_{\ga+\ga'}\overline{\ga}^2+x_{\ga+\ga'}\overline{\ga}=0$ in $A_P$. Thus, for $v=h^1(-n_1-1)\cdots h^r(-n_r-1)e^{\ga+\ga'}\in M_{\hat{\h}}(1,\ga+\ga')$, we have 
		\begin{align*}
			&G(\ga(-1)^2h^1(-n_1-1)\cdots h^r(-n_r-1)e^{\ga+\ga'}+\ga(-1)h^1(-n_1-1)\cdots h^r(-n_r-1)e^{\ga+\ga'})\\
			&=(x_{\ga+\ga'}\overline{\ga}^2+x_{\ga+\ga'}\overline{\ga})(-1)^{n_1+\cdots+n_r} \overline{h^1}\cdots \overline{h^r}\\
			&=0.
		\end{align*}
		Furthermore, by definition \eqref{P-Inverse}, we have $G(M_{\hat{\h}}(1,m\al+n\b))=0$ for $m\al+n\b\in P\setminus\{ 0,\al,-\al, \b,\al+\b\}$. Finally, for $w=h^1(-n_1-1)\cdots h^r(-n_r-1)\vac\in M_{\hat{\h}}(1,0)$, by \eqref{4.22} we have 
		\begin{align*}
			&G(\al(-1)^3 h^1(-n_1-1)\cdots h^r(-n_r-1)\vac-\al(-1)h^1(-n_1-1)\cdots h^r(-n_r-1)\vac)\\
			&=(-1)^{n_1+\cdots+n_r} \overline{h^1}\cdots \overline{h^r} \cdot x^3- (-1)^{n_1+\cdots+n_r} \overline{h^1}\cdots \overline{h^r} \cdot x\\
			&=0.
		\end{align*}
		This shows that $G(O(V_P))=0$, and so $G$ induces a well-defined map $G:A(V_P)=V_P/O(V_P)\to A_P$. It is easy to see from \eqref{P-Inverse} and \eqref{4.27} that $G:A(V_P)\to A_P$ and $F: A_P\to A(V_P)$ are mutually inverse on the generators of $A_P$ and $A(V_P)$. Hence $G$ is an inverse of $F$, and $F$ is an isomorphism of associative algebras.
	\end{proof}
	
	\subsubsection{$A(V_P)$ and the skew-polynomial algebra}
	With \eqref{4.20}, Theorem~\ref{thm4.10}, and the relations \eqref{4.21}--\eqref{4.26}, we obtain the following direct sum decomposition of the Zhu algebra $A(V_P)$: 
	\begin{equation}\label{4.32}
		\begin{aligned}
			A(V_P)=A_P=&\left(\bigoplus_{n=0}^\infty \C(x_{-\al}y^n)\right)\op \C[x,y]/\<x^3-x\>\op \left(\bigoplus_{n=0}^\infty \C (x_{\al}y^n) \right)\\
			&\op \left(\C x_\b\op \C x_\b x\op\C x_\b x^2\right)\op \left(\C x_{\al+\b} \op \C x_{\al+\b}x\op\C x_{\al+\b}x^2\right),
		\end{aligned}
	\end{equation} 
	where the products among the spanning elements are given by \eqref{4.21}--\eqref{4.26}. Using the decomposition \eqref{4.32}, we have the following corollary, which will be used in the next section. 
	
	\begin{coro}\label{coro4.11}
		Let 
		\[
		J=\left(\C x_\b\op \C x_\b x\op\C x_\b x^2\right)\op \left(\C x_{\al+\b} \op \C x_{\al+\b}x\op\C x_{\al+\b}x^2\right).
		\]
		Then $J$ is a two-sided ideal of $A(V_P)$ satisfying $J^2=0$. Moreover, the quotient algebra $A^P=A(V_P)/J$ is isomorphic to $A(V_{A_1})\otimes_{\C} \C[y]$ as vector spaces, where $A_1=\Z\al$, and 
		\[
		A(V_P)=A^P\op J.
		\]
		Furthermore, both $A(V_{A_1})=A(V_{A_1})\otimes \C 1$ and $\C[y]=\C [\vac]\otimes \C[y]$ are subalgebras of the associative algebras $A(V_P)$ and $A^P$. 
	\end{coro}
	
	\begin{proof}
		Let $J_1=\C x_\b\op \C x_\b x\op\C x_\b x^2$ and $J_2=\C x_{\al+\b} \op \C x_{\al+\b}x\op\C x_{\al+\b}x^2$. Then $J=J_1\op J_2$. From \eqref{4.21}--\eqref{4.26}, it is straightforward to verify that $J_1$ and $J_2$ satisfy
		\begin{align*}
			&	xJ_1,\ yJ_1,\ J_1x,\ J_1y \subset J_1,\\
			&	xJ_2,\ yJ_2,\ J_2x,\ J_2y \subset J_2,\\
			&	x_\al J_1,\ J_1 x_\al\subset J_2,\quad x_\ga J_1=J_1x_\ga=0,\quad \ga\in\{-\al,\b,\al+\b\},\\
			&	x_{-\al} J_2,\ J_2 x_{-\al}\subset J_1,\quad x_{\ga'} J_2=J_2x_{\ga'}=0,\quad \ga'\in\{\al,\b,\al+\b\}.
		\end{align*}
		Hence $J=J_1\op J_2$ is a two-sided ideal of $A(V_P)$. Moreover, using $x_{\pm\al}^2=x_\b^2=x_{\al+\b}^2=0$ together with \eqref{4.21} and \eqref{4.24}, we see that $J^2=0$. From the decomposition \eqref{4.32}, we have 
		\begin{align*}
			A(V_P)/J&=\left(\bigoplus_{n=0}^\infty \C(x_{-\al}y^n)\right)\op \C[x,y]/\<x^3-x\>\op \left(\bigoplus_{n=0}^\infty \C (x_{\al}y^n) \right)\\
			&\cong \left( \C x_{-\al}\op \C[x]/\<x^3-x\>\op \C x_\al\right)\otimes_\C \C [y]
		\end{align*}
		as vector spaces. Hence $A(V_P)=A^P\op J$. 
		
		By \eqref{4.21} and \eqref{4.22}, the subspace $\C x_{-\al}\op \C[x]/\<x^3-x\>\op \C x_\al$ is closed under multiplication in $A(V_P)$. Furthermore, the product relations among $x_{ \al}$, $x_{-\al}$, and $x$ coincide with those among $[e^\al]$, $[e^{-\al}]$, and $[\al(-1)\vac]$ in $A(V_{\Z \al})$ (see Section~\ref{sec3.2} and \eqref{3.29}). Thus, the subalgebra $\C x_{-\al}\op \C[x]/\<x^3-x\>\op \C x_\al$ is isomorphic to $A(V_{A_1})$. 
	\end{proof}
	
	In fact, the subalgebra $A^P=A(V_P)/J$ of the Zhu algebra $A(V_P)$ is a skew-polynomial ring over $A(V_{A_1})$. Recall the following definition from \cite{GW04}: 
	
	\begin{df}\label{df4.13}
		Let $R$ be a (not necessarily commutative) ring, $\sigma:R\to R$ a homomorphism, and $\delta:R\to R$ a $\sigma$-derivation, that is, an additive map satisfying $\delta(ab)=\delta(a)b+\sigma(a)\delta(b)$ for all $a,b\in R$. 
		
		Then the {\bf skew-polynomial ring}, or {\bf Ore extension}, $R[x; \sigma;\delta]$ \cite{O33} is the free left $R$-module with basis $\{ 1,x,x^2,x^3,\dots\}$, where the multiplication is determined by 
		\[
		xa=\sigma(a)x+\delta(a), \quad \text{for all } a\in R. 
		\]
	\end{df} 
	
	\begin{lm}\label{lm4.14}
		Let $R$ be the subalgebra $A(V_{A_1})=\spn\{ 1,x_{\al},x_{-\al}, x,x^2\}\subseteq A(V_P)$, and let $\sigma=\Id_R$. Then the derivation $\delta:=[y,\cdot]: A(V_P)\to A(V_P)$ preserves $R$ and satisfies 
		\begin{equation}\label{4.33}
			\delta(1)=\delta(x)=\delta(x^2)=0,\quad \delta(x_\al)=-x_{\al},\quad \text{and}\quad \delta(x_{-\al})=x_{-\al}.
		\end{equation}
		In particular, $\delta$ restricts to a $\sigma$-derivation of $R$. 
	\end{lm}
	
	\begin{proof}
		Since $xy=yx$, we have $[y,1]=[y,x]=[y,x^2]=0$. Moreover, by \eqref{4.22}, $yx_{\pm \al}=x_{\pm \al}y\mp x_{\pm \al}$, hence $[y,x_{\al}]=-x_{\al}$ and $[y,x_{-\al}]=x_{-\al}$. This shows that $\delta=[y,\cdot]$ preserves $R$ and satisfies \eqref{4.33}. Since $\delta(ab)=\delta(a)b+\Id(a)\delta(b)$ for all $a,b\in R$, $\delta$ is a $\sigma=\Id$-derivation. 
	\end{proof}
	
	\begin{coro}\label{coro4.15}
		The quotient $A^P=A(V_P)/J$ is isomorphic to the skew-polynomial algebra $A(V_{A_1})[y; \Id;\delta]$, where $\delta=[y,\cdot]|_{A(V_{A_1})}$. 
	\end{coro}
	
	\begin{proof}
		By Corollary~\ref{coro4.11}, we have $A^P\cong A(V_{A_1})\otimes_{\C} \C[y]=A(V_{A_1})[y]$ as vector spaces. By Lemma~\ref{lm4.14}, the derivation $\delta=[y,\cdot]|_{A(V_{A_1})}$ is a $\Id$-derivation on $A(V_{A_1})$, satisfying $ya=\Id(a)y+\delta(a)$ for all $a\in A(V_{A_1})$. Hence $A^P\cong A(V_{A_1})[y; \Id;\delta]$, in view of Definition~\ref{df4.13}. 
	\end{proof}

	\section{Representation of the rank-two parabolic-type VOA $V_P$}	\label{sec:7}
	
	In this Section, we use our main results in the Section~\ref{sec4} to classify the irreducible modules over the parabolic-type VOA $V_P$.
	
	\subsection{Construction of irreducible modules of $V_P$}
	
	Note that $P=\Z\al\op \Z_{\geq 0}\b$ is also an abelian semigroup. Let $I\leq P$ be sub-semigroup $\Z \al\op \Z_{>0}\b$. In Figure~\ref{fig5}, the dots represent elements in $P$, and the red dots represent the elements in $I$. 
	\begin{figure}
		\centering
		\begin{tikzpicture}
			\coordinate (Origin)   at (0,0);
			\coordinate (XAxisMin) at (-5,0);
			\coordinate (XAxisMax) at (5,0);
			\coordinate (YAxisMin) at (0,0);
			\coordinate (YAxisMax) at (0,5);
			
			\clip (-5.1,-0.7) rectangle (5.1,3.7); 
			\begin{scope} 
				\pgftransformcm{1}{0}{1/2}{sqrt(3)/2}{\pgfpoint{0cm}{0cm}} 
				\coordinate (Bone) at (0,2);
				\coordinate (Btwo) at (2,-2);
				\draw[style=help lines,dashed] (-7,-6) grid[step=2cm] (6,6);
				\foreach \x in {-4,-3,...,4}{
					\foreach \y in {1,...,4}{
						\coordinate (Dot\x\y) at (2*\x,2*\y);
						\node[draw,red,circle,inner sep=2pt,fill] at (Dot\x\y) {};
					}
				}
				\foreach \x in {-4,-3,...,4}{
					\foreach \y in {0}{
						\coordinate (Dot\x\y) at (2*\x,2*\y);
						\node[draw,circle,inner sep=2pt,fill] at (Dot\x\y) {};
					}
				}
				
				\draw [thick,-latex,red] (Origin) 
				-- (Bone) node [above right]  {$\al+\b$};
				\draw [thick,-latex] (Origin) node [xshift=0cm, yshift=-0.3cm] {$0$} 
				-- ($(Bone)+(Btwo)$) node [below right] {$\al$};
				\draw [thick,-latex] (Origin)
				-- ($-1*(Bone)-1*(Btwo)$) node [below left] {$-\al$};
				\draw [thick,-latex,red] (Origin)
				-- ($-1*(Btwo)$) coordinate (B3) node [above left] {$\b$};
			\end{scope} 
			\begin{scope}
				\pgftransformcm{1}{0}{-1/2}{sqrt(3)/2}{\pgfpoint{0cm}{0cm}} 
				\draw[style=help lines,dashed] (-6,-6) grid[step=2cm] (6,6);
			\end{scope}
		\end{tikzpicture}
		\caption{ \label{fig5}}
	\end{figure}

	\begin{lm}\label{lm5.1}
		The subspace $V_I=\bigoplus_{\ga\in I} M_{\hat{\h}}(1,\ga)$ is an ideal of the parabolic-type VOA $V_P$. The quotient VOA $V_P/V_I\cong \bigoplus_{n\in \Z} M_{\hat{\h}}(1,n\al)$ is a subVOA of $V_P$, and 
		\[
		V_P=(V_P/V_I)\op V_I,
		\]
		where $V_P/V_I$ and $V_I$ are both closed for the Jacobi identity of VOAs \eqref{eq:formalJacobi}. 
	\end{lm}
	
	\begin{proof}
		It follows from Figure~\ref{fig5} that $P+I\subseteq I$ and $P=I\op \Z \al$ as abelian semigroups. Then, by \cite[Proposition~3.2]{Liu24}, the subspace $V_I$ is an ideal of the VOA $V_P$. Furthermore, since 
		\[
		Y(M_{\hat{\h}}(1,n\al),z)M_{\hat{\h}}(1,m\al)\subset M_{\hat{\h}}(1,(m+n)\al)((z))
		\]
		for any $m,n\in \Z$, and $M_{\hat{\h}}(1,0)\subset V_P/V_I$, it follows that $V_P/V_I$ is a subVOA of $V_P$, with the same Virasoro element. 
	\end{proof}
	
	\begin{remark}\label{rem5.2}
		The quotient VOA $V_P/V_I$ also admits the following identification as a vector space:
		\begin{equation}\label{5.1}
			\begin{aligned}
				&V_P/V_I\cong \bigoplus_{n\in \Z} M_{\hat{\h}}(1,n\al)
				=\bigoplus_{n\in \Z} M_{\widehat{\C\al}}(1,n\al)\otimes  M_{\widehat{\C\b}}(1,0)
				=V_{\Z\al}\otimes M_{\widehat{\C\b}}(1, 0),\\
				&\al(-n_1)\cdots \al(-n_k)\b(-m_1)\cdots \b(-m_l)e^{n\al}
				\longmapsto \al(-n_1)\cdots \al(-n_k)e^{n\al}\otimes \b(-m_1)\cdots \b(-m_l)\vac.
			\end{aligned} 
		\end{equation}
		However, the identification in \eqref{5.1} is {\bf not} an isomorphism of VOAs between $V_P/V_I$ and the tensor product $V_{\Z\al}\otimes M_{\widehat{\C \b}}(1,0)$ as defined in \cite{FHL}. This is because the operators $E^{+}(-\al,z)$ and $\b(-n)$ (for $n\geq 1$) do not commute when $(\al|\b)\neq 0$. 
		
		On the other hand, it follows directly from the spanning elements \eqref{4.7} and Theorem~\ref{prop4.7} that the Zhu algebra 
		\[
		A(V_P/V_I)=A\!\left(\bigoplus_{n\in \Z} M_{\hat{\h}}(1,n\al)\right)
		\]
		is isomorphic to the skew-polynomial algebra $A^P=A_P/J=A(V_{\Z\al})[y;\Id;\delta]$ as described in Corollaries~\ref{coro4.11} and~\ref{coro4.15}.
	\end{remark}
	
	Note that the rank-one lattice VOA $V_{A_1}=V_{\Z\al}$ is clearly a subVOA of both $V_P$ and $V_P/V_I$ (see \eqref{5.1}). By Theorem~3.1 in \cite{D}, $V_{A_1}$ has two irreducible modules: $V_{\Z\al}$ and $V_{\Z\al+\frac{1}{2}\al}$. We will use these irreducible $V_{A_1}$-modules to construct irreducible $V_P$-modules. 
	
	\subsubsection{Construction of $L^{(0,\la)}$ and $L^{(\frac{1}{2}\al,\la)}$}
	
	Let $A_2=\Z\al\op \Z\b$ be the root lattice of type $A_2$, and recall that $\h=\C\otimes_\Z A_2$ is equipped with a nondegenerate symmetric bilinear form $(\cdot|\cdot):\h\times \h\to \C$.  
	
	\begin{df} \label{df5.3}
		Let $\la\in (\C\al)^\perp\subset \h$. Define $L^{(0,\la)}$ and $L^{(\frac{1}{2}\al,\la)}$ to be the following vector spaces: 
		\begin{align}
			L^{(0,\la)}&:=\bigoplus_{n\in \Z} M_{\hat{\h}}(1,n\al)\otimes \C e^\la
			\cong V_{\Z\al}\otimes M_{\widehat{\C\b}}(1,\la), \label{5.4}\\
			L^{(\frac{1}{2}\al,\la)}&:=\bigoplus_{n\in \Z} M_{\hat{\h}}(1,n\al+\tfrac{1}{2}\al)\otimes \C e^\la
			\cong  V_{\Z\al+\frac{1}{2}\al}\otimes  M_{\widehat{\C\b}}(1,\la),\label{5.5}
		\end{align}
		where $M_{\hat{\h}}(1,n\al)$ and $M_{\hat{\h}}(1,n\al+\frac{1}{2}\al)$ are level-$1$ modules for the Heisenberg Lie algebra $\hat{\h}$.
		
		Define the actions of the operators $h(m)$, $e_{n\al}$, and $z^{n\al}$ for $m,n\in \Z$ on the tensor product spaces $\bigoplus_{n\in \Z} M_{\hat{\h}}(1,n\al)\otimes \C e^\la$ and $\bigoplus_{n\in \Z} M_{\hat{\h}}(1,n\al+\frac{1}{2}\al)\otimes \C e^\la$ by
		\begin{align}
			h(0)&:=h(0)\otimes \Id+\Id \otimes h(0),\qquad 
			h(m):=h(m)\otimes \Id,\quad m\neq 0, \label{5.6}\\
			e_{n\al}&:=e_{n\al}\otimes \Id,\qquad 
			z^{n\al}:=z^{n\al}\otimes \Id,\quad n\in \Z, \label{5.7}
		\end{align}
		where $h(0)e^\la:=(\la|h)e^\la$. In particular, for any $m,n\in \Z$ and $n_1\geq \cdots \geq n_r\geq 1$, we have 
		\begin{align}
			h(m)(h^1(-n_1)\cdots h^r(-n_r)e^{n\al}\otimes e^\la)
			&:=h(m)h^1(-n_1)\cdots h^r(-n_r)e^{n\al}\otimes e^\la,\quad m\neq 0,\label{5.8}\\
			h(0)(e^{m\al}\otimes e^\la)
			&=(m\al+\la|h)e^{m\al}\otimes e^\la,\label{5.9a}\\
			h(0)(e^{m\al+\frac{1}{2}\al}\otimes e^\la)
			&=(m\al+\tfrac{1}{2}\al+\la|h)e^{m\al+\frac{1}{2}\al}\otimes e^\la,\label{5.9}\\
			e_{n\al}(e^{m\al}\otimes e^\la)
			&=\epsilon(n\al,m\al)e^{(m+n)\al}\otimes e^\la,\label{5.10a}\\
			e_{n\al}(e^{m\al+\frac{1}{2}\al}\otimes e^\la)
			&=\epsilon(n\al,m\al)e^{(m+n)\al+\frac{1}{2}\al}\otimes e^\la,\label{5.10}\\
			z^{n\al}(e^{m\al}\otimes e^\la)
			&=z^{(n\al|m\al)}e^{m\al}\otimes e^\la,\label{5.11a}\\
			z^{n\al}(e^{m\al+\frac{1}{2}\al}\otimes e^\la)
			&=z^{(n\al|m\al+\frac{1}{2}\al)}e^{m\al+\frac{1}{2}\al}\otimes e^\la. \label{5.11}
		\end{align}
		
		Define the module vertex operators $Y_M: V_P\to \End(L^{(\epsilon,\la)})[[z,z^{-1}]]$, where $\epsilon=0$ or $\frac{1}{2}\al$, by
		\begin{align}
			Y_M(h^1(-n_1-1)\cdots h^r(-n_r-1)e^{n\al}, z)
			&:={\tiny\begin{matrix}\circ \\\circ\end{matrix}}
			(\partial_{z}^{(n_1)}h^1(z))\cdots (\partial_{z}^{(n_r)}h^r(z))
			Y(e^{n\al},z)
			{\tiny\begin{matrix}\circ \\\circ\end{matrix}},
			\quad n\in \Z,\label{5.12}\\
			Y_M(h^1(-n_1-1)\cdots h^r(-n_r-1)e^{\ga}, z)
			&:=0,\quad \ga\in I=\Z\al\op \Z_{>0}\b, \label{5.13}
		\end{align}
		where $h^i\in \h$ for all $i$, $n_1\geq \cdots \geq n_r\geq 0$, and 
		\[
		Y(e^{n\al},z)=E^-(-n\al,z)E^+(-n\al,z)e_{n\al}z^{n\al}.
		\]
	\end{df}
	
	By \eqref{5.8} and \eqref{5.9}, it is clear that the subspace $M_{\hat{\h}}(1,n\al+\epsilon)\otimes \C e^\la\subset L^{(\epsilon,\la)}$ is a $\hat{\h}$-module, and it is isomorphic to $M_{\hat{\h}}(1,n\al+\epsilon+\la)$, where $\epsilon=0$ or $\tfrac{1}{2}\al$ and $n\in \Z$.

	\subsection{Irreducibility of $L^{(0,\la)}$ and $L^{(\frac{1}{2}\al,\la)}$}
	
	\begin{lm}\label{lm5.4}
		For any $\la\in (\C\al)^\perp\subset \h$, the vector spaces $L^{(\epsilon,\la)}$, where $\epsilon=0$ or $\frac{1}{2}\al$, equipped with $Y_M$ in Definition~\ref{df5.3}, are weak $V_P$-modules. 
	\end{lm}
	\begin{proof}
		We need to show that  the operator $Y_M$ satisfies the truncation property and Jacobi identity.
		For the truncation property, we fix a spanning element $v=h^1(-n_1)\ds h^r(-n_r)e^{n\al+\epsilon}\o e^\la$ of $L^{(\epsilon,\la)}$, where $h^i\in \h$ for all $i$, and $n_1\geq \ds \geq n_r\geq 1$, and show that $Y_M (a,z)v\in L^{(\epsilon,\la)}((z))$ for any $a\in V_P$.
		
		Indeed, if $a\in V_I$, then by \eqref{5.13} we have $a_nv=0$ for any $n\in \Z$, there is noting to prove. Now assume that $a\in M_{\hat{\h}}(1,m\al)$. If $a=e^{m\al}$, then by \eqref{5.6}, \eqref{5.7}, and \eqref{5.12}, we have 
		\begin{equation}\label{5.14}
			Y_M(a,z) (e^{n\al+\epsilon}\o e^\la)=(E^-(-m\al,z)E^+(-m\al,z)e_{m\al}z^{m\al} e^{n\al+\epsilon})\o e^\la\in L^{(\epsilon,\la)}((z)),
		\end{equation}
		since the $V_{A_1}$-modules $V_{\Z\al}$ and $V_{\Z\al+\frac{1}{2}\al}$ satisfy the truncation property. Furthermore, since $[h(-p),E^-(-m\al,z)]=0$ and $[h(-p),E^+(-m\al,z)]=-(h|m\al)z^{-p} E^+(-m\al,z)$ for any $h\in \h$ and $p>0$ (see \cite{FLM,LL}), it follows that 
		\begin{align*}
			Y_M(a,z)v&=h^1(-n_1) Y_M(a,z) (h^2(-n_2)\ds h^r(-n_r)e^{n\al+\epsilon}\o e^\la)\\
			&-(h^1|m\al) z^{-n_1}Y_M(a,z) (h^2(-n_2)\ds h^r(-n_r)e^{n\al+\epsilon}\o e^\la).
		\end{align*}
		Then by an induction on the length $r$ of $v$, with base case given by \eqref{5.14}, we have
		\begin{equation} \label{5.15}
			Y(e^{m\al},z) \left(h^1(-n_1)\ds h^r(-n_r)e^{n\al+\epsilon}\o e^\la\right)\in L^{(\epsilon,\la)}((z)). 
		\end{equation}
		
		Now let $a$ be a general spanning element $a=h_1(-m_1)\ds h_s(-m_s)e^{m\al}$ of $ M_{\hat{\h}}(1,m\al)$, where $h_j\in \h$ for all $j$, and $m_1\geq\ds \geq m_s\geq 1$. We show that $a_kv=0$ if $k\gg 0$. Again, by induction on the length $s$ of $a$, it suffices to consider the case when $s=1$. The proof of the general case is similar, we omit it. Note that 
		\begin{align*}
			&(h_1(-m_1)e^{m\al})_kv\\
			&=\sum_{j\geq 0}\binom{-m_1}{j}(-1)^j h_1(-m_1-j)(e^{m\al})_{k+j}v-\sum_{j\geq 0} \binom{-m_1}{j} (-1)^{m_1+j} (e^{m\al})_{-m_1+k-j}(h_1(j)v)\\
			&=\sum_{j\geq 0}\binom{-m_1}{j}(-1)^j h_1(-m_1-j)(e^{m\al})_{k+j}v\\
			&\ \ -(n\al+\epsilon+\la|h) (e^{m\al})_{-m_1+k} \left(h^1(-n_1)\ds h^r(-n_r)e^{n\al+\epsilon}\o e^\la\right)\\
			&\ \ -\binom{-m_1}{n_r} n_r(h_1|h^r) (-1)^{m_1+n_r} (e^{m\al})_{-m_1+k-n_r} \left(h^1(-n_1)\ds \widehat{h^r(-n_r)}e^{n\al+\epsilon}\o e^\la\right)\\
			&\ \ -\binom{-m_1}{n_{r-1}} n_{r-1}(h_1|h^{r-1})(-1)^{m_1+n_{r-1}} (e^{m\al})_{-m_1+k-n_{r-1}}\left(h^1(-n_1)\ds \widehat{h^{r-1}(-n_{r-1})}h^r(-n_r)e^{n\al+\epsilon}\o e^\la\right)\\
			&\quad \vdots\\
			&\ \ -\binom{-m_1}{n_1}n_1(h_1|h^1)(-1)^{m_1+n_1}(e^{m\al})_{-m_1+k-n_1}\left(\widehat{h^1(-n_1)}\ds h^r(-n_r)e^{n\al+\epsilon}\o e^\la\right).
		\end{align*}
		By \eqref{5.15}, it is clear that we can choose $k\gg 0$ large enough so that each term on the right hand side of the equation above is equal to $0$. This shows the truncation property of $Y_M$. 
		
		It remains to show the Jacobi identity of $Y_M$. Let $a\in M_{\hat{\h}}(1,\ga)$ and $b\in M_{\hat{\h}}(1,\eta)$, where $\ga,\eta\in P$. We need to show that 
		\begin{equation}\label{5.16}
			\begin{aligned}
				&z_0^{-1}\delta\left(\frac{z_1-z_2}{z_0}\right) Y_M(a,z_1)Y_M(b,z_2)-	z_0^{-1}\delta\left(\frac{-z_2+z_1}{z_0}\right) Y_M(b,z_2)Y_M(a,z_1)\\
				&=	z_2^{-1}\delta\left(\frac{z_1-z_0}{z_2}\right) Y_M(Y(a,z_0)b,z_2).
			\end{aligned}
		\end{equation}
		
		Note that $Y(a,z_0)b\in M_{\hat{\h}}(1,\ga+\eta)((z_0))$. If either $\ga$ or $\eta$ are contained in $I\subset P$, then by Lemma~\ref{lm5.1}, \eqref{5.13}, and the fact that $I+P=P+I\ssq I$, both sides of the Jacobi identity \eqref{5.16} are $0$. Now assume $a=h^1(-n_1-1)\ds h^r(-n_r-1)e^{n\al}$ and $b=h_1(-m_1-1)\ds h_s(-m_s-1)e^{m\al}$ for some $m,n\in \Z$, $h^i,h_j\in \h$ for all $i,j$, $n_1\geq \ds\geq n_r\geq 0$, and $m_1\geq \ds\geq m_s\geq 0$. By adopting a similar argument as the proof of Theorem 8.6.1 in \cite{FLM}, we can show that 
		$$[Y_M(a,z_1),Y_M(b,z_2)]=\Res_{z_0} z_2^{-1}Y_M(Y(a,z_0)b,z_2)e^{-z_0(\partial/\partial z_1)}\left((z_1/z_2)^{m\al} \delta(z_1/z_2)\right).$$
		This commutator relation also (essentially) follows from the fact that the $V_{A_1}$-module vertex operators for $V_{\Z\al}$ and $V_{\Z\al+\frac{1}{2}\al}$ satisfy the Jacobi identity. Then by Theorem 8.8.9 in \cite{FLM}, the Jacobi identity \eqref{5.16} holds for $Y_M$. 
	\end{proof}
	

	\begin{lm}\label{lm5.5}
		Given $\la\in (\C\al)^\perp\subset \h$, the weak $V_P$-modules $(L^{(0,\la)},Y_M)$ and $(L^{(\frac{1}{2}\al,\la)},Y_M)$ are irreducible ordinary $V_P$-modules, whose bottom degrees are $\C (\vac\otimes e^\la)$ and $\C (e^{\frac{1}{2}\al}\o e^\la)\op \C (e^{-\frac{1}{2}\al}\o e^\la)$, respectively. 
	\end{lm}
	\begin{proof} 
		Note that $\Res_{z}z Y_M(\om ,z)=L_M(0)=\frac{1}{2}\sum_{i=1}^2 \sum_{s\geq 0}u^i(-s)u^i(s)$, where $\{u^1,u^2\}$ is an orthonormal basis of $\h$. By \eqref{5.6} and the fact that $(\la|\al)=0$, we have 
		\begin{align*}
			L_M(0) (e^{n\al+\epsilon}\o e^\la)&=\frac{1}{2} (n\al+\epsilon+\la|n\al+\epsilon+\la) e^{n\al+\epsilon}\o e^\la\\
			&=\left(\frac{1}{2}(n\al+\epsilon|n\al+\epsilon)+\frac{(\la|\la)}{2}\right) e^{n\al+\epsilon}\o e^\la.
		\end{align*}
		Moreover, by \eqref{5.6} again, it is easy to show that $[L_M(0), h(-n)]=n h(-n)$, for any $h\in \h$ and $n>0$. Hence we have 
		\begin{equation}\label{5.17}
			\begin{aligned}
				&L_M(0)\left(h^1(-n_1)\ds h^r(-n_r)e^{n\al+\epsilon}\o e^\la\right)\\
				&=\left(n_1\ds+n_r+\frac{1}{2}(n\al+\epsilon|n\al+\epsilon)+\frac{(\la|\la)}{2}\right)h^1(-n_1)\ds h^r(-n_r)e^{n\al+\epsilon}\o e^\la,
			\end{aligned}
		\end{equation}
		where $\epsilon=0$ or $\frac{1}{2}\al$, $h^i\in \h$ for all $i$, $n\in \Z$, and $n_1\geq \ds \geq n_r\geq 1$. Since $(\pm \frac{1}{2}\al|\pm \frac{1}{2}\al)=\frac{1}{2}$, then it follows from \eqref{5.4} , \eqref{5.5}, and \eqref{5.17} that $L^{(0,\la)}$ and $L^{(\frac{1}{2}\al,\la)}$ are graded vector spaces, with the grading subspaces given by $L_{M}(0)$-eigenspaces: 
		\begin{equation}\label{5.18}
			L^{(0,\la)}=\bigoplus_{m=0}^\infty \left(L^{(0,\la)}\right)_{\frac{(\la|\la)}{2}+m},\quad L^{(\frac{1}{2}\al,\la)}=\bigoplus_{m=0}^\infty \left(L^{(\frac{1}{2}\al,\la)}\right)_{\frac{(\la|\la)}{2}+\frac{1}{4}+m}.
		\end{equation}
		By \eqref{5.17} and \eqref{5.18}, it is easy to see that the bottom levels $(m=0)$ of $L^{(0,\la)}$ and $L^{(\frac{1}{2}\al,\la)}$ are given by $\C (\vac\otimes e^\la)$ and $\C (e^{\frac{1}{2}\al}\o e^\la)\op \C (e^{-\frac{1}{2}\al}\o e^\la)$, respectively. 
		
		Now we show that $L^{(0,\la)}$ and $L^{(\frac{1}{2}\al,\la)}$ are irreducible. We only prove the irreducibility of $L^{(\frac{1}{2}\al,\la)}$, the other one is similar. Let $W\neq 0$ be a submodule of $L^{(\frac{1}{2}\al,\la)}$. Consider a nonzero element $0\neq u\in W$. By the decomposition~\ref{5.5}, $u$ can be written as follows:
		$$u=u_{-m}+u_{-m+1}+\ds+u_0+\ds +u_n\in \bigoplus_{n\in \Z} M_{\hat{\h}}(1,n\al+\frac{1}{2}\al)\o \C e^\la,$$
		where $u_j\in M_{\hat{\h}}(1,j\al+\frac{1}{2}\al)\o \C e^\la$ for all $-m\leq j\leq n$. By \eqref{5.6} and \eqref{5.9}, we have
		$$\b(0)u_j=\left(j\al+\frac{1}{2}\al+\la|\b\right)u_j=\left((\la|\b)-j-\frac{1}{2}\right)u_j,\quad -m\leq j\leq n.$$
		i.e., $u_j$ with $-m\leq j\leq n$ are eigenvectors of $\b(0)$ of distinct eigenvalues. Since $\b(0)^ku\in W$ for any $k\geq 0$, it follows that $u_j\in W$ for all $j$ (using the Vandermonde determinant). 
		
		Since $u\neq 0$, we may assume that $0\neq u_j\in W$ for some fixed $j$. Since  $M_{\hat{\h}}(1,j\al+\frac{1}{2}\al)\o\C e^\la $ is isomorphic to $\hat{\h}$-module $M_{\hat{\h}}(1,j\al+\frac{1}{2}\al+\la)$ by the remark after Definition~\ref{df5.3}, then by applying $h(m)$, with $h\in \h$ and $m\geq 0$, repeatedly onto $u_j$, we can show that $e^{j\al+\frac{1}{2}\al}\o e^\la\in W$. Hence
		$$e^{(j+n)\al+\frac{1}{2}\al}\o e^\la=\epsilon(j\al,n\al)^{-1} e_{n\al} \left(e^{j\al+\frac{1}{2}\al}\o e^\la\right)\in W,\quad n\in \Z,$$
		in view of \eqref{5.10}. This shows $e^{m\al+\frac{1}{2}\al}\o e^\la\in W$ for all $m\in \Z$. Now it follows from \eqref{5.8} that $M_{\hat{\h}}(1,m\al+\frac{1}{2}\al)\o \C e^\la\ssq W$ for all $m\in \Z$. Hence we have  $L^{(\frac{1}{2}\al,\la)} = W$. 
	\end{proof}

	\subsection{Classification of irreducible modules over $V_P$}\label{Sec:9.2}
	By Lemmas~\ref{lm5.4} and~\ref{lm5.5},
	\begin{equation}\label{5.19}
		\Sigma(P)
		=\left\{ (L^{(0,\la)},Y_M),\, (L^{(\frac{1}{2}\al,\la)},Y_M)
		:\, \la\in (\C\al)^\perp\subset \h \right\}
	\end{equation}
	is a set of irreducible modules over the parabolic-type subVOA $V_P$ of $V_{A_2}$, where $Y_M$ is defined in Definition~\ref{df5.3}.  
	Using the description of the Zhu’s algebra $A(V_P)=A_P$ given in Theorem~\ref{thm4.10} and Corollary~\ref{coro4.11}, we will show that $\Sigma(P)$ forms a complete list of irreducible $V_P$-modules. 
	
	By Lemma~\ref{lm5.5} and Theorem~2.1.2 in \cite{Z}, the spaces 
	\[
	U^{(0,\la)}=\C (\vac\o e^\la)
	\quad\text{and}\quad
	U^{(\frac{1}{2}\al,\la)}
	=\C (e^{\frac{1}{2}\al}\o e^\la)\op \C (e^{-\frac{1}{2}\al}\o e^\la)
	\]
	are irreducible $A(V_P)$-modules.  
	For simplicity, we use the following notations:
	\begin{align}
		&U^{(0,\la)}:=\C e,\quad\text{where } e=\vac\o e^\la, \label{eq:U0la}\\
		&U^{(\frac{1}{2}\al,\la)}:=\C e^+\op\C e^-,\quad\text{where } 
		e^+=e^{\frac{1}{2}\al}\o e^\la,\ \ e^-=e^{-\frac{1}{2}\al}\o e^\la. \label{eq:U1/2alla}
	\end{align}
	By Corollary~\ref{coro4.11} and~\eqref{4.32}, we have 
	\[
	A(V_P)=(A(V_{\Z\al})\o \C[y])\op J
	\]
	as vector spaces, where $A(V_{\Z\al})$ is a subalgebra of $A_P$, and $J$ is a two-sided nilpotent ideal of $A(V_P)$. 
	
	Since the action of $[a]\in A(V_P)$ on $U^{(\epsilon,\la)}$ is given by 
	$o(a)=\Res_z z^{\wt a-1}Y_M(a,z)$, it follows from~\eqref{5.6}--\eqref{5.11} that the spanning elements of $U^{(\epsilon,\la)}$ satisfy
	\begin{align}
		&J.e=J.e^+=J.e^-=0,\label{5.20}\\
		&x_\al.e=x_{-\al}.e=x.e=0,\quad y.e=(\la|\b)e,\label{5.21}\\
		&x_\al.e^+=0,\quad x_\al.e^-=e^+,\quad
		x_{-\al}.e^+=e^-,\quad x_{-\al}.e^-=0,\nonumber\\
		&x.e^{\pm}=\pm e^{\pm},\quad 
		y.e^{\pm}=\big((\la|\b)\mp\tfrac{1}{2}\big)e^{\pm}. \label{5.22}
	\end{align}
	By Theorem~2.2.2 in \cite{Z}, to show that $\Sigma(P)$ is a complete set of irreducible $V_P$-modules, it suffices to prove that
	\begin{equation}\label{5.23}
		\Sigma_0(P)
		=\left\{ U^{(0,\la)},\, U^{(\frac{1}{2}\al,\la)} 
		:\, \la\in (\C\al)^\perp\subset \h \right\}
	\end{equation}
	is a complete list of irreducible $A(V_P)$-modules.
	
	\begin{thm}\label{thm5.6}
		Let $U\neq 0$ be an irreducible $A(V_P)$-module. Then $U$ is isomorphic to either 
		$U^{(0,\la)}$ or $U^{(\frac{1}{2}\al,\la)}$ for some $\la\in (\C\al)^\perp$.
	\end{thm}
	
	\begin{proof}
		Since $J\subset A(V_P)$ is a nilpotent ideal, we must have $J.U=0$.  
		Hence $U$ is an irreducible module over the quotient algebra 
		\[
		A^P=A(V_P)/J
		=\left(\bigoplus_{n=0}^\infty \C(x_{-\al}y^n)\right)
		\op \C[x,y]/\<x^3-x\>
		\op \left(\bigoplus_{n=0}^\infty \C(x_{\al}y^n)\right).
		\]
		By Corollary~\ref{coro4.11}, $A^P\cong A(V_{\Z\al})\o\C[y]$ as a vector space, and $A(V_{\Z\al})$ is a subalgebra of $A^P$. Thus $U$ is also an $A(V_{\Z\al})$-module.  
		Recall that $A(V_{\Z\al})$ is semisimple with two irreducible modules (up to isomorphism):
		$W^0=\C\vac$ and $W^{\frac{1}{2}\al}=\C e^{\frac{1}{2}\al}\op \C e^{-\frac{1}{2}\al}$.
		Hence, as an $A(V_{\Z\al})$-module,
		\[
		U=\bigoplus_{i\in I} W^0 \op \bigoplus_{j\in K} W^{\frac{1}{2}\al},
		\]
		where $|I|$ and $|K|$ denote the multiplicities of $W^0$ and $W^{\frac{1}{2}\al}$, respectively.  
		
		\textbf{Case I.} $I\neq\emptyset$.  
		Then a nonzero copy of $W^0$ is contained in $U$.  
		Let $W=\C[y]\cdot W^0=\C[y]\cdot\vac$.  
		Since $x_{\pm\al}\vac=x\vac=0$, and
		\begin{equation}\label{5.24}
			xy=yx,\quad yx_\al=x_\al y-x_\al,\quad yx_{-\al}=x_{-\al}y-x_{-\al},
		\end{equation}
		it follows from~\eqref{4.22} that $W$ is an $A^P$-submodule.  
		As $U$ is irreducible, we must have $U=W=\C[y]\cdot\vac$.  
		Thus $U$ is an irreducible $\C[y]$-module, and by Hilbert’s Nullstellensatz,
		$U\cong\C[y]/\<y-\la_0\>$ for some $\la_0\in\C$.  
		Choose $\la\in\h$ with $(\la|\al)=0$ and $(\la|\b)=\la_0$; then, by~\eqref{5.21},
		$U\cong U^{(0,\la)}$.
		
		\textbf{Case II.} $K\neq\emptyset$.  
		Then a nonzero copy of $W^{\frac{1}{2}\al}$ is contained in $U$.  
		By~\eqref{5.24}, the subspace $W=\C[y]\cdot W^{\frac{1}{2}\al}\subset U$ 
		is an $A^P$-submodule, hence $U=W=\C[y]\cdot e^{\frac{1}{2}\al}\op \C[y]\cdot e^{-\frac{1}{2}\al}$.  
		We will show that $U\cong U^{(\frac{1}{2}\al,\la)}$ for some $\la\in(\C\al)^\perp$.  
		
		For simplicity, write $e^+=e^{\frac{1}{2}\al}$ and $e^-=e^{-\frac{1}{2}\al}$.  
		As in~\eqref{5.22}, we have
		\[
		x_\al.e^+=0,\quad x_\al.e^-=e^+,\quad 
		x_{-\al}.e^+=e^-,\quad x_{-\al}.e^-=0,\quad
		x.e^{\pm}=\pm e^{\pm}.
		\]
		Furthermore, from~\eqref{5.24} we obtain 
		$x_\al y^n=(y+1)^n x_\al$ and $x_{-\al}y^n=(y-1)^n x_{-\al}$ for all $n\ge0$.  
		Thus for any $f(y),g(y)\in\C[y]$,
		\begin{align}
			&x_\al\big((y-1)f(y)e^-\big)=y f(y+1)e^+,\quad
			x_\al\big(y g(y)e^+\big)=(y+1)g(y+1)e^+=0,\label{5.25}\\
			&x_{-\al}\big(y g(y)e^+\big)=(y-1)g(y-1)e^-,\quad
			x_{-\al}\big((y-1)f(y)e^-\big)=(y-2)f(y-1)x_{-\al}.e^-=0.\label{5.26}
		\end{align}
		Introduce the following subspace in $U$: 
		\[
		N:=y\C[y]\cdot e^+ + (y-1)\C[y]\cdot e^- \subset U.
		\]
		By~\eqref{5.25}, \eqref{5.26}, and the relations $xy=yx$, $x.e^{\pm}=\pm e^{\pm}$,
		$N$ is an $A^P$-submodule of $U$.  
		Since $U$ is irreducible, $N=0$ or $N=U$.
		
		If $N=0$, then $y.e^+=y.e^-=0$, and clearly $U\cong U^{(\frac{1}{2}\al,0)}$.  
		If $N=U$, then there exist $f(y),g(y)\in\C[y]$ such that
		\begin{equation}\label{5.27}
			e^+=y f(y)e^+ + (y-1)g(y)e^-.
		\end{equation}
		Applying $x_\al$ to~\eqref{5.27} and using~\eqref{5.25}, we obtain 
		$0=y g(y+1)e^+$.  
		Applying $x_{-\al}$ and using~\eqref{5.26}, we get $0=(y-1)g(y)e^-$.  
		Hence $e^+=y f(y)e^+$, i.e.
		$$
		0=(y f(y)-1)e^+=(y-\la_k)\cdots(y-\la_1)e^+,
		\quad \la_1,\dots,\la_k\in\C,
		$$
		where $\la_1,\dots,\la_k$ are nonzero.  
		Let $1\le j\le k$ be minimal such that
		\[
		(y-\la_{j-1})\cdots(y-\la_1)e^+\neq0,\quad
		(y-\la_j)\big((y-\la_{j-1})\cdots(y-\la_1)e^+\big)=0.
		\]
		Set $\tilde e^+=(y-\la_{j-1})\cdots(y-\la_1)e^+$.  
		Then $y.\tilde e^+=\la_j\tilde e^+$, 
		$x.\tilde e^+=\tilde e^+$, and $x_\al.\tilde e^+=0$.  
		Let $\tilde e^-:=x_{-\al}.\tilde e^+=(y-1-\la_{j-1})\cdots(y-1-\la_1)e^-$,  
		so that $U=A^P.\tilde e^+=\C\tilde e^+\op\C\tilde e^-$.  
		We also have $x.\tilde e^-=-\tilde e^-$ and $y.\tilde e^-=(\la_j+1)\tilde e^-$.  
		Since $\tilde e^+$ and $\tilde e^-$ are eigenvectors of $x$ with distinct eigenvalues,
		they are linearly independent.  
		Choosing $\la\in\h$ with $(\la|\al)=0$ and $(\la|\b)-\frac{1}{2}=\la_j$, we obtain 
		$U=\C\tilde e^+\op\C\tilde e^-\cong U^{(\frac{1}{2}\al,\la)}$ as $A^P$-modules.  
		By~\eqref{5.20} and Corollary~\ref{coro4.11}, they are also isomorphic as $A(V_P)=A^P\op J$-modules.
	\end{proof}
	
	\begin{coro}\label{coro5.7}
		The set
		\[
		\Sigma(P)
		=\left\{ (L^{(0,\la)},Y_M),\, (L^{(\frac{1}{2}\al,\la)},Y_M)
		:\, \la\in (\C\al)^\perp\subset \h \right\}
		\]
		is a complete list of irreducible modules over the rank-two parabolic-type subVOA $V_P$ of $V_{A_2}$.
	\end{coro}
	
	\begin{remark}
		An alternative approach to the classification of irreducible modules over parabolic-type subVOAs $V_P$ associated with general rank-two lattices VOAs $V_L$ is given in~\cite{LS25}.
	\end{remark}	
	
	\section{Finite inductions for the VOA embedding $V_P\hookrightarrow V_{A_2}$}	\label{sec:8}
	Using the structural theorem of Zhu algebra $A(V_P)$ and the classification theorem of irreducible modules over $V_P$ in the previous sections, we determine the finite induction of irreducible $V_P$-modules under the VOA embedding $V_P\hr V_{A_2}$. 
	
	\subsection{Relations in $\ker (\pi)$}
	For the exact sequence of associative algebras \eqref{eq:seqforA}:
	\begin{equation}\label{eq:8.1}
		\begin{tikzcd}
			0\arrow[r]& \ker(\pi)\arrow[r]& A(V_P)\arrow[r,"\pi"]& A(V_{A_2}),
		\end{tikzcd}
	\end{equation}
	we first describe $\ker(\pi)$. 
	
	By Definition~\ref{def:relations} and Proposition~\ref{prop:presentationofAVA2}, the associative algebras $A(V_P)$ and $A(V_{A_2})$ admit the following presentations:
	\[
	A(V_P)= \C\<x,y,x_{\pm \al},x_\b,x_{\al+\b}\>/R, 
	\qquad 
	A(V_{A_2})=\C\<x,y,x_{\pm \al},x_{\pm \b},x_{\pm(\al+\b)}\>/I,
	\]
	where $R$ and $I$ are the two-sided ideals generated by the relations \eqref{4.21}--\eqref{4.26} and \eqref{rel1}--\eqref{eq:morerel3}, respectively.  
	The sequence \eqref{eq:8.1} extends to the following commutative diagram:
	\[
	\begin{tikzcd}
		&& \C\<x,y,x_{\pm \al},x_\b,x_{\al+\b}\>\arrow[r,"\tilde{\pi}",hook]\arrow[d,two heads]&   
		\C\<x,y,x_{\pm \al},x_{\pm \b},x_{\pm(\al+\b)}\>\arrow[d,two heads]\\
		0\arrow[r]& \ker(\pi)\arrow[r]& A(V_P)\arrow[r,"\pi"]& A(V_{A_2}),
	\end{tikzcd}
	\]
	where $\tilde{\pi}$ denotes the canonical embedding of free associative algebras. Then
	\[
	\ker(\pi)= \C\<x,y,x_{\pm \al},x_\b,x_{\al+\b}\>\cap I/R.
	\]
	Comparing the relations \eqref{rel1}--\eqref{eq:morerel3} with \eqref{4.21}--\eqref{4.26}, it is straightforward to see that $\ker(\pi)$ is the two-sided ideal generated by the following additional relations, in addition to those in \eqref{4.21}--\eqref{4.26}:
	\begin{equation}\label{eq:extrarelations}
		\begin{cases}
			x_{\pm \al}(x+y)^2\pm x_{\pm \al}(x+y)=0,\\
			x_{\pm \al} y^2\mp x_{\pm \al} y=0,\\
			y^3-y=0,\\
			(x+y)^3-(x+y)=0,
		\end{cases}
	\end{equation}
	where we use the same notation for the equivalence classes of elements in $\C\<x,y,x_{\pm \al},x_{\pm \b},x_{\pm(\al+\b)}\>$ under the quotient $A(V_P)$.
	
	\begin{remark}
		Under the VOA isomorphism $L_{\widehat{\sl_3}}(1,0)\cong V_{A_2}$, the equivalence classes of the generators of $A(V_{A_2})$ are given by
		\[
		x = [\al(-1)\vac],\quad 
		y = [\b(-1)\vac],\quad 
		x_{\pm \al} = [e^{\pm \al}],\quad 
		x_{\pm \b} = [e^{\pm \b}],\quad 
		x_{\pm (\al+\b)} = [e^{\pm (\al+\b)}].
		\]
		The relations in \eqref{eq:extrarelations} can also be obtained by computing the following elements in $O(V_{A_2})$:
		\[
		e^{\al+\b}\circ e^{-\b},\quad 
		e^{-\al-\b}\circ e^{\b},\quad 
		e^{-\b}\circ e^{\al+\b},\quad 
		e^{\b}\circ e^{-\al-\b},\quad 
		e^{\b}\circ e^{-\b},\quad 
		e^{\al+\b}\circ e^{-\al-\b}.
		\]
		Note that these are elements of $O(V_{A_2})$ that are not in $O(V_P)$, and
		$
		\ker(\pi) = (V_P \cap O(V_{A_2}))/O(V_P).
		$
	\end{remark}


	
	%
	
	Recall that $\la_1=\frac{1}{3}\al+\frac{2}{3}\b$ and $\la_2=\frac{2}{3}\al+\frac{1}{3}\b$ are the fundamental dominant weights of the $A_2$-weight lattice $\Pi=\Z\la_1\op \Z \la_2$. The irreducible $V_P$-modules are 
	$$L^{(0,\la)}\quad \mathrm{and}\quad L^{(\frac{1}{2},\la)},\quad \mathrm{where}\quad \la\in (\C\al)^\perp,$$ 
	see \eqref{5.4} and \eqref{5.5}. The bottom-degrees of irreducible $V_P$-modules are
	$$\Om_{V_P}(L^{(0,\la)})=	U^{(0,\la)}=\C e\quad \mathrm{and}\quad \Om_{V_P}(L^{(\frac{1}{2}\al,\la)})=U^{(\frac{1}{2}\al,\la)}=\C e^+\op\C e^-,$$
	with $A(V_P)$-actions given by \eqref{5.20}--\eqref{5.22}. 
	
	\begin{prop}\label{lm:inducibleV_Pmodules}
		Denote $A_{V_P}=\pi(A(V_P))\leq A(V_{A_2})$, where $\pi$ is in \eqref{eq:8.1}. Then we have the following description of $A_{V_P}$-modules:
		\begin{equation}\label{eq:inducibleVPmodules}
			\frac{\Om_{V_P}(L^{(0,\la)})}{\ker(\pi).\Om_{V_P}(L^{(0,\la)})}=\begin{cases}
				U^{(0,\la)}&\mathrm{if}\ \la=0\ \mathrm{or}\ \pm\la_2,\\
				0&\mathrm{if}\ \la\in (\C\al)^\perp\bs\{ 0,\pm\la_2\}.
			\end{cases}
		\end{equation}
		Moreover, we have 
		\begin{equation}\label{eq:inducibleVPmodules2}
			\frac{\Om_{V_P}(L^{(\frac{1}{2}\al,\la)})}{\ker(\pi).\Om_{V_P}(L^{(\frac{1}{2}\al,\la)})}=
			\begin{cases}
				U^{(\frac{1}{2}\al,\la)} &\mathrm{if}\ \la=\pm \frac{1}{2}\la_2,\\
				0&\mathrm{if}\ \la\in (\C\al)^\perp\bs\{\pm\frac{1}{2}\la_2\}.
			\end{cases}
		\end{equation}
	\end{prop}
	\begin{proof}
		Consider the $A(V_P)$-module $U^{(0,\la)}=\C e$. By \eqref{5.21} and \eqref{eq:extrarelations}, we have $$\ker(\pi). U^{(0,\la)}=\spn\{(y^3-y).e \}=\spn\{((\la|\b)^3-(\la|\b))\cdot e  \}=\begin{cases}
			0&\mathrm{if}\ \la=0\ \mathrm{or}\ \pm \la_2,\\
			U^{(0,\la)} &\mathrm{if}\  \la\in (\C\al)^\perp\bs\{ 0,\pm\la_2\}.
		\end{cases}$$
		This proves \eqref{eq:inducibleVPmodules}. 
		
		Consider the irreducible $A(V_P)$-module $U^{(\frac{1}{2}\al,\la)}=\C e^+\op\C e^-$. The submodule $\ker(\pi).U^{(\frac{1}{2}\al,\la)}$ is either $0$ or $U^{(\frac{1}{2}\al,\la)}$.
		By  \eqref{5.22}, we have $y.e^{\pm}=((\la|\b)\mp \frac{1}{2})e^\pm$. Then 
		\begin{equation}\label{eq:8.5}
			(y^3-y).e^\pm=\left(\left((\la|\b)\mp (1/2)\right)^3-\left((\la|\b)\mp (1/2)\right)\right)\cdot e^\pm \in \ker(\pi). U^{(\frac{1}{2}\al,\la)}.
		\end{equation}
		If $(\la|\b)\neq \pm \frac{1}{2}$, or equivalently, if $\la\neq \pm \frac{1}{2}\la_2$, then the coefficients of $e^+$ and $e^-$ in \eqref{eq:8.5} are not zero at the same time, and so $\ker(\pi).U^{(\frac{1}{2}\al,\la)}=U^{(\frac{1}{2}\al,\la)}$. 
		
		On the other hand, if $\la=\frac{1}{2}\la_2$, then $y.e^+=0$ and $y.e^-=e^-$. By \eqref{5.22}, 
		\begin{align*}
			(x_{\pm \al}(x+y)^2\pm x_{\pm \al}(x+y)).e^+&=x_{\pm \al}.e^+\pm x_{\pm \al}.e^+=0\quad (\mathrm{as}\ x_\al.e^+=0),\\
			(x_{\pm \al}(x+y)^2\pm x_{\pm \al}(x+y)).e^-&=x_{\pm \al}.(-1+1)^2e^-\pm x_{\pm \al}.(-1+1)e^-=0,\\
			(x_{\pm \al} y^2\mp x_{\pm \al} y).e^+&=0,\\
			(x_{\pm \al} y^2\mp x_{\pm \al} y).e^-&=\begin{cases}
				x_\al.e^--x_\al.e^-=0,\quad \mathrm{or}\\
				x_{-\al}.e^-+x_{-\al}.e^-=0\quad (\mathrm{as}\ x_{-\al}.e^-=0),
			\end{cases}\\
			((x+y)^3-(x+y)).e^+&=x^3.e^+-x.e^+=e^+-e^+=0,\\
			((x+y)^3-(x+y)).e^-&=(-1+1)^3e^--(-1+1)e^-=0.
		\end{align*}
		Thus, $\ker(\pi).U^{(\frac{1}{2}\al, \frac{1}{2}\la_2)}=0$ in view of \eqref{eq:extrarelations}. 
		Finally, if $\la=-\frac{1}{2}\la_2$, then $y.e^+=-e^+$ and $y.e^-=0$. Similarly, we can show that $\ker(\pi).U^{(\frac{1}{2}\al, -\frac{1}{2}\la_2)}=0$. This proves \eqref{eq:inducibleVPmodules2}. 
	\end{proof}
	
	\subsection{Irreducible modules under the finite induction for $V_P\hr V_{A_2}$}

		Recall that the lattice VOA $V_{A_2}\cong L_{\widehat{\sl_3}}(1,0)$ is strongly rational \cite{D,DLM1}. It has three irreducible modules with bottom degree: 
		\begin{align*}
			\Om_{V_{A_2}}(V_{A_2})&=\C \vac\cong L(0),\\
			\Om_{V_{A_2}}(	V_{A_2+\la_1})&=\C e^{\la_1}+\C e^{\la_1-\al}+\C e^{\la_1-\al-\b}\cong L(\la_1),\\ 
			\Om_{V_{A_2}}(V_{A_2+\la_2})&=\C e^{\la_2}+\C e^{\la_2-\b}+\C e^{\la_2-\al-\b}\cong L(\la_2),
		\end{align*}
		as modules over $\sl_3$. 
		\begin{thm}\label{thm:inducedmodulesforV_P}
			The finite induction of irreducible $V_P$-modules under the VOA embedding $V_P \hookrightarrow V_{A_2}$ satisfies
			\[
			\Ind^{V_{A_2}}_{V_P} L^{(\epsilon,\lambda)} = 0 \quad 
			\text{if } (\epsilon,\lambda)\notin \{(0,0),\,(0,\pm \lambda_2),\, (\tfrac{1}{2}\alpha,\pm \tfrac{1}{2}\lambda_2)\}.
			\]
	For the remaining pairs of $(\epsilon,\la)$, we have 
			\begin{align*}
				&\Ind^{V_{A_2}}_{V_P}L^{(0,0)} \cong V_{A_2}, 
				&& \Ind^{V_{A_2}}_{V_P}L^{(0, \lambda_2)} \cong V_{A_2+\lambda_2}, 
				&& \Ind^{V_{A_2}}_{V_P}L^{(0,- \lambda_2)} \cong 0,\\
				&\Ind^{V_{A_2}}_{V_P}L^{(\frac{1}{2}\alpha, \frac{1}{2}\lambda_2)} \cong  V_{A_2+\lambda_1}, 
				&& \Ind^{V_{A_2}}_{V_P}L^{(\frac{1}{2}\alpha, -\frac{1}{2}\lambda_2)} \cong  0,
			\end{align*}
			where $0$ denotes the zero module. 
		\end{thm}
		
		\begin{proof}
			By the definition of finite induction~\eqref{eq:ind}, it suffices to determine the left $A(V_{A_2})$-module 
			\begin{equation}
				\Omega = A(V_{A_2}) \otimes_{A_{V_P}} \frac{\Omega_{V_P}(L^{(\epsilon,\lambda)})}
				{\ker(\pi).\Omega_{V_P}(L^{(\epsilon,\lambda)})}.
			\end{equation}
			By Proposition~\ref{lm:inducibleV_Pmodules}, $\Omega = 0$ if 
			$(\epsilon,\lambda)\notin \{(0,0),\,(0,\pm \lambda_2),\,(\frac{1}{2}\alpha,\pm \frac{1}{2}\lambda_2)\}$.
			We discuss the remaining cases one by one. Note that $	{\ker(\pi).\Omega_{V_P}(L^{(\epsilon,\lambda)})}=0$ for the remaining cases. 
			

			\noindent\textbf{Case I.} $(\epsilon, \lambda) = (0,0)$. 
			By \eqref{5.20} and \eqref{5.21}, we have $U^{(0,0)} = \C \vac$, and
			\[
			x_{\pm \alpha}.\vac = x_{\beta}.\vac = x_{\alpha+\beta}.\vac = x.\vac = y.\vac = 0,
			\]
			where we denote $\vac \otimes e^0$ simply by $\vac$. 
			On the other hand, from \eqref{rel1}--\eqref{rel4}, the following relations hold in 
			$A(V_{A_2}) \otimes_{A_{V_P}} U^{(0,0)}$:
			\[
			x_{-\beta} \otimes \vac 
			= x_{-\beta}y \otimes \vac 
			= x_{-\beta} \otimes y.\vac = 0, 
			\quad 
			x_{-\alpha-\beta} \otimes \vac 
			= x_{-\alpha-\beta}(x+y) \otimes \vac 
			= x_{-\alpha-\beta} \otimes (x+y).\vac = 0.
			\]
			Thus, all the generators of $A(V_{A_2})$ act trivially on $1 \otimes \vac 
			\in A(V_{A_2}) \otimes_{A_{V_P}} U^{(0,0)}$.
			
			We now claim that the $A(V_{A_2})$-module $A(V_{A_2}) \otimes_{A_{V_P}} U^{(0,0)}$ is nonzero. 
			Indeed, by the Hom-tensor duality we have
			\[
			\big(A(V_{A_2}) \otimes_{A_{V_P}} U^{(0,0)}\big)^{\!*}
			\cong 
			\Hom_{A_{V_P}}\!\big(\C \vac, (A(V_{A_2}))^{\!*}\big).
			\]
			Note that $A(V_{A_2}) = \C 1 \oplus A(V_{A_2})_+$ as vector spaces, 
			where $A(V_{A_2})_+$ consists of elements in 
			$A(V_{A_2}) \cong U(\sl_3(\C))/\langle x_{\alpha+\beta}^2\rangle$ of length greater than one.

			Define a linear map 
			$\phi : \C \vac \to (A(V_{A_2}))^{\!*}$ by
			\begin{equation}\label{eq:phiforU00}
				\langle \phi(\vac), 1 \rangle := 1,
				\qquad
				\langle \phi(\vac), a \rangle := 0, 
				\quad \forall\, a \in A(V_{A_2})_+,
			\end{equation}
			where 
			$\langle \cdot, \cdot \rangle : (A(V_{A_2}))^{\!*} \times A(V_{A_2}) \to \C$ 
			is the natural pairing.
			
			We show that $\phi$ is an $A_{V_P}$-homomorphism. 
			Let 
			$S(A_{V_P}) := \{ x_{\pm\alpha}, x_{\beta}, x_{\alpha+\beta}, x, y \}$ 
			be the set of generators of $A_{V_P}$. 
			Consider the left ideal $A(V_{A_2}) \cdot S(A_{V_P}) \subset A(V_{A_2})$. 
			Using the relations \eqref{4.21}--\eqref{4.26}, together with 
			\eqref{eq:extrarelations} and \eqref{rel1}--\eqref{rel4}, 
			it is easy to verify that $1 \notin A(V_{A_2}) \cdot S(A_{V_P})$. 
			Hence, for any $z \in S(A_{V_P})$ and $a \in A(V_{A_2})$,
			\begin{equation}\label{eq:propertyforphi}
				\langle \phi(z.\vac), a \rangle 
				= 0 
				= \langle \phi(\vac), a \cdot z \rangle
				= \langle z.\phi(\vac), a \rangle,
			\end{equation}
			in view of \eqref{eq:phiforU00}. 
			When $z = 1$, the same holds since $a \cdot 1 = a$ for all $a \in A(V_{A_2})$.
			Thus, $\phi$ defined in \eqref{eq:phiforU00} is a nonzero element of 
			$\Hom_{A_{V_P}}\!\big(\C \vac, (A(V_{A_2}))^{\!*}\big)$.
			Therefore,
			\[
			A(V_{A_2}) \otimes_{A_{V_P}} U^{(0,0)} \neq 0.
			\]
			In particular,
			$
			A(V_{A_2}) \otimes_{A_{V_P}} U^{(0,0)} 
			\cong L(0)$ and
			$
			\Ind^{V_{A_2}}_{V_P}L^{(0,0)} 
			\cong V_{A_2}.
			$

			\noindent	\textbf{Case II.} $(\epsilon, \lambda) = (0, \lambda_2)$. 
			By \eqref{5.20} and \eqref{5.21} we have $U^{(0, \lambda_2)} = \C e$, with 
			$x.e = x_{\beta}.e = x_{\alpha+\beta}.e = x_{\pm \alpha}.e = 0$ and $y.e = e$. 
			Again, we first show that $A(V_{A_2}) \otimes_{A_{V_P}} U^{(0, \lambda_2)}$ is nonzero. 
			Similar to Case~I, using the Hom–tensor duality
			\[
			\left(A(V_{A_2}) \otimes_{A_{V_P}} U^{(0, \lambda_2)}\right)^\ast 
			\cong \Hom_{A_{V_P}}(\C e, (A(V_{A_2}))^\ast),
			\]
			we just need to find a nonzero element 
			$\varphi \in \Hom_{A_{V_P}}(\C e, (A(V_{A_2}))^\ast)$. 
			Using the relations \eqref{4.21}--\eqref{4.26}, \eqref{eq:extrarelations}, and \eqref{rel1}--\eqref{rel4}, 
			we can write 
			\[
			A(V_{A_2}) = \mathrm{span}\{1, y, y^2\} \oplus C,
			\]
			where
			\[
			C = \mathrm{span}\{x, x^2, xy, x^2y, x_{\pm \alpha}, x_{\pm \beta}, x_{\pm (\alpha+\beta)}, 
			x_{\pm \alpha}y, x_{\pm \beta}x, x_{\pm (\alpha+\beta)}y\}.
			\]
			Then we define $\varphi : \C e \to (A(V_{A_2}))^\ast$ by letting
			\begin{equation}\label{eq:varphiforU0la2}
				\langle \varphi(e), 1 \rangle = \langle \varphi(e), y \rangle = \langle \varphi(e), y^2 \rangle := 1, 
				\quad 
				\langle \varphi(e), a \rangle := 0, \ \forall a \in C.
			\end{equation}
			
			Let $S := \{x_{\pm \alpha}, x_{\beta}, x_{\alpha+\beta}, x\}$. 
			Since $y.e = e$, to show $\varphi$ is an $A_{V_P}$-homomorphism, we need to verify that
			\begin{equation}\label{eq:needtoshowforvarphi}
				\langle \varphi(e), a \rangle = \langle \varphi(e), a \cdot y \rangle, 
				\quad \text{and} \quad 
				\langle \varphi(e), a \cdot z \rangle = 0, 
				\quad \forall a \in A(V_{A_2}), \ z \in S.
			\end{equation}
			
			By \eqref{rel1}--\eqref{rel4}, the only terms in the left ideal 
			$A(V_{A_2}) \cdot S$ that have nonzero components in the subspace 
			$\mathrm{span}\{1, y, y^2\}$ are spanned by 
			$\frac{1}{2}y^2 - \frac{1}{2}y$ and 
			$\frac{1}{2}(x+y)^2 - \frac{1}{2}(x+y)$. 
			By \eqref{eq:varphiforU0la2}, we have 
			$\langle \varphi(e), y^2 - y \rangle = 0$ and 
			$\langle \varphi(e), (x+y)^2 - (x+y) \rangle 
			= \langle \varphi(e), y^2 - y \rangle = 0$. 
			Thus $\langle \varphi(e), A(V_{A_2}) \cdot S \rangle = 0$. 
			Moreover, it is easy to show that $C \cdot y \subseteq C$, 
			hence $\langle \varphi(e), a \rangle = 0 = \langle \varphi(e), a \cdot y \rangle$ 
			for any $a \in C$, in view of \eqref{eq:varphiforU0la2}. 
			Finally, for $a = \lambda + \mu y + \gamma y^2 \in \mathrm{span}\{1, y, y^2\}$, 
			we have 
			\[
			\langle \varphi(e), a \rangle 
			= \lambda + \mu + \gamma 
			= \langle \varphi(e), \lambda y + \mu y^2 + \gamma \rangle 
			= \langle \varphi(e), a \cdot y \rangle,
			\]
			in view of \eqref{eq:needtoshowforvarphi}. 
			This shows that $\varphi$ is a nonzero element in 
			$$\Hom_{A_{V_P}}(\C e, (A(V_{A_2}))^\ast) 
			\cong \left(A(V_{A_2}) \otimes_{A_{V_P}} U^{(0, \lambda_2)}\right)^\ast.$$ 
			
			Now $1 \otimes e \in A(V_{A_2}) \otimes_{A_{V_P}} U^{(0, \lambda_2)}$ 
			is a highest-weight vector for $\mathfrak{sl}_3$ of highest weight $\lambda_2$, 
			since 
			$x.(1 \otimes e) = 1 \otimes x.e = 0 = (\alpha | \lambda_2) \cdot (1 \otimes e)$ 
			and 
			$y.(1 \otimes e) = 1 \otimes y.e = (\beta | \lambda_2) \cdot (1 \otimes e)$. 
			Then 
			$A(V_{A_2}) \otimes_{A_{V_P}} U^{(0, \lambda_2)} 
			= U(\mathfrak{sl}_3)/\langle x^2_{\alpha+\beta} \rangle \cdot (1 \otimes e)$ 
			is a finite-dimensional highest-weight $\mathfrak{sl}_3$-module 
			of highest weight $\lambda_2$. 
			Hence it is isomorphic to $L(\lambda_2)$~\cite{V68, Hum2}. 
			This shows that 
			\[
			\mathrm{Ind}^{V_{A_2}}_{V_P} L^{(0, \lambda_2)} 
			= \Phi^\L_{V_{A_2}}(L(\lambda_2)) \cong V_{A_2 + \lambda_2}.
			\]

			
			\noindent	\textbf{Case III.} $(\epsilon, \lambda) = (0, -\lambda_2)$. 
			We claim that $A(V_{A_2}) \otimes_{A_{V_P}} U^{(0, -\lambda_2)} = 0$. 
			Indeed, since $U^{(0, -\lambda_2)} = \C e$, with 
			$x.e = x_{\beta}.e = x_{\alpha+\beta}.e = x_{\pm \alpha}.e = 0$ 
			and $y.e = -e$, we have 
			\[
			A(V_{A_2}) \otimes_{A_{V_P}} U^{(0, -\lambda_2)} 
			= \mathrm{span}\{1 \otimes e,\, x_{-\beta} \otimes e,\, x_{-\alpha-\beta} \otimes e\}.
			\]
			
			By \eqref{rel1}--\eqref{rel4}, we have $x_{-\beta}x_{\beta} = \frac{1}{2}y^2 - \frac{1}{2}y$. 
			Note that 
			\[
			\Bigl(\tfrac{1}{2}y^2 - \tfrac{1}{2}y\Bigr).e 
			= \tfrac{1}{2}(-1)^2 e + \tfrac{1}{2}e = e,
			\]
			and so
			\[
			1 \otimes e 
			= 1 \otimes \Bigl(\tfrac{1}{2}y^2 - \tfrac{1}{2}y\Bigr).e
			= \Bigl(\tfrac{1}{2}y^2 - \tfrac{1}{2}y\Bigr).(1 \otimes e)
			= x_{-\beta} \otimes x_{\beta}.e
			= x_{-\beta} \otimes 0
			= 0.
			\]
			Moreover, $x_{-\beta} \otimes e = x_{-\beta}.(1 \otimes e) = 0$ and 
			$x_{-\alpha-\beta} \otimes e = x_{-\alpha-\beta}.(1 \otimes e) = 0$. 
			Thus $A(V_{A_2}) \otimes_{A_{V_P}} U^{(0, -\lambda_2)} = 0$, and hence
			$
			\mathrm{Ind}^{V_{A_2}}_{V_P} L^{(0, -\lambda_2)}  \cong 0.
			$
			
			\noindent\textbf{Case IV.} $(\epsilon, \lambda) = \left(\tfrac{1}{2}\alpha, \tfrac{1}{2}\lambda_2\right)$. 
			Similar to Cases~I and~II, we first show that the 
			$A(V_{A_2})$-module $A(V_{A_2}) \otimes_{A_{V_P}} U^{(\frac{1}{2}\alpha, \frac{1}{2}\lambda_2)}$ 
			is nonzero. By \eqref{5.20} and \eqref{5.22}, we have 
			$U^{(\frac{1}{2}\alpha, \frac{1}{2}\lambda_2)} = \C e^+ \oplus \C e^-$, 
			on which the generators of $A_{V_P}$ act as
			\begin{equation}\label{eq:generatorsact}
				\begin{aligned}
					&x_{\alpha}.e^+ = 0, && x_{-\alpha}.e^+ = e^-, && x.e^{\pm} = \pm e^{\pm}, &&
					x_{\beta}.e^{\pm} = x_{\alpha+\beta}.e^{\pm} = 0,\\
					&x_{\alpha}.e^- = e^+, && x_{-\alpha}.e^- = 0, &&
					y.e^+ = 0, && y.e^- = e^-.
				\end{aligned}
			\end{equation}

			Write $A(V_{A_2}) = \mathrm{span}\{1, x, x^2\} \oplus D$, where 
			\[
			D = \mathrm{span}\{y, y^2, xy, x^2y, x_{\pm\alpha}, x_{\pm\beta}, x_{\pm(\alpha+\beta)}, 
			x_{\pm\alpha}y, x_{\pm\beta}x, x_{\pm(\alpha+\beta)}y\}.
			\]
			Define a linear map 
			$\psi : \C e^+ \oplus \C e^- \to (A(V_{A_2}))^{\ast}$ by letting
			\begin{equation}\label{eq:defforpsi}
				\begin{aligned}
					&\langle \psi(e^+), 1 \rangle = \langle \psi(e^+), x \rangle 
					= \langle \psi(e^+), x^2 \rangle := 1, 
					\quad \langle \psi(e^+), D \rangle := 0,\\
					&\langle \psi(e^-), a \rangle := \langle \psi(e^+), a \cdot x_{-\alpha} \rangle, 
					\quad \forall\, a \in A(V_{A_2}).
				\end{aligned}
			\end{equation}
			We want to show that 
			$\psi \in \Hom_{A_{V_P}}(U^{(\frac{1}{2}\alpha, \frac{1}{2}\lambda_2)}, (A(V_{A_2}))^{\ast})$, 
			i.e.,
			\begin{equation}\label{eq:needtoshowfinal}
				\langle \psi(z.e^{\pm}), a \rangle 
				= \langle \psi(e^{\pm}), a \cdot z \rangle, 
				\quad \forall a \in A(V_{A_2}),\ z \in \{x_{\pm \alpha}, x_{\beta}, x_{\alpha+\beta}, x, y\}.
			\end{equation}

			Indeed, by \eqref{eq:generatorsact}, 
			$S = \{x_{\alpha}, x_{\beta}, x_{\alpha+\beta}, y, (x - 1)\}$ annihilates $e^+$. 
			It is easy to show that the left ideal $A(V_{A_2}) \cdot S$ 
			is contained in $D + \mathrm{span}\{(x - 1), (x^2 - x)\}$ 
			by relations \eqref{4.21}--\eqref{4.26}. 
			Since $\langle \psi(e^+), D \rangle 
			= \langle \psi(e^+), x - 1 \rangle 
			= \langle \psi(e^+), x^2 - x \rangle = 0$ by \eqref{eq:defforpsi}, 
			we have
			\[
			\langle \psi(u.e^+), A(V_{A_2}) \rangle = 0 
			= \langle \psi(e^+), A(V_{A_2}) \cdot u \rangle, 
			\quad \forall u \in S.
			\]
			Thus \eqref{eq:needtoshowfinal} holds for all $z \in \{x_{\alpha}, x_{\beta}, x_{\alpha+\beta}, y, x\}$. 
			For $z = x_{-\alpha}$, by \eqref{eq:defforpsi} and $x_{-\alpha}.e^+ = e^-$, 
			we have 
			\[
			\langle \psi(x_{-\alpha}.e^+), a \rangle 
			= \langle \psi(e^-), a \rangle 
			= \langle \psi(e^+), a \cdot x_{-\alpha} \rangle.
			\]
			Hence \eqref{eq:needtoshowfinal} holds for all $z \in A_{V_P}$ acting on $e^+$.

			On the other hand, 
			$T = \{x_{-\alpha}, x_{\beta}, x_{\alpha+\beta}, (x + 1), (y - 1)\}$ annihilates $e^-$ 
			by \eqref{eq:generatorsact}. 
			Noting $x_{-\alpha}x_{-\alpha} = (x + 1)x_{-\alpha} = x_{\beta}x_{-\alpha} = 0$, 
			we get from \eqref{eq:defforpsi}
			\begin{align*}
				\langle \psi(x_{-\alpha}.e^-), a \rangle &= 0 = \langle \psi(e^-), a \cdot x_{-\alpha} \rangle,\\
				\langle \psi((x + 1).e^-), a \rangle &= 0 = \langle \psi(e^-), a \cdot (x + 1) \rangle,\\
				\langle \psi(x_{\beta}.e^-), a \rangle &= 0 = \langle \psi(e^-), a \cdot x_{\beta} \rangle.
			\end{align*}
			Moreover, by \eqref{rel1}--\eqref{rel4}, 
			$(y - 1)x_{-\alpha} = x_{-\alpha}y$ and 
			$x_{\alpha+\beta}x_{-\alpha} = -x_{\beta}x$. 
			Since $\langle \psi(e^+), A(V_{A_2}) \cdot y \rangle = 0$, 
			\begin{align*}
				\langle \psi((y - 1).e^-), a \rangle 
				= 0 
				= \langle \psi(e^-), a \cdot (y - 1) \rangle.
			\end{align*}
			Also, since $\psi$ is invariant for $e^+$, we have 
			$\langle \psi(e^+), A(V_{A_2}) \cdot x_{\beta}x \rangle = 0$, 
			and hence 
			\[
			\langle \psi(x_{\alpha+\beta}.e^-), a \rangle 
			= 0 
			= \langle \psi(e^-), a \cdot x_{\alpha+\beta} \rangle.
			\]
			Finally, for $z = x_{\alpha}$, since 
			$\langle \psi(e^+), a \cdot (x^2 - 1) \rangle 
			= \langle \psi(e^+), a \cdot (x - 1) \rangle = 0$, 
			we have
			\[
			\langle \psi(e^+), a \rangle 
			= \langle \psi(e^-), a \cdot x_{\alpha} \rangle,
			\]
			while $\langle \psi(e^+), a \rangle 
			= \langle \psi(x_{\alpha}.e^-), a \rangle$ 
			since $x_{\alpha}.e^- = e^+$. 
			This completes the proof of \eqref{eq:needtoshowfinal}. 
			Hence $\psi$ is a nonzero element of
			\[
			\Hom_{A_{V_P}}\bigl(U^{(\frac{1}{2}\alpha, \frac{1}{2}\lambda_2)}, 
			(A(V_{A_2}))^{\ast}\bigr)
			\cong 
			\left(A(V_{A_2}) \otimes_{A_{V_P}} 
			U^{(\frac{1}{2}\alpha, \frac{1}{2}\lambda_2)}\right)^{\ast}.
			\]

			From \eqref{eq:generatorsact}, 
			$1 \otimes e^+$ is a highest-weight vector in 
			$A(V_{A_2}) \otimes_{A_{V_P}} U^{(\frac{1}{2}\alpha, \frac{1}{2}\lambda_2)}$ 
			of highest weight $\lambda_1$, 
			since $y.e^+ = 0 = (\lambda_1 | \beta)e^+$ and 
			$x.e^+ = e^+ = (\lambda_1 | \alpha)e^+$. 
			Moreover, $1 \otimes e^- \in U(\mathfrak{sl}_3) \cdot (1 \otimes e^+)$, 
			and hence 
			\[
			A(V_{A_2}) \otimes_{A_{V_P}} U^{(\frac{1}{2}\alpha, \frac{1}{2}\lambda_2)} 
			= U(\mathfrak{sl}_3) \cdot (1 \otimes e^{\pm}) 
			\cong L(\lambda_1),
			\]
			as $\mathfrak{sl}_3$-modules~\cite{V68}. 
			Therefore,
			\[
			\mathrm{Ind}^{V_{A_2}}_{V_P} L^{(\frac{1}{2}\alpha, \frac{1}{2}\lambda_2)} 
			= \Phi^\L_{V_{A_2}}(L(\lambda_1)) 
			\cong V_{A_2 + \lambda_1}.
			\]

			\noindent	\textbf{Case V.} $(\epsilon, \lambda) = \left(\tfrac{1}{2}\alpha, -\tfrac{1}{2}\lambda_2\right)$. 
			Similar to Case~III, we claim that 
			$A(V_{A_2}) \otimes_{A_{V_P}} U^{(\frac{1}{2}\alpha, -\frac{1}{2}\lambda_2)} = 0$. 
			Indeed, $U^{(\frac{1}{2}\alpha, -\frac{1}{2}\lambda_2)} = \C e^+ \oplus \C e^-$, 
			with all relations in \eqref{eq:generatorsact} holding except that 
			$y.e^+ = -e^+$ and $y.e^- = 0$. 
			Then by \eqref{rel1}--\eqref{rel4},
			\[
			1 \otimes e^+ 
			= 1 \otimes \Bigl(\tfrac{1}{2}y^2 - \tfrac{1}{2}y\Bigr).e^+ 
			= \Bigl(\tfrac{1}{2}y^2 - \tfrac{1}{2}y\Bigr).(1 \otimes e^+) 
			= x_{-\beta} \otimes x_{\beta}.e^+ = 0,
			\]
			and 
			$1 \otimes e^- = 1 \otimes x_{-\alpha}.e^+ 
			= x_{-\alpha}.(1 \otimes e^+) = 0$. 
			Hence
			\[
			A(V_{A_2}) \otimes_{A_{V_P}} 
			U^{(\frac{1}{2}\alpha, -\frac{1}{2}\lambda_2)} 
			= U(\mathfrak{sl}_3) \cdot (1 \otimes e^+) 
			+ U(\mathfrak{sl}_3) \cdot (1 \otimes e^-) 
			= 0,
			\]
			and therefore
			$
			\mathrm{Ind}^{V_{A_2}}_{V_P} L^{(\frac{1}{2}\alpha, -\frac{1}{2}\lambda_2)}  \cong 0.
			$
		\end{proof}
		


	\end{document}